\documentclass[12pt]{amsart}
\usepackage{amssymb}

\numberwithin{equation}{section}  
\theoremstyle{plain}
\newtheorem{thm}{Theorem}[section]

\newtheorem{cor}[thm]{Corollary}
\newtheorem{lem}{Lemma}[section]

\theoremstyle{definition}

\newtheorem{example}{Example}[section]

\theoremstyle{remark}

\def\R{{\mathbb R}}
\def\C{{\mathbb C}}
\def\Rn{{{\mathbb R}^n}}

\def\Rnx{{\R^n_x}}
\def\N{{\mathbb N}}

\def\Sph{{\mathbb S}}
\def\FT{{\mathcal F}}
\def\L2tx{{L^2(\R_t\times\R^n_x)}}
\def\Lx{{L^2(\R^n_x)}}
\def\p#1{{\left({#1}\right)}}
\def\b#1{{\left\{{#1}\right\}}}
\def\br#1{{\left[{#1}\right]}}
\def\n#1{{\left\|{#1}\right\|}}
\def\abs#1{{\left|{#1}\right|}}
\def\jp#1{{\left\langle{#1}\right\rangle}}
\def\Im{\operatorname{Im}}

\def\supp{\operatorname{supp}}
\def\tr{\operatorname{Tr}}
\def\rank{\operatorname{rank}}
\def\va{\varphi}
\def\ka{\kappa}

\title{Smoothing properties of evolution
equations via canonical transforms and comparison}

\author[]{Michael Ruzhansky and Mitsuru Sugimoto}
\address{
  Michael Ruzhansky:
  \endgraf
  Department of Mathematics
  \endgraf
  Imperial College London
  \endgraf
  180 Queen's Gate, London SW7 2AZ, UK
  \endgraf
  {\it E-mail address} {\rm m.ruzhansky@imperial.ac.uk}
  \endgraf
  \medskip
  Mitsuru Sugimoto:
  \endgraf
  Department of Mathematics, Graduate School of Science
  \endgraf
  Osaka University
  \endgraf
  Machikaneyama-cho 1-16, Toyonaka, Osaka 560-0043, Japan
  \endgraf
  {\it E-mail address} {\rm sugimoto@math.wani.osaka-u.ac.jp}
  }
\thanks{This work was completed with the aid of
\lq\lq UK-Japan Joint Project Grant" by
\lq\lq The Royal Society (UK)" and
\lq\lq Japan Society for the Promotion of Science". The first
author was also supported by the Leverhulme Research
Fellowship.
}

\date{\today}

\begin{document}

\begin{abstract}
The paper describes a new approach to global smoothing problems
for dispersive and non-dispersive evolution equations 
based on the global canonical transforms
and the underlying global microlocal analysis.
For this purpose, the Egorov--type theorem is established
with canonical
transformations in the form of a class of 
Fourier integral operators, and their weighted 
$L^2$--boundedness properties are derived.
This allows us to globally reduce general dispersive 
equations to normal forms in one or two dimensions.
Then, several new comparison techniques for evolution
equations are introduced. In particular, they allow us
to relate different smoothing estimates by comparing
certain expressions involving their symbols.
As a result, it is shown that the majority of smoothing
estimates for different equations are equivalent
to each other.
Moreover, new estimates as well as 
several refinements of known results are obtained.
The proofs are considerably simplified.
A comprehensive 
analysis is presented of smoothing estimates for 
homogeneous and inhomogeneous, dispersive 
and also non-dispersive equations with constant coefficients.
Results are presented also for equations with time dependent
coefficients.
Applications are given to the detailed description of smoothing
properties of the Schr\"odinger,
relativistic Schr\"odinger, wave, Klein-Gordon, and other equations.
Critical cases of some estimates and their relation to the
trace estimates are discussed.
\end{abstract}

\maketitle

\tableofcontents

\section{Introduction}
In the last two decades, since the independent pioneering works
by Ben-Artzi and Devinatz
\cite{BD2}, Constantin and Saut \cite{CS},
Sj\"olin \cite{Sj} and Vega \cite{V}, 
the local, and then global smoothing effects of 
Schor\"odinger equations,
or more generally, those of dispersive equations has been 
intensively investigated.
Similar smoothing effects have been observed for different
equations of great importance in mathematical physics
(for example, smoothing for generalized Korteweg-de Vries equations
was already studied by Kato \cite{Ka2}, several other equations
were studied in a series of papers by
Kenig, Ponce and Vega \cite{KPV1}--\cite{KPV5}), etc).
Over the years, several techniques to understand these
smoothing properties through 
the Fourier analysis, functional analysis,
spectral theory and harmonic analysis
have been developed. The analysis of such smoothing estimates
is particularly important in applications to nonlinear
evolution equations, especially to those
with derivatives in the potential or in the nonlinearity.
Over the last three decades two major approaches,
Strichartz and smoothing estimates, proved to be
two extremely efficient tools for dealing with nonlinear
equations.
The smoothing effect
is crucial in allowing to recover the loss of
derivatives in the equation making
these estimates a very good substitute for the 
Strichartz estimates that are
normally used for semilinear equations.
\par
The objective of this paper is to provide a new approach 
leading to
a comprehensive understanding of the effect of 
global smoothing, 
together with new results, through two novel ideas.
It will allow us not only to recover existing and to prove new
estimates, but to effectively show that the smoothing phenomenon
for equations describing often completely different physical
processes (like wave, Klein-Gordon,
Schr\"odinger, relativistic Schr\"odinger,
KdV, Benjamin-Ono, Davey-Stewartson,
Shrira, Zakharov-Schulman, and many
other equations) is of essentially the same nature.
For this, we will provide a way to show the equivalence of
smoothing properties for very different equations by
introducing two new ideas for the subject. 

The first is the idea of
{\it canonical transformations}.
Let us mention immediately that
will not concentrate specifically on the local smoothing
since it is contained in its global version.
Although canonical transforms are
well known in microlocal analysis of local problems, 
we will apply them in a
global setting here to globally reduce problems to 
normal forms in lower dimensions. As it turns out, 
it is then possible to carry out 
the pointwise analysis in these model
problems. 
However, there are some essential differences 
with the microlocal case.
On one hand, we will still be able to reduce elliptic operators
to one dimensional models. On the other hand, in the case
of dispersive operators (or operators of real principal 
type) the global
reduction will be made to models in two dimensions,
in difference with the well-known microlocal constructions 
of Duistermaat and H\"ormander \cite{DH}. 
\par
Another important idea presented and explored in 
this paper is 
an introduction of a certain 
{\it comparison principle} for evolution equations
which will allow us to derive new estimates for dispersive 
(and non-dispersive) equations
from known ones, as well as compare estimates for different
equations.
The idea here is that we can compare certain
expressions involving symbols and weights for different
estimates and conclude that one estimate implies the other
if an inequality between these expressions holds. In particular,
it will imply that smoothing estimates are equivalent if
certain expression involving symbols are equivalent. 
This will, for example, show that a variety of
global smoothing estimates
for Schr\"odinger equations are simply
equivalent to the corresponding
estimates for the relativistic Schr\"odinger, wave, 
Klein--Gordon, linearised KdV, Benjamin--Ono, and
other equations. In addition, it will show that the local smoothing
effect for Schr\"odinger equations that was established by
Sj\"olin \cite{Sj} and Vega \cite{V} is equivalent to the
energy conservation of a travelling wave in one dimension.
The gain of $1/2$-derivative corresponds to the Jacobian of the
frequency transformation between Schr\"odinger and a one-dimensional
wave in the radial direction, and the $1/2$-derivative smoothing
for Schr\"odinger is the energy estimate for this wave, 
which in
turn is just the translation invariance property
of the Lebesgue measure on the real line
(see \eqref{eq:trivial} and the discussion around it).
The local gain of one derivative for Korteweg-de Vries
equation was also observed by Kato \cite{Ka2}, whose proof
used the algebraic properties of the symbol and the fact 
that the situation is one-dimensional. Again, by the comparison
principle we will immediately recover this result (and its
global version) from the $1/2$-smoothing for Schr\"odinger,
or from the energy conservation for the wave equation.

The methods and ideas described in the paper seem to 
apply very efficiently to the area of smoothing estimates.
These estimates  are usually immediately applied to
the well-posedness of nonlinear equations and many methods are
developed for this purpose.  From this point of view, results
presented here have immediate consequences for the well-posedness
problems for rather general classes of nonlinear evolution
equations. 
Moreover, it can be expected that they can be also applied to
a variety of other problems where weighted estimates and phase
space analysis are of importance.

Let us mention that there has already been a lot of
literature on the subject of global smoothing estimates
from different points of view.
See, Ben-Artzi and Devinatz \cite{BD1, BD2},
Ben-Artzi and Klainerman \cite{BK}, Chihara \cite{Ch}, 
Hoshiro \cite{Ho2},
Kato and Yajima \cite{KY},
Kenig, Ponce and Vega \cite{KPV1, KPV2, KPV3, KPV4}, 
Linares and Ponce \cite{LP}, Simon \cite{Si}, 
Sugimoto \cite{Su1, Su2},
Walther \cite{Wa1, Wa2}, and many others. The
two most commonly used techniques are based on the
limiting absorption principle from the spectral theory
or on the restriction theorems in harmonic analysis.
They are both better adapted to study elliptic equations
and often do not work in the dispersive non-elliptic cases
(let alone non-dispersive equations). 
Indeed, there may be lack of information about the
spectrum, or level surfaces of a non-elliptic 
symbol may be non-compact, may
have singularities or vanishing Gaussian curvature.
The described method of
canonical transforms works equally well in all dispersive
situations (elliptic and non-elliptic).
As such, methods developed here may be used to producing
their counterpart restriction theorems for corresponding
level surfaces.
Moreover, since one can make full use of
the microlocal analysis, most results can be extended to
different types of symbols, e.g. to
quasi-homogeneous symbols, just by microlocalising in
appropriate directions and applying results of this paper.
Because of this we will not give a full treatment of
quasi-homogeneous symbols here, mostly restricting ourselves
to homogeneous symbols with lower order terms for the
clarity of the exposition. 

Moreover, we will obtain the 
corresponding results on the global smoothing for
solutions to inhomogeneous problems. There are considerably
less results on this topic available in the literature.
Mostly the Schr\"odinger equation was treated
(e.g. Linares and Ponce \cite{LP}, Kenig, Ponce and
Vega \cite{KPV5}), or the one dimensional case
(Kenig, Ponce and Vega \cite{KPV3, KPV4} or Laurey \cite{La}).
Some more general results on the local
smoothing for dispersive operators were obtained by
Chichara \cite{Ch} and Hoshiro \cite{Ho2}, and for 
dispersive differential operators by Koch and Saut \cite{KoSa}.
In this paper we will extend these results in two directions:
we will establish the global smoothing for rather general
dispersive equations of different orders in all dimensions.
Moreover, using the presented methods further results can
be obtained for some classes of non-dispersive equations
as well. In turn, many of the established results can
be also further extended to include small perturbations by
terms of the order up to the established smoothing
(e.g. for small magnetic potentials in the case of the
Schr\"odinger equation as in Georgiev and Tarulli \cite{GT},
or for the wave or Dirac equations as in e.g.
D'Ancona and Fanelli \cite{DF},
etc.) In non-dispersive cases they may involve certain
structural conditions on perturbations corresponding to the
invariant form of smoothing estimates established in
this paper. Moreover, global versions of the analytic
microlocal smoothing 
(as considered locally by using FBI transforms
by e.g. Robbiano and Zuily \cite{RZ} for Schr\"odinger, 
or Takuwa \cite{Ta} for dispersive equations, with 
a preceding work in the smooth setting by
Craig, Kappeler and Strauss \cite{CKS})
can be also expected to be obtained
by the presented methods. Smoothing of solutions to
Schr\"odinger equations has been also analysed by
Kapitanski and Safarov \cite{KaSa}, where a relation between
propagators of Schr\"odinger and wave equations was established.
In fact, in our setting
that relation can be also viewed as a composition of
a canonical transform and a comparison argument of this paper
(in particular using relation \eqref{prop:dim1ex-intro} below).
In this paper we will concentrate on the case of the case
of equations with ``constant coefficients", for the clarity and
comprehensiveness of the exposition. Further applications to
evolution equations for time and space dependent operators
will appear elsewhere.

\par
Now we will explain the essence of the approach and 
relate it to the known techniques. The main idea is that
instead of establishing smoothing estimates for different
classes of equations, we will instead relate such estimates
to each other. In particular, we will be able to relate
the majority of such estimates to a simple pointwise estimate
for a travelling wave in lower dimensions, which will follow
by a simple Fourier analysis argument. This will shed 
some light
on the nature of such estimates in a wider context, 
exhibiting a
quantitative smoothing phenomenon of the same type for
equations describing very different physical processes.
A comprehensive understanding of smoothing for equations
with constant coefficients should also 
make an impact on
problems with potentials and problems of perturbations of
such equations, an area which saw an amazingly 
rapid development over the recent years.
In particular, having an explicit relation between smoothing
estimates for constant coefficient equations can suggest
the corresponding admissible classes of potentials, etc.
Moreover, it should also influence the understanding of
smoothing--Strichartz estimates, a combined technique which
also proved to be very effective in nonlinear problems.

We will also suggest some invariant forms of the smoothing
estimates which we expect to continue to
hold in non-dispersive
cases as well, and we will give several results to justify
this expectation (these are estimates
\eqref{EQ:main-invariant}--\eqref{EQ:main-invariant4}).
We note that despite their natural appearance in many problems,
quite limited results are available for
non-dispersive equations while it is known that some usual
estimates fail in those case. In Section 
\ref{SECTION:invariant} we will discuss what smoothing estimates
are natural for such equations
and Section \ref{SECTION:nondisp} will be
devoted to the analysis of the non-dispersive cases.

As one of the simplest cases,
let us first consider the following Schr\"odinger equation:
\begin{equation}\label{eq:eq1}
\left\{
\begin{aligned}
\p{i\partial_t+\Delta_x}\,u(t,x)&=0\quad\text{in $\R_t\times\R^n_x$},\\
u(0,x)&=\varphi(x)\quad\text{in $\R^n_x$}.
\end{aligned}
\right.
\end{equation}
We know that the solution
operator $e^{it\Delta_x}$ preserves the $L^2$-norm 
for each fixed $t\in\R$.
On the other hand, the extra gain of regularity of order 
$1/2$ in $x$ can
be observed if we integrate the solution in $t$.
For example, in the case $n=1$, we have
\begin{equation}\label{eq:eq2}
   \n{|D_x|^{1/2}u(\cdot,x)}_{L^2(\R_t)}
  \le C\|\varphi\|_{L^2(\R)},
\end{equation}
for all $x\in\R$.
This result was given by e.g. Kenig, Ponce and Vega \cite{KPV1}.
Again, an application of the comparison principle in
Section \ref{SECTION:comparison} will allow us to
compare this estimate to the one-dimensional wave equation and
will show that  \eqref{eq:eq2}
is nothing else but the energy conservation
for the wave. 
In the one dimensional case this is again just the 
translation invariance of the Lebesgue measure.

In the higher dimensional case $n\geq2$, 
similar global smoothing properties are of importance:
\begin{equation}\label{eq:eq4}
\n{A u}_{L^2\p{\R_t\times\R^n_x}}\le C\n{\varphi}_{L^2\p{\R^n_x}},
\end{equation}
where $A$ is one of the following:
\begin{align*}
&(1)\quad A=\langle x\rangle^{-s}|D_x|^{1/2};\,\, s>1/2,
\\
&(2)\quad A=|x|^{\alpha-1}|D_x|^{\alpha};\,\,
 1-n/2<\alpha<1/2,
\\
&(3)\quad A=\langle x\rangle^{-s}\langle D_x\rangle^{1/2};\,\,
 s\geq1\quad (\text{$s>1$ if $n=2$}).
\end{align*}
We use the standard notation
\[
\jp{x}=\p{1+|x|^2}^{1/2} \quad\textrm{and}\quad
\jp{D_x}=\p{1-\Delta_x}^{1/2}.
\]
The type (1) was given by Ben-Artzi and Klainerman 
\cite{BK} ($n\geq3$), and
Chihara \cite{Ch} ($n\geq2$).
The type (2) was given by Kato and Yajima \cite{KY}
($n\geq3$, $0\leq\alpha<1/2$ or $n=2$, $0<\alpha<1/2$ ), 
and Sugimoto \cite{Su1}
($n\geq2$).
Watanabe \cite{W} showed that it is not true for $\alpha=1/2$.
The type (3) was given by Kato and Yajima \cite{KY} ($n\geq3$),
and Walther \cite{Wa1} ($n\geq2$) 
who also showed that it is not true 
for $s<1$ ($s\leq1$ if $n=2$). The type (4) estimate with
homogeneous weight and $\jp{D_x}^{1/2}$ was not
considered much but it
will be discussed in Sections \ref{SECTION:main} and \ref{SECTION:invariant}.
\par
Each proof was carried out by proving one of the following estimates
(or their variants):
\begin{equation}\label{eq:eq5}
\n{
\widehat{A^*f}_{|\rho \Sph^{n-1}}
}_{L^2\p{\rho \Sph^{n-1}}}
\le 
C\sqrt{\rho}
\n{f}_{L^2(\R^n)}
\qquad\text{(Restriction theorem)},
\end{equation}
where, $\rho \Sph^{n-1}=\b{\xi\in\Rn:\,|\xi|=\rho}$, ($\rho>0$),
or
\begin{equation}\label{eq:eq6}
 \sup_{\Im\zeta>0}\abs{\p{R(\zeta)A^*f,A^*f}} 
\leq
C \n{f}^2_{L^2(\R^n)}
\qquad\text{(Resolvent estimate)},
\end{equation}
where $R(\zeta)=\p{-\Delta_x-\zeta}^{-1}$.
Estimate \eqref{eq:eq5} implies the dual one of estimate \eqref{eq:eq4}.
Estimate \eqref{eq:eq6} implies \eqref{eq:eq4} 
since the resolvent $R(\zeta)$
is the Laplace transform of the solution operator 
$e^{it\Delta_x}$ of equation \eqref{eq:eq1}:
\[
R(\zeta)=\frac1i\int_0^\infty e^{it\Delta_x}e^{i\zeta t}\,dt\quad(\Im\zeta>0).
\]
The fact that \eqref{eq:eq6} implies \eqref{eq:eq5}
is due to the formula
\[
\Im\p{R(\rho^2+i0)f,f}
=\frac1{4(2\pi)^{n-1}\rho}\n{
\widehat{f}_{|\rho \Sph^{n-1}}
}_{L^2\p{\rho \Sph^{n-1}}}^2,
\]
see e.g. H\"ormander \cite[Corollary 14.3.10]{H}.
\par
In this paper we introduce several new ideas
to prove estimate \eqref{eq:eq4}. The main two proposed 
methods 
(canonical transforms and comparison principles) are centred
at comparing different estimates rather than looking
at them individually. This approach will allow us to actually
relate most of estimates to each other as well as to their normal
forms. For example, we will show that estimates \eqref{eq:eq4}
with $A$ as in (1), (2), or (3), are equivalent
to some simple one dimensional estimates.
To explain this idea, let us first
recall that operators other than the Schr\"odinger
operator have also attracted much attention for their
smoothing properties. For example, relativistic Schr\"odinger
equations have been investigated in \cite{BN} and \cite{Wa2},
wave and Klein--Gordon equations in \cite{Be},
Korteveg--de Vries equations in \cite{KPV2}, Benjamin--Ono
equations in \cite{KPV4}, Davey--Stewartson systems in
\cite{LP}, certain dispersive polynomial equations in \cite{BD2},
third order differential equations in \cite{KoSa},
to mention a few, and they can be 
expressed in the general form
\begin{equation}\label{EQ:Disp}
\left\{
\begin{aligned}
\p{i\partial_t+a(D_x)}\,u(t,x)&=0\quad\text{in $\R_t\times\R^n_x$},\\
u(0,x)&=\varphi(x)\quad\text{in $\R^n_x$},
\end{aligned}
\right.
\end{equation}
where $a(\xi)$ is a real-valued
function of $\xi=(\xi_1,\ldots,\xi_n)$
with the growth of order $m$,
and $a(D_x)$ is the corresponding operator.
Equations of this type have been extensively
studied under the ellipticity ($a(\xi)\not=0$ for $\xi\not=0$)
or the dispersiveness 
($\nabla a(\xi)\not=0$ for $\xi\not=0$) conditions.
Under such conditions, various global smoothing estimates
have been established for solutions $u(t,x)=e^{ita(D_x)}\varphi(x)$
in many papers, in both differential
and pseudo-differential cases (\cite{BN}, \cite{BD2},
\cite{Ch}, \cite{CS}, \cite{Ho1}, \cite{Ho2}, \cite{KY},
\cite{KPV1}, \cite{RS1},  \cite{Wa2}, etc.).
The dispersiveness condition was shown 
to be necessary for certain
types of estimates (see Hoshiro \cite{Ho2}), 
but we will show how
to get around that.
Now, suppose that we want to establish a weighted smoothing
estimate of the form 
\begin{equation}\label{EQ:comp-idea1}
\n{w(x)\rho(D_x) e^{ita(D_x)}\va(x)}_\L2tx\leq C
\n{\va}_{L^2(\R^n_x)},
\end{equation}
giving a smoothing of type $\rho(D_x)$ with some weight $w(x)$.
The rough idea of the canonical transform method is to use
certain operators $T$ for which we have the relations
\[
a(D_x)\circ T=T\circ\widetilde{a}(D_x)
\quad \textrm{and} \quad 
\rho(D_x)\circ T=T\circ \widetilde{\rho}(D_x),
\] 
for some other operators $\widetilde{a}(D_x)$ and 
$\widetilde{\rho}(D_x)$.
Then we also have
 $e^{it a(D_x)}\circ T=T\circ e^{it\widetilde{a}(D_x)}$.
We now substitute $T\va$ for $\va$ in estimate \eqref{EQ:comp-idea1},
and have
$$
\n{w(x)\rho(D_x) e^{ita(D_x)}T\va(x)}_\L2tx\leq C
\n{T\va}_{L^2(\R^n_x)}.
$$
Using the above identities we can conclude that estimate
\eqref{EQ:comp-idea1} is equivalent to the estimate
\begin{equation}\label{EQ:comp-idea4}
\n{w(x)T\widetilde{\rho}(D_x) e^{it\widetilde{a}(D_x)}
\va(x)}_\L2tx\leq C
\n{T\va}_{L^2(\R^n_x)}.
\end{equation}
If now operators $T$ and $T^{-1}$ are bounded in $L^2(\R^n_x)$
and in weighted $L^2(\R^n_x)$ with weight $w(x)$ respectively, we can remove
them from \eqref{EQ:comp-idea4} to finally conclude that
weighted smoothing estimate \eqref{EQ:comp-idea1} is equivalent
to 
\begin{equation}\label{EQ:comp-idea5}
\n{w(x)\widetilde{\rho}(D_x) e^{it\widetilde{a}(D_x)}
\va(x)}_\L2tx\leq C
\n{\va}_{L^2(\R^n_x)}.
\end{equation}
This idea can be used in a variety of ways. Not only can we derive
one smoothing estimate from another, but we can also consider
equivalent classes of smoothing estimates and find
their normal forms where the analysis would follow from
some straightforward argument. Thus, Section 
\ref{SECTION:model} will be devoted to estimates in such
model cases, while Section \ref{SECTION:canonical} will be
devoted to weighted estimates for necessary operators $T$.
Moreover, in Sections \ref{SECTION:comparison}
we will develop comparison principles
which will also allow us to relate the model estimates
from different classes
among each other.  All of the arguments will be invertible,
thus establishing a more or less complete set of relations
among different types of smoothing estimates
(dispersive in Section \ref{SECTION:main} and
non-dispersive in Section \ref{SECTION:nondisp}). In addition,
in Section \ref{SECTION:model} we will also relate estimates
with different weights.
\par
As for transformation operators $T$ and $T^{-1}$, we will
consider Fourier integral operators, or rather operators
which can be globally
written in the form
\begin{equation}\label{EQ:FIO}
      Tu(x)=\p{2\pi}^{-n}\int_{\R^n}\int_{\R^n}
      e^{i\Phi(x,y,\xi)}p(x,y,\xi)u(y)\, dy d\xi
\quad
(x\in\R^n),
\end{equation}
where $p(x,y,\xi)$ is an amplitude function and 
$\Phi(x,y,\xi)$ is a real phase function (not always
positively homogeneous in $\xi$ in our applications). 
Especially, if $p(x,y,\xi)=1$ and  $\Phi(x,y,\xi)$ satisfies
the graph condition
\begin{align*}
\Lambda
&=\b{(x,\Phi_x,y,-\Phi_y);\,\Phi_\xi=0}
\\
&=\b{(x,\xi),\chi(x,\xi)}
\subset T^*\R^n\times T^*\R^n,
\end{align*}
we have the relation
\begin{align*}
&T\circ A(X,D_x)\circ T^*=B(X,D_x)+(lower\; order\; terms),
\\
&B(x,\xi)= (A\circ\chi)(x,\xi),
\end{align*}
for pseudo-differentiable operators 
$A(X,D_x)$ and $B(X,D_x)$.
In this way, Fourier integral operators are recognised 
as a tool of the
realisation of the canonical transformation.
This fact is well-known microlocally as Egorov's theorem,
and by taking phase function appropriately, properties of the
operator $B(X,D_x)$ can be extracted from those of the
operator $A(X,D_x)$.
By using Egorov's theorem, many {\em qualitative} 
properties of solutions
of partial differential equations 
(propagation of singularities, construction of parametrises,
etc.)
have been investigated.
Our main interest is to establish {\em quantitative} 
properties as well
(global $L^2$-property for example) by the same idea.
In this paper, we take
\begin{equation}\label{EQ:PHASE}
\Phi(x,y,\xi)=x\cdot\xi-y\cdot\psi(\xi)
\end{equation}
and use the exact relation
\begin{equation}\label{relation}
T\circ\sigma(D_x)=a(D_x)\circ T,\qquad
a(D_x)=\p{\sigma\circ\psi}(D_x),
\end{equation}
for translation invariant
pseudo-differential operators $\sigma(D_x)$ and $a(D_x)$.
For example, the Laplacian $\Delta_x=\partial^2_{x_1}+\cdots+\partial^2_{x_n}$
can be transformed to $\partial^2_{x_n}$
by choosing an appropriate $\psi(\xi)$, and
hence we will be able to reduce the smoothing estimate 
for Schor\"odinger equation \eqref{eq:eq1}
to the one dimensional estimate \eqref{eq:eq2}.
We note that since we will be working with operators with
constant coefficients we are able to perform the exact
global calculus, 
in comparison to the calculus modulo lower order or
smoothing terms provided by the Egorov's theorem. Moreover,
we will be using the exact inverse $T^{-1}$ rather than 
the adjoint $T^*$. The global $L^2$--boundedness of operators
\eqref{EQ:FIO} has been investigated before, for example
by Asada and Fujiwara \cite{AF}, Kumano-go \cite{Ku} and
Boulkhemair \cite{Bo1, Bo2}. 
Unfortunately, in all these papers an
assumption was made for the second order derivatives
matrix $\nabla^2_\xi\Phi(x,y,\xi)$ to be 
globally bounded in all variables, which clearly
fails for the phase \eqref{EQ:PHASE}. However, the global $L^2$
and also weighted $L^2$ boundedness theorems for Fourier integral
operators without such assumption are required for our
analysis. Some of these 
results have been established by the authors in \cite{RS2}
and some will be proved in Section \ref{SECTION:canonical}.
\par
It is remarkable that 
the method of canonical transformations described above
allows us to carry out a global microlocal reduction
of equation \eqref{EQ:Disp} to the
model cases $|\xi_n|^m$ (elliptic case)
or $\xi_1|\xi_n|^{m-1}$ (non-elliptic case)
under the dispersiveness condition.
For example, for equation
\eqref{EQ:Disp}
Chihara \cite{Ch} used involved spectral and harmonic
analysis and established
the estimate
\begin{equation}\label{EQ:chihara}
\n{\jp{x}^{-s}|D_x|^{(m-1)/2}e^{ita(D_x)}\varphi(x)}_{L^2\p{\R_t\times\R^n_x}}
\leq C\n{\varphi}_{L^2\p{\R^n_x}}
\qquad(s>1/2)
\end{equation}
in the case when $a(\xi)$ is positively homogeneous of order $m>1$.
With canonical transforms, this estimate
is easily reduced to low dimensional
pointwise estimates
\begin{align}
&\n{|D_x|^{(m-1)/2}e^{it|D_x|^m}\varphi(x)}_{L^2(\R_t)}
\leq C\n{\varphi}_{L^2(\R_x)},
\label{eq:ell}
\\
&\n{|D_y|^{(m-1)/2}e^{itD_x|D_y|^{m-1}}\varphi(x,y)}_{L^2(\R_t\times\R_y)}
\leq C\n{\varphi}_{L^2\p{\R^2_{x,y}}},
\label{eq:nonell}
\end{align}
for all $x\in\R$, respectively.
Note that estimate \eqref{eq:ell} with $m=2$ is estimate \eqref{eq:eq2}
for the Schr\"odinger equation in one dimension.
By establishing \eqref{eq:ell} and \eqref{eq:nonell} directly,
we will be able to immediately obtain \eqref{EQ:chihara}
for $m>0$, thus also including the hyperbolic case $m=1$, which
will be important for further analysis, in particular for the
understanding of the meaning of various estimates in terms of
the finite speed of propagation of singularities, etc.
The results which will be thus obtained on this path 
generalise and extend many known results
in the literature mentioned above. Moreover, 
this new idea gives us a clear comprehensive understanding
of the smoothing effects of dispersive equations.
\par
In addition, we will introduce another technique 
with which we can show that the comparison
of the symbols implies the same comparison of corresponding operators.
For example, in the one dimensional case,
if we have
\[
\frac{|\sigma(\xi)|}{|f^\prime(\xi)|^{1/2}}
\leq A \frac{|\tau(\xi)|}{|g^\prime(\xi)|^{1/2}}
\]
then we have automatically estimate
\[
\|\sigma(D_x)e^{it f(D_x)}\varphi(x)\|_{L^2(\R_t)}\leq A
\|\tau(D_x)e^{it g(D_x)}\varphi(x)\|_{L^2(\R_t)},
\]
for all $x\in\R$. 
This will, in turn, imply a variety of
weighted estimates.
It will also allow us to relate normal forms of estimates for
operators of different orders. As an example, let
us mention the following consequence for $n=1$ and $l,m>0$:
\begin{equation}\label{prop:dim1ex-intro}
\n{|D_x|^{(m-1)/2}e^{it|D_x|^{m}}
\varphi(x)}_{L^2(\R_t)}=
\sqrt{\frac{l}{m}}
\n{|D_x|^{(l-1)/2}e^{it|D_x|^{l}}
\varphi(x)}_{L^2(\R_t)}
\end{equation}
for every $x\in\R$,
assuming that $\supp\widehat{\varphi}\subset [0,+\infty)$ or
$(-\infty,0]$.
We will introduce this kind of comparison principles
in more general settings, which will prove to be 
another strong tool to induce
general estimates from simple ones.
Particularly, if we use the comparison principle in both
directions, we can show
the equivalence of many different smoothing estimates.
For example, 
using \eqref{prop:dim1ex-intro} with $l=1$,
we can show that estimate \eqref{eq:ell} or 
\eqref{eq:nonell}
is equivalent to the same estimate but just in 
the special case $m=1$.
This fact means that these two standard estimates can
in turn be 
derived from the
equality
\begin{equation}\label{eq:trivial}
\n{e^{itD_x}\varphi(x)}_{L^2(\R_t)}
=\n{\varphi}_{L^2\p{\R_{x}}}
\end{equation}
in the case $n=1$, which is just the conservation of energy 
for the
travelling wave in one dimension.
In this way, smoothing estimates
for dispersive equation \eqref{EQ:Disp} 
can be surprisingly reduced to
just a simple equality \eqref{eq:trivial}, 
which is a straightforward consequence of
the trivial fact $e^{itD_x}\varphi(x)=\varphi(x+t)$.
Thus, we can immediately recover the gain of $1/2$-derivatives
for the Schr\"odinger and of one derivative for the Korteweg-
de Vries equations (as in e.g. Kato \cite{Ka2}).
In this way we can actually reduce all dispersive smoothing
estimates to those for model hyperbolic, Schr\"odinger,
relativistic, KdV,  or other equations (whichever we prefer), 
or we can show that they
are all equivalent to each other. In addition, we will 
find some explicit best constants based on a constant found
by Simon \cite{Si} using Kato's theory \cite{Ka1}. In general,
we will concentrate on smoothing estimates with $L^2$--norms,
but the idea of comparison principle can be extended to
$L^p$--norms as well useful to Strichartz estimates
(see e.g. Corollary \ref{cor:Strichartz}). For example, 
it will immediately follow that for all $0<p\leq\infty$,
quantities
$||e^{it\sqrt{-\Delta}}\va||_{L^p(\Rnx,L^2(\R_t))}$,
$|||D_x|^{1/2}e^{-it\Delta}\va||_{L^p(\Rnx,L^2(\R_t))}$,
and $|||D_x|e^{it(-\Delta)^{3/2}}\va||_{L^p(\Rnx,L^2(\R_t))}$
for propagators of the wave, Schr\"odinger, and KdV 
type equations are all equivalent.
\par
On the other hand, 
coupled dispersive equations are of immense importance in 
applications while with only limited analysis available.
To give an example,
let $v(t,x)$ and $w(t,x)$
solve the following coupled system of Schr\"odinger equations:
\begin{equation}\label{eq2}
\left\{
\begin{aligned} 
i\partial_t v & =\Delta_x v+b(D_x)w, \\
i\partial_t w & =\Delta_x w+c(D_x)v,  \\
v(0,x) & =v_0(x), w(0,x)=w_0(x). 
\end{aligned}
\right.
\end{equation}
This is the simplest example of Schr\"odinger equations coupled
through linearised operators $b(D_x), c(D_x)$. 
Such equations appear
in many areas in physics. For example, this is a model of wave
packets with two modes (in the presence of resonances),
see Tan and Boyd \cite{TB}.
In fiber optics they appear to describe certain 
types of a pair 
of coupled modulated wave-trains (see e.g.
Manganaro and Parker \cite{MP}).
They also describe the field of optical solitons in fibres
(see Zen and Elim \cite{ZE}) as well as Kerr dispersion and stimulated
Raman scattering for ultrashort pulses transmitted through fibres.
In these cases the linearised operators $b$ and $c$ 
would be of
zero order. In models of optical pulse propagation of birefringent
fibres and in wavelength-division-multiplexed systems they are
of the first order (see Pelinovsky and Yang \cite{PY}). 
They may be of higher orders as well, for example
in models of optical solitons with higher order effects 
(see Nakkeeran \cite{Na}). A by now standard way to tackle nonlinear
versions of \eqref{eq2} are Strichartz and
smoothing estimates. We will give some examples
of such approach based on the critical case of one of the
smoothing estimates established in this paper. For example,
we will apply it to the global in time well-posedness of
derivative nonlinear Schr\"odinger equations with some
structural conditions. The details fall outside the
scope of this paper and will appear elsewhere.
\par
Suppose now that we are in the simplest 
situation when system 
\eqref{eq2} can be diagonalised. Its eigenvalues are 
$a_\pm(\xi)=-|\xi|^2\pm\sqrt{b(\xi)c(\xi)}$ and the system uncouples
into scalar equations of type \eqref{EQ:Disp}
with operators $a(D_x)=a_\pm(D_x)$.
Since the structure of operators $b(D_x), c(D_x)$ 
may be quite involved,
this motivates the study of scalar equations 
\eqref{EQ:Disp} 
with 
operators $a(D_x)$ of rather general form. Not only the presence
of lower order terms is important in time global problems, 
the principal part may be
rather general since we may have $\nabla a_\pm=0$ 
at some points.
In such situation we microlocalise around such points and
lose the structure (but not the properties) of the symbol
completely.

The combination of the proposed two new methods
(canonical transformations and the comparison principles)
however has a good power on the occasion of this analysis.
Besides the simplification of the proofs of smoothing estimates
for standard dispersive equations,
we have an advantage in treating rather general dispersive equations
where $a(\xi)$ admits lower order terms,
and also non-dispersive equations where the dispersiveness condition
$\nabla a(\xi)\neq 0$ breaks (Section \ref{SECTION:nondisp}).
The inclusion of lower order terms in the analysis is 
essential here since time global properties are dominated by
the low frequency part $|\xi|\leq R$ ($R>0$) of $a(\xi)$.
This fact is true even for the ordinary Schr\"odinger equation,
with the homogeneous Laplacian.
In such low frequency 
case we can not talk about the principal part of
an operator, so operators with lower order terms appear
naturally, and will be considered as
condition (L) in Section \ref{SECTION:main}. 
\par
We also suggest an invariant form of smoothing estimates
which remain valid also in some
areas without dispersion, where standard
smoothing estimates are known to fail.
Let us observe the following form of estimates of 
types (1)--(3) for \eqref{eq:eq4}.
The first estimate may be rewritten in the form
\begin{equation}\label{EQ:main-invariant}
\n{\jp{x}^{-s}|\nabla a(D_x)|^{1/2}
e^{it a(D_x)}\varphi(x)}_{L^2\p{\R_t\times\R^n_x}}
\leq C\n{\varphi}_{L^2\p{\R^n_x}} \quad {\rm(}s>1/2{\rm)}.
\end{equation}
An analogous invariant forms for the other smoothing estimates
are estimate
\begin{equation}\label{EQ:main-invariant3}
\n{\abs{x}^{\alpha-m/2}|\nabla a(D_x)|^{\alpha/(m-1)}e^{ita(D_x)}\varphi(x)}_
{L^2\p{\R_t\times\R^n_x}}
\leq C\n{\varphi}_{L^2\p{\R^n_x}}
\quad(m\neq1),
\end{equation}
which we expect to hold for
$(m-n)/2<\alpha<(m-1)/2$, with $m>0$ being the order of $a(D_x)$,
and estimate
\begin{equation}\label{EQ:main-invariant4}
\n{\jp{x}^{-m/2}\jp{\nabla a(D_x)}^{1/2}e^{ita(D_x)}\varphi(x)}_
{L^2\p{\R_t\times\R^n_x}}
\leq C\n{\varphi}_{L^2\p{\R^n_x}},
\end{equation}
for $n>m>1$
because of the low frequency contribution. 
Such estimates have a number of advantages which we would
like to point out taking \eqref{EQ:main-invariant} as an
example:
\begin{itemize}
\item in the dispersive case it is equivalent to the usual
estimate \eqref{EQ:main1-0} below;
\item it does continue to hold for a variety of non-dispersive
equations, where $\nabla a(\xi)$ may become zero on some set
and when \eqref{EQ:main1-0} fails;
\item it does take into account zeros of the gradient
$\nabla a(\xi)$, which is also responsible for the interface
between dispersive and non-dispersive zone (e.g. how quickly
the gradient vanishes);
\item it is invariant under canonical transforms of the equation;
\item the proposed estimates are scaling 
invariant or ``almost'' invariant;
\item the estimates are also ``sharp'' (see Section 
\ref{SECTION:invariant}).
\end{itemize}
We observe the estimate
\begin{equation}\label{EQ:main1-0}
\n{\jp{x}^{-s}|D_x|^{(m-1)/2}e^{ita(D_x)}
\varphi(x)}_{L^2\p{\R_t\times\R^n_x}}
\leq C\n{\varphi}_{L^2\p{\R^n_x}} 
\end{equation}
which is known from e.g. Theorems \ref{M:H1} and 
Corollary \ref{M:L5}.
An additional
advantage of using estimate \eqref{EQ:main-invariant}
rather than \eqref{EQ:main1-0} is that \eqref{EQ:main-invariant}
takes into account possible zeros of the gradient
$\nabla a(\xi)$ is the non-dispersive case.
Thus,
in the one dimensional or in the 
radially symmetric cases of Section
\ref{SECTION:comparison} we will see that estimate
\eqref{EQ:main-invariant} is still valid in non-dispersive
cases. In Section \ref{SECTION:nondisp} we will show
examples of this estimate in other non-dispersive cases.
In particular, we will justify the invariant estimates
above in several situations using the following ideas:
\begin{itemize}
\item we can microlocalise around non-dispersive points and
apply the canonical transforms there
(Theorem \ref{THM:isolated-critical});
\item in radially symmetric cases we can use the comparison
principle (Theorems \ref{prop:dim1eqmod},
\ref{Th:nondisprad}, and Corollary \ref{COR:RStype}) ;
\item in some quasi-homogeneous cases or when  the symbol can
be represented as a sum of one dimensional monomials,
we can use the comparison principle as well
(Examples \ref{nondispersive1}, \ref{nondispersive2}
and \ref{nondispersive-shrira});
\item in the homogeneous case with some information on the
Hessian we can use canonical transforms to reduce the
general case to the previous situation 
(Theorem \ref{th1}).
\end{itemize}
In several situations estimate \eqref{EQ:main-invariant4} will
be weaker than \eqref{EQ:main-invariant} and
\eqref{EQ:main-invariant3}, so we may mostly
concentrate on these two. Moreover, in Sections
\ref{SECTION:canonical},
\ref{SECTION:main}, and \ref{SECTION:invariant}, we will
argue that in various complicated situations (like
in general non-dispersive cases) estimates
\eqref{EQ:main-invariant}--\eqref{EQ:main-invariant4} are
the (sharp) smoothing estimates that one can hope to obtain.
In addition, we will derive estimates for equations with time
dependent coefficients. In general, the dispersive estimates
for equations with time dependent coefficients may be a
delicate problem, with decay rates heavily depending on the
oscillation in coefficients (for a survey of different
results for the wave equation with lower order terms see,
e.g. Reissig \cite{Rei}). However, we will show in Section
\ref{SECTION:invariant} that the smoothing estimates still
remain valid if we introduce an appropriate factor into the
estimate. Such estimates become a natural extension of the
invariant estimates to the time dependent setting.
\par
We will explain the organisation of this paper.
In Section \ref{SECTION:comparison}, we give the
precise statements of the comparison principle.
There we will also give an example for Strichartz type norms
in Corollary \ref{cor:Strichartz}.
In Section \ref{SECTION:model}, 
we prove important model estimates and
also the equivalence of them by using the comparison principle.
We will also apply the comparison principle again 
to compare many estimates
with the estimates given here,
and get secondary comparison results.
In Section \ref{SECTION:canonical}, we introduce and show
the fundamental properties of our main tools which
originate in the idea of canonical transformation.
In Section \ref{SECTION:main}, we list results
which extend and explain estimate \eqref{eq:eq4} 
with types (1)--(3), which were partially announced 
by the authors in
\cite{RS1} and \cite{RS4}.
Especially, these kinds of time-global estimate 
for the operator $a(D_x)$
with lower order terms are new results
provided by the new method.
We also explain how general cases can be reduced
to the model estimates given in Section \ref{SECTION:model}.
Additional arguments 
with the idea of canonical transformation are also
presented there.
A second comparison result for radially symmetric case is
also given there.
In Section \ref{SECTION:invariant} we will propose and
discuss an invariant form of smoothing estimates which 
remains to hold in non-dispersive situations as well,
and we also discuss the sharpness of all the estimates.
The case of time--dependent coefficients will be treated
in Section \ref{SECTION:time-dependent}.
In Section \ref{SECTION:nondisp}, we will establish invariant
estimates for several case of non-dispersive equations
by using the second comparison results.
In Section \ref{SECTION:relative} we apply the second comparison
result further to the relativistic Schr\"odinger,
Klein--Gordon, and wave equations.
Sections \ref{SECTION:model-inh} and 
\ref{SECTION:inhomogeneous} are
devoted to non-homogeneous problems as a counterpart of Sections
\ref{SECTION:model} and \ref{SECTION:main}, respectively.
Section \ref{SECTION:other} is devoted to related problems,
including critical cases of some of the estimates, and the
corresponding trace theorems.
\par
Finally we comment on the notation used in this paper.
As usual, we will denote $D_{x_j}=-i\partial_{x_j}$ and
view operators $a(D_x)$ as Fourier multipliers.
Constants denoted by letter $C$ 
in estimates are always positive and
may differ on different 
occasions, but will still be denoted by the same letter.

\section{Comparison principle}
\label{SECTION:comparison}
In this section we will introduce a useful tool to derive 
new smoothing estimates
from known ones and to relate different estimates for solutions to
different equations with each other. We will concentrate 
on smoothing estimates with $L^2$--norms, 
and then will also give an application
to Strichartz type norms in Corollary \ref{cor:Strichartz}.

Thus, we will present a comparison principle for solutions
$u(t,x)=e^{it f(D_x)}\varphi(x)$ and $v(t,x)=e^{it g(D_x)}\varphi(x)$
to evolution equations with operators $f(D_x)$ and
$g(D_x)$, where $t\in\R$ and $x\in\R^n$:
\[
\left\{
\begin{aligned}
\p{i\partial_t+f(D_x)}\,u(t,x)&=0,\\
u(0,x)&=\varphi(x),
\end{aligned}
\right.
\quad {\rm and}\quad
\left\{
\begin{aligned}
\p{i\partial_t+g(D_x)}\,v(t,x)&=0,\\
v(0,x)&=\varphi(x).
\end{aligned}
\right.
\]
In the sequel, we write 
$x=(x_1,\ldots,x_n)$, $\xi=(\xi_1,\ldots,\xi_n)$, and
$D_x=(D_1,D_2\ldots,D_n)$ where
$D_j$ denotes $D_{x_j}=\frac{1}{i}\frac{\partial}
{\partial x_j}$, ($j=1,2,\ldots,n$).
\par
First we note the following fundamental result:
\medskip
\begin{thm}\label{prop:basiceq}
Let $f\in C^1(\R^n)$ be a real-valued function
such that, for almost all $\xi'=(\xi_2,\ldots,\xi_n)\in\R^{n-1}$,
$f(\xi_1,\xi')$ is strictly monotone in $\xi_1$ 
on the support of a measurable function $\sigma$ on $\R^n$.
Then we have
\begin{equation}\label{EQ:basiceq}
\n{\sigma(D_x)
e^{it f(D_x)}\varphi(x_1,x')}_{L^2(\R_t\times\R_{x'}^{n-1})}^2
= (2\pi)^{-n}
 \int_\Rn |\widehat{\varphi}(\xi)|^2 \frac{|\sigma(\xi)|^2}
 {|\partial f/\partial\xi_1(\xi)|}\, d\xi
\end{equation}
for all $x_1\in\R$, where $x'=(x_2,\ldots,x_n)\in\R^{n-1}$.
\end{thm}
\medskip
\begin{proof}
Let $\eta=\Phi(\xi)$ and $\xi=\Phi^{-1}(\eta)$ be changes
of variables defined by
$$
\Phi(\xi)=(f(\xi),\xi'); \quad
\Phi^{-1}(\eta)=(s(\eta),\eta'),
$$
where
we write $\eta=(\eta_1,\eta')$, $\eta'=(\eta_2,\ldots,\eta_n)$.
We assume that all the integrals below make sense which can be justified
in an usual manner using the assumption and Sard's theorem.
In view of this we perform calculations
on the set $|\partial \Phi(\xi)|=|\partial f/\partial \xi_1(\xi)|\not=0$.
We have
\begin{align*}
 & \sigma(D_x)e^{itf(D_x)}\varphi(x) \\
 = & (2\pi)^{-n}\int_\Rn e^{itf(\xi)} e^{ix\cdot\xi}
 \sigma(\xi)\widehat{\varphi}(\xi) d\xi \\
 = &(2\pi)^{-n}
\int_{\Phi(\Rn)} e^{i(t\eta_1+x'\cdot\eta')}e^{ix_1 s(\eta)}
 \sigma(\Phi^{-1}(\eta))
  \widehat{\varphi}(\Phi^{-1}(\eta))
   |\partial\Phi^{-1}(\eta)|\, d\eta,
\end{align*}
where we used the substitution $\xi=\Phi^{-1}(\eta)$ on the
support of $\chi$. Using Plancherel's identity, we get
\begin{align*}
& \|\sigma(D_x)e^{itf(D_x)}\varphi(x)\|^2
_{L^2(\R_t\times\R_{x'}^{n-1})} \\
= & (2\pi)^{-n}
\int_{\Phi(\Rn)} \abs{\sigma(\Phi^{-1}(\eta))
  \widehat{\varphi}(\Phi^{-1}(\eta))}^2
   \abs{\partial\Phi^{-1}(\eta)}^2\, d\eta
\\
= & (2\pi)^{-n}
\int_\Rn \abs{\sigma(\xi)
  \widehat{\varphi}(\xi)}^2
   \abs{\partial\Phi^{-1}(\Phi(\xi))}^2|\partial\Phi(\xi)|\, d\xi
\\
 = &(2\pi)^{-n}
 \int_\Rn |\widehat{\varphi}(\xi)|^2 \frac{|\sigma(\xi)|^2}
 {|\partial f/\partial\xi_1(\xi)|}\, d\xi,
\end{align*}
where we have used the substitution $\eta=\Phi(\xi)$
and the identity
$|\partial\Phi^{-1}(\Phi(\xi))|=|\partial\Phi(\xi)|^{-1}=
|\partial f/\partial\xi_1(\xi)|^{-1}$.
Note that this quantity is independent of $x_1$,
finishing the proof of \eqref{EQ:basiceq}. 
\end{proof}
\medskip
The following comparison principle
is a straightforward consequence of Theorem \ref{prop:basiceq}:
\medskip
\begin{cor}\label{prop:dimneq}
Let $f,g\in C^1(\R^n)$ be real-valued functions
such that, for almost all $\xi'=(\xi_2,\ldots,\xi_n)\in\R^{n-1}$,
$f(\xi_1,\xi')$ and $g(\xi_1,\xi')$
are strictly monotone in $\xi_1$ 
on the support of a measurable function $\chi$ on $\R^n$.
Let $\sigma,\tau\in C^0(\R^n)$ be such that, for some $A>0$, we have
\begin{equation}\label{EQ:compassdimn}
\frac{|\sigma(\xi)|}{\abs{\partial_{\xi_1} f(\xi)}^{1/2}}
\leq A \, \frac{|\tau(\xi)|}{\abs{\partial_{\xi_1} g(\xi)}^{1/2}}
\end{equation} 
for all $\xi\in\supp\chi$ satisfying
$D_1 f(\xi)\not=0$ and 
$D_1 g(\xi)\not=0$.
Then we have
\begin{multline}\label{EQ:dimnest}
\n{\chi(D_x)\sigma(D_x)
e^{it f(D_x)}\varphi(x_1,x')}_{L^2(\R_t\times\R_{x'}^{n-1})} 
\\ \leq  A 
\|\chi(D_x)\tau(D_x)e^{it g(D_x)}\varphi(\widetilde x_1,x')\|_
{L^2(\R_t\times\R_{x'}^{n-1})}
\end{multline}
for all $x_1,\widetilde x_1\in\R$, where $x'=(x_2,\ldots,x_n)\in\R^{n-1}$.
Consequently, for any  measurable function $w$ on $\R$ we have
\begin{multline}\label{EQ:dimnestwgt}
\n{w(x_1)\chi(D_x)\sigma(D_x)
e^{it f(D_x)}\varphi(x)}_{L^2(\R_t\times\R_{x}^{n})} 
\\ \leq  A 
\|w(x_1)\chi(D_x)\tau(D_x)e^{it g(D_x)}\varphi(x)\|_
{L^2(\R_t\times\R_{x}^{n})}
\end{multline}
Moreover, if $\chi\in C^0(\Rn)$ and $w\neq0$ on a set
of $\R$ with positive measure , the converse is true,
namely, if we have estimate \eqref{EQ:dimnest} for all $\va$,
for some $x_1, \widetilde{x}_1\in\R$,
or if we have estimate \eqref{EQ:dimnestwgt} for all $\va$,
and the norms are finite, then we also have
inequality \eqref{EQ:compassdimn}.
\end{cor}
\medskip
We remark that the last inequality in Corollary \ref{prop:dimneq}
gives the comparison between
different weighted estimates.
The reason to introduce
function $\chi$ into the estimates is that the relation
between symbols may be different for different regions
of the frequencies $\xi$,
(for example this is the case for the relativistic Schr\"odinger
and for the Klein-Gordon equations which will be discussed in Section
\ref{SECTION:relative}), so we have freedom to
choose different $\sigma$ for different types of behaviour
of $f^\prime$.
The assumption $\sigma,\tau\in C^0(\R^n)$ made there is for
the clarity of the exposition and can clearly be
relaxed. We will not need it in this paper, but if 
$\sigma$ and $\tau$ are simply measurable, satisfy
\eqref{EQ:compassdimn} almost everywhere, and if all the
integrals make sense, the conclusion of Corollary
\ref{prop:dimneq} and subsequent results continue to hold.
\par
In the case $n=1$, we neglect $x'=(x_2,\ldots,x_n)$ 
in a natural way and
just write $x=x_1$, $\xi=\xi_1$, and $D_x=D_1$.
Similarly in the case $n=2$, we use the notation $(x,y)=(x_1,x_2)$,
$(\xi,\eta)=(\xi_1,\xi_2)$, and $(D_x,D_y)=(D_1,D_2)$.
In both cases, we write $\widetilde x=\widetilde x_1$ 
in notation of
Corollary \ref{prop:dimneq}.
Then we have the following corollaries:
\medskip
\begin{cor}\label{prop:dim1eq}
Suppose $n=1$.
Let $f,g\in C^1(\R)$ be real-valued and strictly
monotone on the support of a measurable function $\chi$ on $\R$.
Let $\sigma,\tau\in C^0(\R)$ be such that, for some $A>0$, we have
\begin{equation}\label{dim1eq:1}
\frac{|\sigma(\xi)|}{|f^\prime(\xi)|^{1/2}}
\leq A \frac{|\tau(\xi)|}{|g^\prime(\xi)|^{1/2}}
\end{equation}
for all $\xi\in\supp\chi$ satisfying
$f^\prime(\xi)\not=0$ and $g^\prime(\xi)\not=0$.
Then we have
\begin{equation}\label{dim1eq:2}
\|\chi(D_x)\sigma(D_x)e^{it f(D_x)}\varphi(x)\|_{L^2(\R_t)}\leq A
\|\chi(D_x)\tau(D_x)e^{it g(D_x)}\varphi(\widetilde x)\|_{L^2(\R_t)}
\end{equation}
for all $x,\widetilde x\in\R$.
Consequently, for general $n\geq1$ and
for any measurable function $w$ on $\Rn$, we have
\begin{multline}\label{dim1eq:3}
\|w(x)\chi(D_j)\sigma(D_j)e^{it f(D_j)}\varphi(x)\|_
{L^2(\R_t\times\R_x^n)} \\ \leq A
\|w(x)\chi(D_j)\tau(D_j)e^{it g(D_j)}\varphi(x)\|_
{L^2(\R_t\times\R_x^n)},
\end{multline}
where $j=1,2,\ldots,n$.
Moreover, if $\chi\in C^0(\R)$ and $w\neq0$ on a set
of $\R^n$ with positive measure,
the converse is true,
namely, if we have estimate \eqref{dim1eq:2} for all $\va$,
for some $x, \widetilde{x}\in\R$,
or if we have estimate \eqref{dim1eq:2} for all $\va$,
and the norms are finite, then we also have
inequality \eqref{dim1eq:1}.
\end{cor}
\medskip
\begin{cor}\label{prop:dim2eq}
Suppose $n=2$.
Let $f,g\in C^1(\R^2)$ be real-valued functions
such that, for almost all $\eta\in\R$,
$f(\xi,\eta)$ and $g(\xi,\eta)$ are strictly
monotone in $\xi$ on the support of a measurable function $\chi$ on $\R^2$.
Let $\sigma,\tau\in C^0(\R^2)$ be such that, for some $A>0$, we have
\begin{equation}\label{dim2eq:1}
\frac{|\sigma(\xi,\eta)|}
{\abs{\partial f/\partial \xi(\xi,\eta)}^{1/2}}
\leq A \frac{|\tau(\xi,\eta)|}
{\abs{\partial g/\partial \xi(\xi,\eta)}^{1/2}}
\end{equation}
for all $(\xi,\eta)\in\supp\chi$ satisfying
$\partial f/\partial \xi(\xi,\eta)\not=0$ and 
$\partial g/\partial \xi(\xi,\eta)\not=0$.
Then we have
\begin{multline}\label{dim2eq:2}
\n{\chi(D_x,D_y)\sigma(D_x,D_y)
e^{it f(D_x,D_y)}\varphi(x,y)}_{L^2(\R_t\times\R_{y})} \\ \leq A
\|\chi(D_x,D_y)\tau(D_x,D_y)e^{it g(D_x,D_y)}\varphi(\widetilde x,y)\|_
{L^2(\R_t\times\R_{y})}
\end{multline}
for all $x,\widetilde x\in\R$.
Consequently, for general $n\geq2$ and 
for any  measurable function $w$ on $\R^{n-1}$ we have
\begin{multline}\label{dim2eq:3}
\|w(\check{x}_k)\chi(D_j,D_k)\sigma(D_j,D_k)
e^{it f(D_j,D_k)}\varphi(x)\|_
{L^2(\R_t\times\R_x^n)} \\ \leq A
\|w(\check{x}_k)\chi(D_j,D_k)\tau(D_j,D_k)
e^{it g(D_j,D_k)}\varphi(x)\|_
{L^2(\R_t\times\R_x^n)},
\end{multline}
where $j\neq k$ and $\check{x}_k=(x_1,\ldots,x_{k-1},x_{k+1},\ldots,x_n)$.
Moreover, if $\chi\in C^0(\R^2)$ and $w\neq0$ on a set
of $\R^{n-1}$ with positive measure, the converse is true,
namely, if we have estimate \eqref{dim2eq:2} for all $\va$,
for some $x, \widetilde{x}\in\R$,
or if we have estimate \eqref{dim2eq:2} for all $\va$,
and the norms are finite, then we also have
inequality \eqref{dim2eq:1}.
\end{cor}
By the same argument as used in the proof of Theorem \ref{prop:basiceq}
and Corollary \ref{prop:dimneq},
we have a comparison result for radially symmetric case.
Below, we denote the set of the positive real numbers $(0,\infty)$
by $\R_+$.
\begin{thm}\label{prop:dim1eqmod}
Let $f,g\in C^1(\R_+)$ be real-valued and strictly
monotone on the support of a measurable function $\chi$ on $\R_+$.
Let $\sigma,\tau\in C^0(\R_+)$ be such that, for some $A>0$, we have
\begin{equation}\label{EQ:compassdimmod}
\frac{|\sigma(\rho)|}{|f^\prime(\rho)|^{1/2}}
\leq A \frac{|\tau(\rho)|}{|g^\prime(\rho)|^{1/2}}
\end{equation}
for all $\rho\in\supp\chi$ satisfying
$f^\prime(\rho)\not=0$ and $g^\prime(\rho)\not=0$.
Then we have
\begin{equation}\label{EQ:dim1estmod}
\|\chi(|D_x|)\sigma(|D_x|)e^{it f(|D_x|)}\varphi(x)\|_{L^2(\R_t)}
\leq A
\|\chi(|D_x|)\tau(|D_x|)e^{it g(|D_x|)}\varphi(x)\|_{L^2(\R_t)}
\end{equation}
for all $x\in\R^n$.
Consequently, for any measurable function $w$ on $\Rn$, we have
\begin{multline}\label{EQ:dim1estxcor}
\|w(x)\chi(|D_x|)\sigma(|D_x|)e^{it f(|D_x|)}\varphi(x)\|_
{L^2(\R_t\times\R_x^n)} \\ \leq A
\|w(x)\chi(|D_x|)\tau(|D_x|)e^{it g(|D_x|)}\varphi(x)\|_
{L^2(\R_t\times\R_x^n)}.
\end{multline}
Moreover, if $\chi\in C^0(\R_+)$ and $w\neq0$ on a set
of $\R^n$ with positive measure, the converse is true,
namely, if we have estimate \eqref{EQ:dim1estmod} for all $\va$,
for some $x\in\Rn$,
or if we have estimate \eqref{EQ:dim1estxcor} for all $\va$,
and the norms are finite, then we also have
inequality \eqref{EQ:compassdimmod}.
\end{thm}
\begin{proof}
Below, we will write $\xi=\rho\omega$,
where $\rho>0$ and $\omega\in \Sph^{n-1}$.
As usual we perform calculations
on the set $f^\prime(\rho)\not=0$, where the
inverse of $f$ is differentiable.
We have
\begin{align*}
 & \chi(|D_x|)\sigma(|D_x|)e^{itf(|D_x|)}\va(x) \\
 =&(2\pi)^{-n} \int_\Rn e^{itf(|\xi|)} e^{ix\cdot\xi}
 (\chi\sigma)(|\xi|)\widehat{\varphi}(\xi)\, d\xi \\
 =&(2\pi)^{-n} 
 \int_{\R_+} \int_{\Sph^{n-1}}  
 e^{itf(\rho)} e^{i\rho x\cdot\omega}(\chi\sigma)(\rho)
  \widehat{\varphi}(\rho\omega) \rho^{n-1}\,d\rho
 d\omega \\
 =&(2\pi)^{-n} 
 \int_{f(\R_+)}   \int_{\Sph^{n-1}} e^{it\eta} e^{if^{-1}(\eta)x\cdot\omega}
 (\chi\sigma)(f^{-1}(\eta))
  \widehat{\varphi}(f^{-1}(\eta)\omega) 
  f^{-1}(\eta)^{n-1}|(f^{-1})^\prime(\eta)|\, d\omega
 d\eta,
\end{align*}
where we used a substitution $\rho=f^{-1}(\eta)$ on the
support of $\chi$. Using Plancherel's identity, we get
\begin{equation}\label{EQ:bacsiceqrod}
\begin{aligned}
& \|\chi(|D_x|)\sigma(|D_x|)e^{itf(|D_x|)}\va(x)
\|^2_{L^2(\R_t)} \\
= & (2\pi)^{-2n+1} \int_{f(\R_+)}\,d\eta\;\times \\
 &\qquad\times
\left|
 \int_{\Sph^{n-1}} e^{if^{-1}(\eta)x\cdot\omega}
(\chi\sigma)(f^{-1}(\eta)) \widehat{\varphi}(f^{-1}(\eta)\omega)
f^{-1}(\eta)^{n-1}|(f^{-1})^\prime(\eta)|
d\omega \right|^2
       \\
= & (2\pi)^{-2n+1} 
\int_{\R_+} \left|
 \int_{\Sph^{n-1}} e^{i\rho x\cdot\omega}
(\chi\sigma)(\rho) \widehat{\varphi}(\rho\omega)
\rho^{n-1}|(f^{-1})^\prime(f(\rho))|
d\omega \right|^2 |f^\prime(\rho)|
  \,d\rho     \\
= & (2\pi)^{-2n+1} 
\int_{\R_+} \left|
 \int_{\Sph^{n-1}} e^{i\rho x\cdot\omega}
 \widehat{\varphi}(\rho\omega)
d\omega \right|^2\rho^{2(n-1)} |\chi(\rho)|^2
 \frac{|\sigma(\rho)|^2}{|f^\prime(\rho)|}
  \,d\rho,
\end{aligned}
\end{equation}
where we have used the substitution $\eta=f(\rho)$ again
and the identity
$(f^{-1})^\prime(f(\rho))=f^\prime(\rho)^{-1}$.
From assumption \eqref{EQ:compassdimmod}
it follows that
\begin{align*}
&\|\chi(|D_x|\sigma(|D_x|)e^{itf(|D_x|)}\va(x)\|^2_{L^2(\R_t)} \\
 \leq &
 (2\pi)^{-2n+1}A^2 \int_{\R_+} \left|
 \int_{\Sph^{n-1}} e^{i\rho x\cdot\omega}
 \widehat{\varphi}(\rho\omega)
d\omega \right|^2\rho^{2(n-1)} |\chi(\rho)|^2
 \frac{|\tau(\rho)|^2}{|g^\prime(\rho)|}
  d\rho  \\
 = & A^2
\|\chi(|D_x|)\tau(|D_x|)e^{itg(|D_x|)}\va(x)\|^2_{L^2(\R_t)},
\end{align*}
finishing the proof of \eqref{EQ:dim1estmod}. Estimate
\eqref{EQ:dim1estxcor} follows from it immediately.
The converse is also obtained from equality \eqref{EQ:bacsiceqrod}
which holds for any (radially symmetric) function $\va$.
\end{proof}

\medskip
{\bf Strichartz type norms.}
In fact, once we have estimate \eqref{EQ:dim1estmod},
we can take any further norm with respect to $x$. For example,
with Strichartz estimates in mind, we can take $L^p$ norms
as well.

\begin{cor}\label{cor:Strichartz}
Let functions $f,g,\sigma,\tau$ be as in Theorem 
{\rm \ref{prop:dim1eqmod}} and satisfy relation
{\rm \eqref{EQ:compassdimmod}}. Let $0<p\leq\infty$. Then, 
for any measurable function $w$ on $\Rn$, we have the
estimate
\begin{multline}\label{EQ:dim1estxcor-Str}
\|w(x)\chi(|D_x|)\sigma(|D_x|)e^{it f(|D_x|)}\varphi(x)\|_
{L^p(\Rnx,L^2(\R_t))} \\ \leq A
\|w(x)\chi(|D_x|)\tau(|D_x|)e^{it g(|D_x|)}\varphi(x)\|_
{L^p(\Rnx,L^2(\R_t))}.
\end{multline}
\end{cor}
We also note that if expressions on both sides of
\eqref{EQ:compassdimmod} are equivalent, we obtain the
equivalence of norms in \eqref{EQ:dim1estxcor-Str}.
For example, it immediately follows that for all $0<p\leq\infty$,
quantities
$||e^{it\sqrt{-\Delta}}\va||_{L^p(\Rnx,L^2(\R_t))}$,
$|||D_x|^{1/2}e^{-it\Delta}\va||_{L^p(\Rnx,L^2(\R_t))}$,
and 
$|||D_x|e^{it(-\Delta)^{3/2}}\va||_{L^p(\Rnx,L^2(\R_t))}$
for propagators of the wave, Schr\"odinger, and KdV 
type equations
are equivalent.

By an easy application of Minkowski's inequality for integrals,
we have inequalities 
$$  ||f||_{L^2(\R_t,L^{p_1}(\Rnx))} \leq
   C||f||_{L^{p_1}(\Rnx,L^2(\R_t))}, \quad
    ||f||_{L^{p_2}(\Rnx,L^2(\R_t))}\leq 
   C ||f||_{L^2(\R_t,L^{p_2}(\Rnx))},
$$
for $p_1\leq 2\leq p_2$, relating norms in 
\eqref{EQ:dim1estxcor-Str} to the usual Strichartz norms.
We also note that the $L^2$--norm in time is critical for a
variety of equations, and Strichartz estimates with $p=\infty$
may fail, so the smaller $L^{\infty}(\Rnx,L^2(\R_t))$--norms 
may be a good substitute in some situations. Among other
things this shows the equivalence of 
$L^{p}(\Rnx,L^2(\R_t))$--norms for different equations, similar
to the situation with smoothing estimates exhibited in this paper.
We will address these issues in more detail elsewhere.

\section{Equivalent model estimates}
\label{SECTION:model}
Let us now give important examples of the use of the comparison
principle described in Section \ref{SECTION:comparison}.
We still use the same notation as in Section \ref{SECTION:comparison}.
That is, denoting the dimension of the variable $x$ by $n$, we write
$x=(x_1,\ldots,x_n)$ and $D_x=(D_1,D_2\ldots,D_n)$.
We just write $x=x_1$, $D_x=D_1$ in the case $n=1$, 
and $(x,y)=(x_1,x_2)$, $(D_x,D_y)=(D_1,D_2)$ in the case $n=2$.
\par
If both sides in expression \eqref{EQ:compassdimn} 
in Corollary \ref{prop:dimneq}
are equivalent,
we can use the comparison in two directions, from which it
follows that norms on both sides in \eqref{EQ:dimnest}
are equivalent.
The same is true for Corollaries \ref{prop:dim1eq}, \ref{prop:dim2eq}
and Theorem \ref{prop:dim1eqmod}.
In particular, we can
conclude that many smoothing estimates for the Schr\"odinger
type equations of different orders are equivalent to each other.
Indeed, applying Corollary \ref{prop:dim1eq} in two directions, we
immediately obtain that for 
$n=1$ and $l,m>0$, we have
\begin{equation}\label{prop:dim1ex}
\n{|D_x|^{(m-1)/2}e^{it|D_x|^{m}}
\varphi(x)}_{L^2(\R_t)}=
\sqrt{\frac{l}{m}}
\n{|D_x|^{(l-1)/2}e^{it|D_x|^{l}}
\varphi(x)}_{L^2(\R_t)}
\end{equation}
for every $x\in\R$,
assuming that $\supp\widehat{\varphi}\subset [0,+\infty)$ or
$(-\infty,0]$.
Applying Corollary \ref{prop:dim2eq}, we
similarly obtain that for $n=2$ and $l,m>0$, we have 
\begin{multline}\label{prop:dim2ex}
\n{|D_y|^{(m-1)/2}e^{itD_x|D_y|^{m-1}}\varphi(x,y)}_{L^2(\R_t\times\R_y)}
\\ =
\n{|D_y|^{(l-1)/2}e^{itD_x|D_y|^{l-1}}\varphi(x,y)}_{L^2(\R_t\times\R_y)}
\end{multline}
for every $x\in\R$.
On the other hand, in the case $n=1$, we have easily
\begin{equation}\label{core}
\n{e^{itD_x}\varphi(x)}_{L^2(\R_t)}
=\n{\varphi}_{L^2\p{\R_{x}}}
\quad \text{for all $x\in\R$},
\end{equation}
which is a straightforward consequence of the fact
$e^{itD_x}\varphi(x)=\varphi(x+t)$.
By using equality \eqref{core},
we can estimate the right hand sides of equalities \eqref{prop:dim1ex}
and \eqref{prop:dim2ex} with $l=1$, and as a result, we have easily
the following variety of pointwise estimates in low dimensions:
\medskip
\begin{thm}\label{prop:basic}
Suppose $n=1$ and $m>0$.
Then we have
\begin{equation}\label{basic:1dim}
\n{|D_x|^{(m-1)/2}e^{it|D_x|^m}\varphi(x)}_{L^2(\R_t)}
\leq C\n{\varphi}_{L^2(\R_x)}
\end{equation}
for all $x\in\R$.
Suppose $n=2$ and $m>0$.
Then we have
\begin{equation}\label{basic:2dim}
\n{|D_y|^{(m-1)/2}e^{itD_x|D_y|^{m-1}}\varphi(x,y)}_{L^2(\R_t\times\R_y)}
\leq C\n{\varphi}_{L^2\p{\R^2_{x,y}}}
\end{equation}
for all $x\in\R$.
Each estimate above is equivalent to itself with $m=1$
which is a direct consequence of equality \eqref{core}.
In particular, we have equalities \eqref{prop:dim1ex} and \eqref{prop:dim2ex}.
\end{thm}
\medskip
Estimates \eqref{basic:1dim} and \eqref{basic:2dim}
in Theorem \ref{prop:basic} in the special case $m=2$
were shown by Kenig, Ponce and Vega \cite[p.56]{KPV1} and by
Linares and Ponce \cite[p.528]{LP}, respectively.
Theorem \ref{prop:basic} shows that
these results, together with their generalisation to
other orders $m$,
are in fact just corollaries of the elementary one 
dimensional fact
$e^{itD_x}\varphi(x)=\varphi(x+t)$ once we apply the comparison principle.
\par
By using the comparison principle in the radially symmetric case,
we have also another type of equivalence of
smoothing estimates.
In fact, by Theorem \ref{prop:dim1eqmod}, we immediately obtain
\begin{align*}
\n{|x|^{\beta-1} |D_x|^{\beta} e^{it |D_x|^2}\varphi}_\L2tx
&=
\sqrt{\frac{m}{2}}
\n{|x|^{\beta-1} |D_x|^{m/2+\beta-1} e^{it |D_x|^m}\varphi}_\L2tx
\\
&=
\sqrt{\frac{m}{2}}
\n{|x|^{\alpha-m/2} |D_x|^{\alpha} e^{it |D_x|^m}\varphi}_\L2tx,
\end{align*}
where $m>0$ and $\alpha=m/2+\beta-1$.
On the other hand, we know the estimate
\begin{equation}\label{core2'}
\n{\abs{x}^{\beta-1}|D_x|^{\beta}e^{it|D_x|^2}
\varphi(x)}_
{L^2\p{\R_t\times\R^n_x}}
\leq C\n{\varphi}_{L^2\p{\R^n_x}} \qquad
(1-n/2<\beta<1/2),
\end{equation}
which was given by Sugimoto \cite[Theorem 1.1]{Su1}.
Noticing that $1-n/2<\beta<1/2$ is equivalent to $(m-n)/2<\alpha<(m-1)/2$,
we have the estimate
\begin{multline}\label{core2g}
\n{\abs{x}^{\alpha-m/2}|D_x|^{\alpha}e^{it|D_x|^{m}}
\varphi(x)}_
{L^2\p{\R_t\times\R^n_x}}
\leq C\n{\varphi}_{L^2\p{\R^n_x}} \\
(m>0,\quad (m-n)/2<\alpha<(m-1)/2).
\end{multline}
We note that estimate \eqref{core2'} is a special case ($m=2$) of
estimate \eqref{core2g}, but the comparison principle of
Section \ref{SECTION:comparison} shows that
they are equivalent to each other.
\par
We remark that estimate \eqref{core2'} 
is implied from its restricted version
\begin{equation}\label{core2''}
\n{\abs{x}^{\beta-1}|D_x|^{\beta}e^{it|D_x|^2}
\varphi(x)}_
{L^2\p{\R_t\times\R^n_x}}
\leq C\n{\varphi}_{L^2\p{\R^n_x}} \qquad
(1/2-\varepsilon\leq\beta<1/2),
\end{equation}
where $\varepsilon>0$ is sufficiently small.
(The case $0<\varepsilon<1/2$ is the result of 
Kato and Yajima \cite{KY}, and the critical 
case of this estimate with $\epsilon=0$ was given in
\cite{Su2} and explained geometrically in \cite{RS3}).
In fact, estimate \eqref{core2'} with $1-n/2<\beta<1/2-\varepsilon$
can be reduced to
the one with $\beta=1/2-\varepsilon$
if we use the estimate
\[
\n{|x|^{\beta-1}|D_x|^\beta v}_{L^2(\R^n)}\leq C
\n{|x|^{(1/2-\varepsilon)-1}|D_x|^{1/2-\varepsilon} v}_{L^2(\R^n)}
\]
which is a consequence of the following lemma.
\begin{lem}[\cite{SW}, Theorem B$^*$]\label{Lem:SW}
Suppose $k<n/2$, $l <n/2$,
$0<m<n$, and
$k+l+m=n$.
Then the operator
$|x|^{-l}|D_x|^{m-n}|x|^{-k}$ is $L^2(\R^n)$--bounded.
\end{lem}
\par
Furthermore, we can show that in fact estimate \eqref{core2g}
is also equivalent to estimate
\begin{multline}\label{core2}
\n{
 \jp{x}^{-m/2}
  e^{it|D_x|^m}\varphi(x)
}_{L^2(\R_t\times\R^n_x)}
\leq C\n{\varphi}_{L^2(\R_x^n)}\qquad 
(n>m>1)
\\
\text{for all $\varphi$ such that $\supp\widehat{\varphi}
\subset\{\xi\in\R^n:\,|\xi|\leq1\}$}
\end{multline}
(given by Walther in \cite[Theorem 4.1]{Wa2}).
In fact,  estimate \eqref{core2} is
a direct consequence of estimate \eqref{core2g} with $\alpha=0$
if we notice a trivial inequality $\jp{x}^{-m/2}\leq|x|^{-m/2}$.
Note also that the assumption $n>m>1$ assures
$(m-n)/2<\alpha=0<(m-1)/2$.
On the other hand,
by Theorem \ref{prop:dim1eqmod} we have
\begin{align*}
&\n{
 \jp{x}^{\alpha-m/2}\chi(|D_x|)|D_x|^\alpha
  e^{it|D_x|^m}\varphi(x)
}_{L^2(\R_t\times\R^n_x)}
\\
= &
\sqrt{\frac{\mu}{m}}
\n{
 \jp{x}^{\alpha-m/2}\chi(|D_x|)|D_x|^{\alpha+(\mu-m)/2}
  e^{it|D_x|^\mu}\varphi(x)
}_{L^2(\R_t\times\R^n_x)}
\\
= &
\sqrt{\frac{m-2\alpha}{m}}
\n{
 \jp{x}^{-\mu/2}\chi(|D_x|)
  e^{it|D_x|^\mu}\varphi(x)
}_{L^2(\R_t\times\R^n_x)},
\end{align*}
where $m>0$ and $\mu=m-2\alpha>0$.
Hence, from estimate \eqref{core2} with $m=\mu$, we obtain
\[
\n{
 \jp{x}^{\alpha-m/2}\chi(|D_x|)|D_x|^\alpha
  e^{it|D_x|^m}\varphi(x)
}_{L^2(\R_t\times\R^n_x)}
\leq C\n{\varphi}_{L^2(\R_x^n)},
\]
where $n>m-2\alpha>1$, or equivalently $(m-n)/2<\alpha<(m-1)/2$.
Here we take a cut-off function 
$\chi(\rho)\in C^\infty_0([0,1))$ 
such that $\chi(\rho)\equiv1$ for $\rho\leq1/2$.
From this estimate, we obtain estimate \eqref{core2g}.
In fact,
we have the equality
\begin{multline*}
\n{
 |x|^{\alpha-m/2}
|D_x|^\alpha
  e^{it|D_x|^m}\varphi(x)}_{L^2(\R_t\times\R^n_x)}
\\
=
\lim_{\lambda\searrow0}
\n{
 \lambda^{\alpha-m/2}\jp{x/\lambda}^{\alpha-m/2}
\chi(\lambda|D_x|)|D_x|^\alpha
  e^{it|D_x|^m}\varphi(x)}_{L^2(\R_t\times\R^n_x)},
\end{multline*}
and noticing the identities
$\n{g(t,x)}_{L^2(\R_t\times\R^n_x)}
=\lambda^{m/2+n/2}\n{g(\lambda^m t, \lambda x)}_{L^2(\R_t\times\R^n_x)}$
and $(m(\lambda D_x)\varphi)(\lambda x)=m(D_x)(\varphi(\lambda\, \cdot))(x)$,
we have
\begin{multline*}
\n{
 |x|^{\alpha-m/2}
|D_x|^\alpha
  e^{it|D_x|^m}\varphi(x)}_{L^2(\R_t\times\R^n_x)}
\\
\leq
\sup_{\lambda>0}
\n{
 \jp{x}^{\alpha-m/2}
\chi(|D_x|)|D_x|^\alpha
  e^{it|D_x|^m}\varphi_\lambda(x)}_{L^2(\R_t\times\R^n_x)},
\end{multline*}
where $\varphi_\lambda(x)=\lambda^{n/2}\varphi(\lambda x)$.
Note also that $\n{\varphi_\lambda}_{L^2(\R^n_x)}=\n{\varphi}_{L^2(\R^n_x)}$.
\par
Finally we remark that the last inequality implies 
\[
\n{
 |x|^{\alpha-m/2}
|D_x|^\alpha
  e^{it|D_x|^m}\varphi(x)}_{L^2(\R_t\times\R^n_x)}
\leq
\sup_{\lambda>0}
\n{
 \jp{x}^{\alpha-m/2}
|D_x|^\alpha
  e^{it|D_x|^m}\varphi_\lambda(x)}_{L^2(\R_t\times\R^n_x)}
\]
by the comparison principle Theorem \ref{prop:dim1eqmod}.
Thus we can conclude the following:
\medskip
\begin{thm}\label{prop:basic2}
We have equivalent estimates \eqref{core2'}, \eqref{core2g}, 
and \eqref{core2}.
Furthermore, they are equivalent to estimate \eqref{core2''}
with sufficiently small $\varepsilon>0$.
In particular, for $m>0$ {\rm(}and any $\alpha,\beta${\rm)}
we have the following relations {\rm(}which are finite
for $\alpha,\beta$ as in the above estimates{\rm)}
\begin{align*}
&\n{|x|^{\beta-1} |D_x|^{\beta} e^{it |D_x|^2}\varphi}_\L2tx
 =
\sqrt{\frac{m}{2}}
\n{|x|^{\beta-1} |D_x|^{m/2+\beta-1} e^{it |D_x|^m}\varphi}_\L2tx,
\\
&\begin{aligned}
\n{
 \jp{x}^{\alpha-m/2}|D_x|^\alpha
  e^{it|D_x|^m}\varphi(x)
}_{L^2(\R_t\times\R^n_x)}
&\leq 
\n{
 |x|^{\alpha-m/2}|D_x|^\alpha
  e^{it|D_x|^m}\varphi(x)
}_{L^2(\R_t\times\R^n_x)} \\
& \leq 
\sup_{\lambda>0}\n{
 \jp{x}^{\alpha-m/2}|D_x|^\alpha
  e^{it|D_x|^m}\varphi_\lambda(x)
}_{L^2(\R_t\times\R^n_x)},
\end{aligned}
\end{align*}
where $\varphi_\lambda(x)=\lambda^{n/2}\varphi(\lambda x)$, and
we take $\alpha\leq m/2$ in the last estimate.
The operator norms of operators
$\jp{x}^{\alpha-m/2} |D_x|^\alpha e^{it|D_x|^m}$ and
$|x|^{\alpha-m/2} |D_x|^\alpha e^{it|D_x|^m}$ as mappings from
$L^2(\Rn)$ to $\L2tx$ are equal.
\end{thm}
As a nice consequence, for $n\geq 3$ and $m>0$
we can conclude also the estimate
\begin{equation}\label{EQ:Simon-const}
\n{|x|^{-1} |D_x|^{m/2-1}
  e^{it|D_x|^m}\varphi(x)
}_{L^2(\R_t\times\R^n_x)}\leq
\sqrt{\frac{2\pi}{m(n-2)}}\n{\va}_{L^2(\R_x^n)},
\end{equation}
where the constant $\sqrt{\frac{2\pi}{m(n-2)}}$ is sharp.
This follows from the first equality in Theorem \ref{prop:basic2}
with $\beta=0$ and the fact that the constant
$C=\sqrt{\frac{\pi}{n-2}}$ is sharp in \eqref{core2'} with
$\beta=0$, as shown by Simon \cite{Si} as a 
consequence of constants in Kato's theory \cite{Ka1}.

In general, best constants in the radially symmetric case
can be obtained by changing to spherical harmonics and
looking at the appearing one dimensional integral. Thus, 
if $n\geq 2$ and $f$ is injective and
differentiable on $(0,\infty)$,
the best constant in the inequality
$$
\n{w(|x|)\sigma(|D_x|)e^{it f(|D_x|)}\va(x)}_\L2tx\leq
C\n{\va}_{L^2(\R^n_x)}
$$
is given by
$$
C=(2\pi)^{(n+1)/2}\p{\mathop{\sup_{\rho>0}}_{k\in\N}
\b{\rho\sigma(\rho)^2 f^\prime(\rho)^{-1}
\int_0^\infty J_{\nu(k)}(r\rho)^2 w(r)^2 r dr}}^{1/2},
$$
where for $\lambda>-1/2$ the Bessel function 
$J_\lambda$ of order
$\lambda$ is given by
$$
J_\lambda(\rho)=\frac{\rho^\lambda}{2^\lambda
\Gamma(\lambda+1/2) \Gamma(1/2)} \int_{-1}^1
e^{i\rho r} (1-r^2)^{\lambda-1/2} dr,
$$
and $\nu(k)=n/2+k-1$.
This expression was obtained by Walther \cite{Wa2}, and it
can be used to analyse estimates for radially symmetric
equations by carefully looking at the asymptotic
behaviour of Bessel functions and subsequent integrals.

\bigskip
The estimates listed in Theorems 
\ref{prop:basic} and \ref{prop:basic2}
will act as model ones later.
In the subsequent sections, further smoothing results will
be derived from them, hence from simple estimates \eqref{core}
and \eqref{core2''}, by the 
(introduced further) method of canonical transformations
or some combination use of it and the comparison principle.
The following are straightforward results of Theorems \ref{prop:basic}
and \ref{prop:basic2}:
\medskip
\begin{cor}\label{Th:typeI}
Suppose $n\geq1$, $m>0$, and $s>1/2$.
Then we have
\begin{equation}\label{model:1}
\n{\jp{x_n}^{-s}|D_n|^{(m-1)/2}e^{it|D_n|^m}\varphi(x)}_{L^2(\R_t\times\R^n_x)}
\leq
 C\n{\varphi}_{L^2(\R_x^n)}.
\end{equation}
Suppose $n\geq2$, $m>0$, and $s>1/2$.
Then we have
\begin{equation}\label{model:2}
\n{\jp{x_1}^{-s}|D_n|^{(m-1)/2}e^{itD_1|D_n|^{m-1}}\varphi(x)}
_{L^2(\R_t\times\R^n_x)}
\leq
 C\n{\varphi}_{L^2(\R_x^n)}.
\end{equation}
\end{cor}
\medskip
\begin{proof}
Use first the square integrability of $\jp{x_n}^{-s}$ in one dimension,
then apply estimate \eqref{basic:1dim}
in $x_n$ to obtain estimate \eqref{model:1}.
Similarly estimate \eqref{model:2} is obtained from
estimate \eqref{basic:2dim}.
\end{proof}
\medskip
\begin{cor}\label{Th:typeII}
Suppose $m>0$ and $(m-n)/2<\alpha<(m-1)/2$.
Then we have
\begin{equation}\label{model:5}
\n{\abs{x}^{\alpha-m/2}|D_x|^{\alpha}e^{it|D_x|^{m}}
\varphi(x)}_
{L^2\p{\R_t\times\R^n_x}}
\leq C\n{\varphi}_{L^2\p{\R^n_x}}.
\end{equation}
Suppose $m>0$ and $(m-n+1)/2<\alpha<(m-1)/2$.
Then we have
\begin{equation}\label{model:6}
\n{\abs{x}^{\alpha-m/2}|D'|^{\alpha} 
  e^{it\p{|D_1|^m-|D'|^m}}\varphi(x)
}_\L2tx\leq C\n{\varphi}_\Lx,
\end{equation}
where $D'=(D_2,\ldots,D_n)$.
\end{cor}
\medskip
\begin{proof}
Estimate \eqref{model:5} is the same one as estimate \eqref{core2g}.
From estimate \eqref{model:5} in $x'\in\R^{n-1}$,
where $x'=(x_2,\ldots,x_n)$,
and Plancherel's theorem in $x_1$,
we obtain estimate \eqref{model:6}
if we notice the trivial inequality
$|x|^{\alpha-m/2}\leq|x'|^{\alpha-m/2}$.
\end{proof}
\bigskip
By using the comparison principle again,
we can compare many estimates
with the model estimates stated above,
which have been also induced by the comparison principle
from the trivial estimate \eqref{core} and so on.
For example, in notation of Corollary \ref{prop:dim1eq}, setting
$\tau(\xi)=|\xi|^{(m-1)/2}$ and
$g(\xi)=|\xi|^{m}$, we have
$|\tau(\xi)|/
|g'(\xi)|^{1/2}=m^{-1/2}.$
Similarly in notation of Corollary \ref{prop:dim2eq}, setting
$\tau(\xi,\eta)=|\eta|^{(m-1)/2}$ and
$g(\xi,\eta)=\xi|\eta|^{m-1}$, we have
$|\tau(\xi,\eta)|/|\partial g/\partial \xi (\xi,\eta)|^{1/2}=1$.
Hence, noticing that $\chi(D_x)$ is 
$L^2$--bounded for $\chi\in L^\infty$,
we obtain the following {\it secondary comparison}
results from Corollary \ref{Th:typeI}.
\medskip
\begin{cor}\label{COR:dim1ex}
Suppose $n\geq 1$ and $s>1/2$.
Let $\chi\in L^\infty(\R)$. 
Let $f\in C^1(\R)$ be real-valued and strictly monotone
on $\supp\chi$.
Let $\sigma\in C^0(\R)$ be such that, for some $A>0$, we have
\[
|\sigma(\xi)|
\leq A \left|f'(\xi)\right|^{1/2}
\]
for all $\xi\in\supp\chi$.
Then we have
\[
\n{\jp{x_j}^{-s}\chi(D_j)\sigma(D_j)
e^{it f(D_j)}\varphi(x)}_{L^2(\R_t\times\R_x^n)}\leq C
\n{\varphi}_\Lx,
\]
where $j=1,2,\ldots,n$.
\end{cor}
\medskip
\begin{cor}\label{COR:dim2ex}
Suppose $n\geq 2$ and $s>1/2$.
Let $\chi\in L^\infty(\R^2)$. 
Let $f\in C^1(\R^2)$ be a real-valued function such that,
for almost all $\eta\in\R$,
$f(\xi,\eta)$ is strictly monotone in $\xi$
on $\supp\chi$.
Let $\sigma\in C^0(\R^2)$ be such that
for some $A>0$ we have
\[
|\sigma(\xi,\eta)|
\leq A \left|{\frac{\partial f}{\partial \xi}(\xi,\eta)}
\right|^{1/2}
\]
for all $(\xi,\eta)\in\supp\chi$.
Then we have
\[
\n{\jp{x_j}^{-s}\chi(D_j,D_k)\sigma(D_j,D_k)
e^{it f(D_j,D_k)}\varphi(x)}_{L^2(\R_t\times\R_x^n)}\leq C
\n{\varphi}_\Lx,
\]
where $j\neq k$.
\end{cor}
\medskip
\par
Applications of these secondary comparison results
will be given in Section \ref{SECTION:nondisp}.
We can also have results
for radially symmetric operators if we
compare them with estimate \eqref{model:5}
of Corollary \ref{Th:typeII} by using
Theorem \ref{prop:dim1eqmod}.
We will discuss it later in the last part of
Section \ref{SECTION:main} together with more different type of 
results.

\section{Canonical transforms}
\label{SECTION:canonical}
Based on the argument in the introduction,
we will now introduce the main tool to reduce general operators to
normal forms.
That is the canonical transformation which changes
the equation
\[
\left\{
\begin{aligned}
\p{i\partial_t+a(D_x)}\,u(t,x)&=0,\\
u(0,x)&=\varphi(x),
\end{aligned}
\right.
\quad {\rm to}\quad
\left\{
\begin{aligned}
\p{i\partial_t+\sigma(D_x)}\,v(t,x)&=0,\\
v(0,x)&=g(x),
\end{aligned}
\right.
\]
where $a(D_x)$ and $\sigma(D_x)$ are related with each other
as in the relation \eqref{relation} in the introduction,
i.e. we have
$
a(\xi)=\p{\sigma\circ\psi}(\xi).
$
If the initial data $\varphi(x)$ is the corresponding
transform of $g(x)$,
then the solution $u(t,x)=e^{it a(D_x)}\varphi(x)$
is the corresponding transform of $v(t,x)=e^{it \sigma(D_x)}g(x)$.
In this way, we will reduce general smoothing estimates to
model ones listed in Section \ref{SECTION:model}.
\par
\medskip
Now we will describe this more precisely.
Let $\Gamma$, $\widetilde{\Gamma}\subset\R^n$ be open sets and
$\psi:\Gamma\to\widetilde{\Gamma}$
be a $C^\infty$-diffeomorphism (we do not assume them to be cones
since we do not require homogeneity of phases).
We always assume that
\begin{equation}\label{infty}
C^{-1}\leq\abs{\det \partial\psi(\xi)}\leq C\quad(\xi\in\Gamma),
\end{equation}
for some $C>0$.
We set formally
\begin{equation}\label{DefI}
\begin{aligned}
&I_\psi u(x)
=\FT^{-1}\left[\FT u\p{\psi(\xi)}\right](x)
=(2\pi)^{-n}\int_{\R^n}\int_{\R^n}
   e^{i(x\cdot\xi-y\cdot\psi(\xi))}u(y)\, dyd\xi,
\\
&I^{-1}_\psi u(x)
=\FT^{-1}\left[\FT u\p{\psi^{-1}(\xi)}\right](x)
=(2\pi)^{-n}\int_{\R^n}\int_{\R^n}
  e^{i(x\cdot\xi-y\cdot\psi^{-1}(\xi))}u(y)\, dyd\xi.
\end{aligned}
\end{equation}
The operators $I_\psi$ and $I^{-1}_\psi$ can be justified
by using cut-off functions
$\gamma\in C^\infty(\Gamma)$ and
$\widetilde{\gamma}=\gamma\circ\psi^{-1}\in C^\infty(\widetilde{\Gamma})$
which satisfy $\supp\gamma\subset\Gamma$,
$\supp\widetilde{\gamma}\subset\widetilde{\Gamma}$.
We set
\begin{equation}\label{DefI0}
\begin{aligned}
I_{\psi,\gamma} u(x)
&=\FT^{-1}\left[\gamma(\xi)\FT u\p{\psi(\xi)}\right](x)
\\
&=(2\pi)^{-n}\int_{\R^n}\int_{\Gamma}
 e^{i(x\cdot\xi-y\cdot\psi(\xi))}\gamma(\xi)u(y) dyd\xi,
\\
I_{\psi,\gamma}^{-1} u(x)
&=\FT^{-1}\left[\widetilde{\gamma}(\xi)\FT
u\p{\psi^{-1}(\xi)}\right](x)
\\
&=(2\pi)^{-n}\int_{\R^n}
\int_{\widetilde{\Gamma}}
 e^{i(x\cdot\xi-y\cdot\psi^{-1}(\xi))}
 \widetilde{\gamma}(\xi)u(y) dyd\xi.
\end{aligned}
\end{equation}
In the case that $\Gamma$, $\widetilde{\Gamma}\subset\R^n\setminus0$
are open cones,
we may consider the homogeneous $\psi$ and $\gamma$ which satisfy
$\supp\gamma\cap \Sph^{n-1}\subset\Gamma\cap \Sph^{n-1}$ and
$\supp\widetilde{\gamma}\cap 
\Sph^{n-1}\subset\widetilde{\Gamma}\cap \Sph^{n-1}$,
where $\Sph^{n-1}=\b{\xi\in\Rn: |\xi|=1}$.
Then we have the expressions for compositions
\begin{equation}\label{eq:cut}
I_{\psi,\gamma}=\gamma(D_x)\cdot I_\psi=I_\psi\cdot\widetilde{\gamma}(D_x),\quad
I_{\psi,\gamma}^{-1}
=\widetilde{\gamma}(D_x)\cdot I^{-1}_\psi=I^{-1}_\psi\cdot\gamma(D_x),
\end{equation}
and the identities
\begin{equation}\label{eq:id}
I_{\psi,\gamma}\cdot I_{\psi,\gamma}^{-1}=\gamma(D_x)^2,\quad
I_{\psi,\gamma}^{-1}\cdot I_{\psi,\gamma}=\widetilde{\gamma}(D_x)^2.
\end{equation}
We have also the formula
\begin{equation}\label{eq:cnon}
I_{\psi,\gamma}\cdot\sigma(D_x)
=\p{\sigma\circ\psi}(D_x)\cdot I_{\psi,\gamma},\quad
I_{\psi,\gamma}^{-1}\cdot \p{\sigma\circ\psi}(D_x)
=\sigma(D_x)\cdot I_{\psi,\gamma}^{-1}.
\end{equation}
\par
\medskip
We also introduce the weighted $L^2$-spaces.
For the weight function $w(x)$, let $L^2_{w}(\R^n;w)$ be
the set of measurable functions $f:\Rn\to\C$ 
such that the norm
\[
\n{f}_{L^2(\R^n;w)}
=\p{\int_{\R^n}\abs{w(x) f(x)}^2\,dx}^{1/2}
\]
is finite.
Then we have the following fundamental theorem:
\par
\medskip
\begin{thm}\label{Th:reduction}
Assume that the operator $I_{\psi,\gamma}$ defined by \eqref{DefI0}
is $L^2(\R^n;w)$--bounded.
Suppose that we have the estimate
\begin{equation}\label{red}
\n{w(x)\rho(D_x)e^{it\sigma(D_x)}\varphi(x)}_{L^2\p{\R_t\times\R^n_x}}
\leq C\n{\varphi}_{L^2\p{\R^n_x}}
\end{equation}
for all $\varphi$ such that $\supp\widehat{\varphi}\subset\supp\widetilde\gamma$, where
$\widetilde\gamma=\gamma\circ\psi^{-1}$.
Assume also that the function
\begin{equation}\label{bdd}
q(\xi)=\frac{\gamma\cdot\zeta}{\rho\circ \psi}(\xi)
\end{equation}
is bounded.
Then we have
\begin{equation}\label{org}
\n{w(x)\zeta(D_x)e^{ita(D_x)}\varphi(x)}_{L^2\p{\R_t\times\R^n_x}}
\leq C\n{\varphi}_{L^2\p{\R^n_x}}
\end{equation}
for all $\varphi$
such that $\supp\widehat{\varphi}\subset\supp\gamma$,
where $a(\xi)=(\sigma\circ\psi)(\xi)$.
\end{thm}
\medskip
\begin{proof}
Substituting 
$I_{\psi,q}^{-1}\varphi$ for $\varphi$ in \eqref{red},
where $I_{\psi,q}^{-1}=I^{-1}_\psi\cdot q(D_x)$,
we have
\[
\n{w(x)I_{\psi,q}^{-1}(\rho\circ\psi)(D_x)e^{ita(D_x)}\varphi(x)}
  _{L^2\p{\R_t\times\R^n_x}}
\leq C\n{I_{\psi,q}^{-1}\varphi}_{L^2\p{\R^n_x}}
\]
for $\varphi$ such that $\supp\widehat{\varphi}\subset\supp\gamma$.
Here we have noticed \eqref{eq:cnon}.
Then we have
\[
\n{w(x)I_{\psi,\gamma}^{-1}\zeta(D_x)e^{ita(D_x)}\varphi(x)}
  _{L^2\p{\R_t\times\R^n_x}}
\leq C\n{I_{\psi,q}^{-1}\varphi}_{L^2\p{\R^n_x}}.
\]
By Plancherel's theorem, we have the $L^2$--boundedness of
$I_{\psi,q}^{-1}$
if we notice the assumption \eqref{infty} 
and the boundedness of $q(\xi)$
given by \eqref{bdd}.
On the other hand, $I_{\psi,\gamma}$ is $L^2(\R^n;w)$--bounded
by the assumption, and we obtain \eqref{org} if we notice \eqref{eq:id}.
\end{proof}
\medskip
As for the $L^2(\R^n;w)$--boundedness of the operator $I_{\psi,\gamma}$,
we have criteria for some special weight functions.
For $\ka\in\R$, let $L^2_\ka(\R^n)$, $\Dot{L}^2_\ka(\R^n)$ be
the set of measurable functions $f$ such that the norm
\[
\n{f}_{L^2_\ka(\R^n)}
=\p{\int_{\R^n}\abs{\langle x\rangle^\ka f(x)}^2\,dx}^{1/2},
\qquad
\n{f}_{\Dot{L}^2_\ka(\R^n)}
=\p{\int_{\R^n}\abs{|x|^\ka f(x)}^2\,dx}^{1/2}
\]
is finite, respectively.
\par
The following theorem is a simplified version of
\cite[Theorem 1.1]{RS2} given by the authors, where
the $L^2_\ka$--boundedness for more general
$x$-dependent Fourier integral operators was treated under less
restrictive conditions, with exact expressions for
the numbers of derivatives, etc.
These weighted boundedness results played an important role 
in the critical case of some of the smoothing
estimates in \cite{RS3}. They will be of crucial importance here
as well.
\par
\medskip
\begin{thm}\label{Th:L2k}
Suppose $\ka\in\R$.
Assume that all the derivatives of entries of the 
$n\times n$ matrix
$\partial\psi$ and those of $\gamma$ are bounded.
Then the operators $I_{\psi,\gamma}$ and $I^{-1}_{\psi,\gamma}$ 
defined by
\eqref{DefI0} are $L^2_{\ka}(\R^n)$--bounded.
\end{thm}
\medskip
For homogeneous $\psi$ and $\gamma$, 
we have another type of weighted
boundedness result:
\medskip
\begin{thm}\label{Th:L'2k}
Let $\Gamma$, $\widetilde{\Gamma}\subset\R^n\setminus0$ be open cones.
Suppose $|\ka|< n/2$.
Assume $\psi(\lambda\xi)=\lambda\psi(\xi)$,
$\gamma(\lambda\xi)=\gamma(\xi)$ for all $\lambda>0$ and $\xi\in\Gamma$.
Then the operators $I_{\psi,\gamma}$ and $I^{-1}_{\psi,\gamma}$
defined by \eqref{DefI0} are $L^2_{\ka}(\R^n)$--bounded
and $\Dot{L}^2_{\ka}(\R^n)$--bounded.
\end{thm}
\medskip
We remark that the boundedness in Theorem \ref{Th:L'2k}
with the case $\ka\leq0$ is
equivalent to the one with $\ka\geq0$
by the duality argument.
In fact, the formal adjoint of $I_\psi$ can be given by
\begin{align*}
I^*_\psi u(x)
&=(2\pi)^{-n}\int_{\R^n}\int_{\R^n}
   e^{-i(y\cdot\xi-x\cdot\psi(\xi))}u(y)\, dyd\xi,
\\
&=(2\pi)^{-n}\int_{\R^n}\int_{\R^n}
   e^{i(x\cdot\xi-y\cdot\psi^{-1}(\xi))}
\abs{\det \partial\psi^{-1}(\xi)}
u(y)\, dyd\xi,
\\
&=(2\pi)^{-n}\int_{\R^n}\int_{\R^n}
   e^{i(x\cdot\xi-y\cdot\psi^{-1}(\xi))}
\abs{\det \partial\psi(\psi^{-1}(\xi))}^{-1}
u(y)\, dyd\xi,
\\
&=I^{-1}_\psi\cdot\abs{\det \partial\psi(D_x)}^{-1}u(x),
\end{align*}
from which we obtain the formula
\[
I_{\psi,\gamma}^*=I^{-1}_{\psi,d}\,;\qquad
 d(\xi)=\abs{\det\partial\psi(\xi)}^{-1}\gamma(\xi).
\]
Note that $d(\xi)$ satisfies the same property as that of $\gamma(\xi)$
in virtue of \eqref{infty}.
\par
We also remark that the $L^2_\ka(\R^n)$--boundedness in Theorem \ref{Th:L'2k}
is equivalent to the $\Dot{L}^2_\ka(\R^n)$--boundedness.
In fact, the $L^2_{\ka}(\R^n)$--boundedness is a 
straightforward consequence
of the $\Dot{L}^2_{\ka}(\R^n)$--boundedness in the case $\ka\geq0$.
On the other hand, the $L^2_{\ka}(\R^n)$--boundedness induces the
$\Dot{L}^2_{\ka}(\R^n)$--boundedness by the scaling argument
because we have $I_{\psi,\gamma}D_\lambda=D_\lambda I_{\psi,\gamma}$,
and also have
\[
\lambda^{n/2+k}\n{D_\lambda u}_{L^2_\ka(\R^n)}
=\n{(\lambda^2+|x|^2)^{k/2}u(x)}_{L^2(\R^n)}
\to\n{u}_{\Dot{L}^2_\ka(\R^n)}\qquad (\lambda\searrow0),
\]
where $D_\lambda$ denotes the dilation
operator $D_\lambda:u(x)\mapsto u(\lambda x)$.
\medskip
\par
We prepare a few lemmas which will be used to prove Theorem \ref{Th:L'2k}.
The following two results are due to Kurtz and Wheeden 
\cite[Theorem 3]{KW}, and
Stein and Weiss \cite[Theorem B$^{*}$]{SW} 
(see also Lemma \ref{Lem:SW}),
respectively.
\medskip
\begin{lem}\label{Prop:wtbdd}
Suppose $|\ka|<n/2$.
Assume that $m(\xi)\in C^n(\R^n\setminus0)$
and all the derivative of $m(\xi)$ satisfies
$|\partial^\gamma m(\xi)|\leq C_\gamma
|\xi|^{-|\gamma|}$ for all $\xi\not=0$ and $|\gamma|\leq n$.
Then $m(D_x)$ is $L^2_{\ka}(\R^n)$ and 
$\Dot{L}^2_{\ka}(\R^n)$--bounded.
\end{lem}
\medskip
\begin{lem}\label{Prop:sw}
Suppose $1-n/2<\ka<n/2$.
Then the operator
$|D_x|^{-1}$ is $L^2_\ka(\R^n)$-$L^2_{\ka-1}(\R^n)$--bounded and
$\Dot{L}^2_\ka(\R^n)$-$\Dot{L}^2_{\ka-1}(\R^n)$--bounded.
\end{lem}
\medskip
We remark that, in Lemma \ref{Prop:wtbdd},
the $L^2_{\ka}(\R^n)$--boundedness 
is equivalent to the $\Dot{L}^2_{\ka}(\R^n)$--boundedness,
and the $L^2_\ka(\R^n)$-$L^2_{\ka-1}(\R^n)$--boundedness 
in Lemma \ref{Prop:sw}
is also equivalent to
the $\Dot{L}^2_\ka(\R^n)$-$\Dot{L}^2_{\ka-1}(\R^n)$--boundedness,
by essentially the same argument as in the above remark.
\medskip
\begin{proof}[Proof of Theorem \ref{Th:L'2k}]
In view of the remarks below Theorem \ref{Th:L'2k},
it suffices to show the $L_\ka^2$--boundedness of $I_{\psi,\gamma}$
in the case $0\leq \ka<n/2$.
\par
First we assume $n\geq3$.
If we note
\[
e^{ix\cdot\xi}=\frac{1-ix\cdot\partial_\xi}{\jp{x}^2}e^{ix\cdot\xi},
\]
we can justify, by integration by parts,
\begin{align*}
I_{\psi,\gamma} u(x)
&=(2\pi)^{-n}\int\int e^{i(x\cdot\xi-y\cdot\psi(\xi))}
\gamma(\xi) u(y) dyd\xi
\\
&=(2\pi)^{-n}\int\int e^{i(x\cdot\xi-y\cdot\psi(\xi))}
\p{
\frac
{\gamma(\xi)+x\gamma(\xi){}^t\partial\psi(\xi){}^ty+ix\cdot\partial\gamma(\xi)}
{\jp{x}^2}
}
u(y) dyd\xi,
\end{align*}
and have the formula
\begin{equation}\label{eq:formula}
I_{\psi,\gamma}
=\frac1{\jp{x}^2}I_{\psi,\gamma}
+\frac{x}{\jp{x}^2}{}^t\partial\psi(D_x)I_{\psi,\gamma} {}^tx
+i\frac{x}{\jp{x}^2}\cdot I_{\psi,\eta}|D_x|^{-1},
\end{equation}
where $\eta(\xi)=|\psi(\xi)|\partial\gamma(\xi)$, and it satisfies
the same assumption of the theorem as that of $\gamma(\xi)$.
Assume that $I_{\psi,\gamma}$ is $L^2_{\ka-1}$--bounded
under the assumption of the theorem.
Then, by the formula \eqref{eq:formula} and Lemmas \ref{Prop:wtbdd}
and Lemma \ref{Prop:sw}, $I_{\psi,\gamma}$ is also 
$L^2_{\ka}$--bounded if $1-n/2<\ka<n/2$.
On the other hand, by Plancherel's theorem and assumption \eqref{infty},
we have the $L^2$--boundedness of $I_{\psi,\gamma}$
under the assumption of the theorem.
Then, by induction and the interpolation,
we have the $L^2_{\ka}$--boundedness of $I_{\psi,\gamma}$
with $0\leq \ka\leq k_0$, where $k_0$ is the largest integer 
less than $n/2$.
As for $k_0<\ka<n/2$, we have $0< \ka-1< k_0$ in the case $n\geq3$.
Hence, from the $L^2_{\ka-1}$--boundedness of 
$I_{\psi,\gamma}$, we obtain
the $L^2_{\ka}$--boundedness.
\par
In the cases $n=1,2$, we can construct a ($C^1$-)diffeomorphism
$\psi_e:\R^n\setminus0\to\R^n\setminus0$ which is an extension
of $\psi:\Gamma\to\widetilde{\Gamma}$ satisfying
$C^{-1}\leq\abs{\det \partial\psi_e(\xi)}\leq C$
($\xi\in\R^n\setminus0$) for some $C>0$.
(In fact, it is trivial in the case $n=1$.
In the case $n=2$, because of the homogeneity of $\psi(\xi)$,
we have only to extend the function on the arc $\Gamma\cap \Sph^1$
to $\Sph^1$ keeping the diffeomorphism.
It can be carried out by an elementary 
argument and we will omit the
details.)
Then, instead of \eqref{eq:formula},
we have
\[
I_{\psi,\gamma}=\gamma(D_x)I_{\psi_e},
\qquad
I_{\psi_e}
=\frac1{\jp{x}^2}I_{\psi_e}
+\frac{x}{\jp{x}^2}{}^t\partial\psi_e(D_x)I_{\psi_e}{}^tx.
\]
From this formula, together with the $L^2$--boundedness of $I_{\psi_e}$
and that of all the entries of $\partial\psi_e(D_x)$,
we obtain similarly the $L^2_\ka$--boundedness of
$I_{\psi_e}$ with $0\leq \ka\leq1$.
Since we have the $L^2_\ka$--boundedness of $\gamma(D_x)$
for $|\ka|<n/2$ by Lemma \ref{Prop:wtbdd},
we can conclude that $I_{\psi,\gamma}$ 
is $L^2_\ka$--bounded with $0\leq \ka<n/2$.
\end{proof}

\section{Smoothing estimates for dispersive equations}
\label{SECTION:main}
As an application of the canonical transformations described in Section
\ref{SECTION:canonical}, we can derive smoothing estimates for general
dispersive equations from model estimates listed
in Section \ref{SECTION:model}.
Note that the estimates that
we will present are derived from just
two simple estimates \eqref{core} and \eqref{core2''} in virtue of
the comparison principle.
The results which will be thus obtained in this section
generalise many known results of the form
\eqref{eq:eq4} in the introduction.
For the optimality of orders, see Section \ref{SECTION:invariant}.
%
%
\par
Let us consider the solution
\[
u(t,x)=e^{ita(D_x)}\varphi(x)
\]
to the equation
\[
\left\{
\begin{aligned}
\p{i\partial_t+a(D_x)}\,u(t,x)&=0\quad\text{in $\R_t\times\R^n_x$},\\
u(0,x)&=\varphi(x)\quad\text{in $\R^n_x$},
\end{aligned}
\right.
\]
where we always assume that 
function $a(\xi)$ is real-valued.
Let $a_m(\xi)\in C^\infty(\R^n\setminus0)$,
the {\it principal} part of $a(\xi)$, be 
a positively homogeneous function
of order $m$, that is, satisfy
$a_m(\lambda\xi)=\lambda^m a_m(\xi)$ for all $\lambda>0$ and $\xi\neq0$.
\par
We sometimes decompose the initial data $\varphi$ into the sum of
the {\it low frequency} part
$\varphi_{l}$ and the {\it high frequency} part $\varphi_{h}$, where
$\supp\widehat{\varphi_l}\subset\b{\xi:|\xi|< 2R}$
and $\supp\widehat{\varphi_h}\subset\b{\xi:|\xi|> R}$
with sufficiently large $R>0$.
Each part can be realised by multiplying $\chi(D_x)$ or $(1-\chi)(D_x)$
to $\varphi(x)$, hence to $u(t,x)$,
where $\chi\in C_0^\infty(\Rn)$ is an appropriate cut-off function.
\par
First we consider the case that $a(\xi)$
has no lower order terms,
and assume that $a(\xi)$ is {\it dispersive}:
\medskip
\begin{equation}\tag{{\bf H}}
a(\xi)=a_m(\xi),\qquad\nabla a_m(\xi)\neq0 \quad(\xi\in\R^n\setminus0),
\end{equation}
\medskip
\par\noindent
where $\nabla=(\partial_1,\ldots,\partial_n)$ and
$\partial_j=\partial_{\xi_j}$.
A typical example is
$a(\xi)=a_m(\xi)=|\xi|^m$.
Especially, $a(\xi)=a_2(\xi)=|\xi|^2$
is the case of the Schr\"odinger equation.
\par
The following result is derived from Corollary \ref{Th:typeI}
and it is a generalisation of the result by
Ben-Artzi and Klainerman \cite{BK} which treated the case
$a(\xi)=|\xi|^2$ and $n\geq 3$ (using spectral methods):
\medskip
\begin{thm}\label{M:H1}
Assume {\rm{(H)}}.
Suppose $n\geq 1$, $m>0$, and $s>1/2$.
Then we have
\begin{equation}\label{EQ:main1}
\n{\jp{x}^{-s}|D_x|^{(m-1)/2}e^{ita(D_x)}\varphi(x)}_{L^2\p{\R_t\times\R^n_x}}
\leq C\n{\varphi}_{L^2\p{\R^n_x}}.
\end{equation}
\end{thm}
\medskip
Chihara \cite{Ch} proved Theorem \ref{M:H1} in the case $m>1$,
by proving the restriction theorem \eqref{eq:eq5}
or the resolvent estimates \eqref{eq:eq6}.
We will, however, give a simpler proof by reducing estimate
\eqref{EQ:main1}
for elliptic $a(\xi)$
to one dimensional model estimate \eqref{model:1}
and non-elliptic $a(\xi)$ to two dimensional \eqref{model:2}
in Corollary \ref{Th:typeI}.
Recall that these model estimates are a corollary of
estimates \eqref{basic:1dim} and \eqref{basic:2dim}
in Theorem \ref{prop:basic},
which is a direct consequence of just a trivial
estimate \eqref{core}. We also note that $m=1$ is
the case of the wave equation and is important for reducing
the estimates to the model energy conservation case
\eqref{core}.
\par
We also get a scaling invariant estimate for homogeneous weights
$|x|^{-s}$ instead of non-homogenous ones $\jp{x}^{-s}$.
The following result is derived from Corollary \ref{Th:typeII}
and it is a generalisation of the result by
Kato and Yajima \cite{KY} which treated the case $a(\xi)=|\xi|^2$
with $n\geq 3$ and $0\leq \alpha<1/2$, or with $n=2$ and $0<\alpha<1/2$.
Ben-Artzi and Klainerman \cite{BK} gave an alternative
proof of the case $a(\xi)=|\xi|^2$ with $n\geq 3$ and
$0\leq \alpha<1/2$, based on the estimate with 
a non-homogeneous weight and spectral decompositions. 
Our extension of these results is
as follows:
\medskip
\begin{thm}\label{M:H2}
Assume {\rm{(H)}}.
Suppose $m>0$ and $(m-n+1)/2<\alpha<(m-1)/2$,
or $m>0$ and $(m-n)/2<\alpha<(m-1)/2$ in the elliptic case $a(\xi)\neq0$
$(\xi\neq0)$.
Then we have
\begin{equation}\label{EQ:main5}
\n{\abs{x}^{\alpha-m/2}|D_x|^{\alpha}e^{ita(D_x)}\varphi(x)}_
{L^2\p{\R_t\times\R^n_x}}
\leq C\n{\varphi}_{L^2\p{\R^n_x}}.
\end{equation}
\end{thm}
\medskip
Sugimoto \cite{Su1}
proved Theorem \ref{M:H2} for elliptic $a(\xi)$ 
of order $m=2$ and
$1-n/2<\alpha<1/2$, $n\geq 2$.
We note that in general we can not allow 
$\alpha=(m-1)/2$ in estimate \eqref{EQ:main5}, see
Section \ref{SECTION:invariant}. However, a sharp version of
this estimate is still possible if one cut-off the main
global singularity of the solution $u(t,x)=e^{ita(D_x)}\varphi(x)$.
The location of this 
singularity is at the set of all classical trajectories
corresponding to the operators $a(D_x)$. Such results
and their sharpness have been discussed in authors'
paper \cite{RS3}. We note that this case has deep
implications clarifying the null-form structure
for derivative nonlinear Schr\"odinger equations and
equations of similar type.
\par
We have another type of smoothing estimate replacing
$|D_x|^{(m-1)/2}$ by $\jp{D_x}^{(m-1)/2}$.
The following result is a direct consequence of 
Theorems \ref{M:H1} and
\ref{M:H2},
and it also extends the result by Kato and Yajima \cite{KY}
which treated the case $a(\xi)=|\xi|^2$ and $n\geq 3$:
\medskip
\begin{cor}\label{M:H3}
Assume {\rm (H)}.
Suppose $n-1>m>1$, or $n>m>1$ in the elliptic case $a(\xi)\neq0$
$(\xi\neq0).$
Then we have
\begin{equation}\label{EQ:main3}
\n{\jp{x}^{-m/2}\jp{D_x}^{(m-1)/2}e^{ita(D_x)}\varphi(x)}_
{L^2\p{\R_t\times\R^n_x}}
\leq C\n{\varphi}_{L^2\p{\R^n_x}}.
\end{equation}
\end{cor}
\medskip
\begin{proof}[Proof of Corollary \ref{M:H3}]
Theorem \ref{M:H1} implies the stronger estimate
for the high frequency part of estimate \eqref{EQ:main3}
replacing the weight $\jp{x}^{-m/2}$ by $\jp{x}^{-s}$ with $s>1/2$.
Theorem \ref{M:H2} with $\alpha=0$ also implies the stronger estimate
for the low frequency part replacing the weight
$\jp{x}^{-m/2}$ by $|x|^{-m/2}$.
\end{proof}
\medskip
We remark that Walther \cite{Wa2} used spherical harmonics
and asymptotics of Bessel functions
to prove the result of Corollary \ref{M:H3} directly
in the radially symmetric case of
$a(\xi)=|\xi|^m$ (this satisfies assumption (H)
and the ellipticity).
In the elliptic case with $m=2$, Walther's result was extended
to the non-radially symmetric case by the authors \cite{RS2}.
Corollary \ref{M:H3} is the development of that analysis allowing
non-elliptic operators as well.
We may also look at the other type of global smoothing of the
form \eqref{EQ:main3}, but with the weight $\jp{x}^{-m/2}$ replaced
by homogeneous ones.
However, this follows from the previous types.
For example, we can observe that
estimate \eqref{EQ:main5} trivially implies
$$
\n{|x|^{\alpha-m/2}\jp{D_x}^{\alpha}e^{ita(D_x)}\varphi_h(x)}_
{L^2\p{\R_t\times\R^n_x}}
\leq C\n{\varphi_h}_{L^2\p{\R^n_x}},
$$
for high frequency parts,
while for low frequency part we get 
$$
\n{|x|^{-m/2}
e^{ita(D_x)}\varphi_l(x)}_
{L^2\p{\R_t\times\R^n_x}}
\leq C\n{\varphi_l}_{L^2\p{\R^n_x}}
$$
as a special case of \eqref{EQ:main5} with $\alpha=0$.
\medskip
\par
The main idea to prove Theorems \ref{M:H1} and 
\ref{M:H2} is to reduce them
to Corollaries \ref{Th:typeI} and \ref{Th:typeII}
by using Theorem \ref{Th:reduction}.
If some estimate for $e^{it\sigma(D_x)}$ is listed there,
then all our task is to find $\psi(\xi)$
such that $a(\xi)=(\sigma\circ\psi)(\xi)$
and verify all the boundedness assumptions we need.
We will use the notation $\xi=(\xi_1,\ldots,\xi_n)$,
$\eta=(\eta_1,\ldots,\eta_n)$, and $D_x=(D_1,\ldots,D_n)$ 
as used there.
\par
We assume (H).
Let $\Gamma\subset\R^n\setminus0$ 
be a sufficiently small conic neighbourhood of
$e_n=(0,\ldots0,1)$, and take a cut-off function
$\gamma(\xi)\in C^\infty(\Gamma)$ which is 
positively homogeneous of
order $0$ and satisfies $\supp\gamma\cap \Sph^{n-1} 
\subset\Gamma\cap \Sph^{n-1}$.
By the microlocalisation and the rotation of the initial data
$\varphi$,
we may assume 
$\supp\widehat{\varphi}\subset\supp\gamma$.
The dispersive assumption $\nabla a_m(e_n)\neq0$ in this direction
implies the following two possibilities:
\medskip
\begin{description}
\item[(i)]
$\partial_n a_m(e_n)\neq0$.
Then, by Euler's identity $a_m(\xi)=(1/m)\nabla a_m(\xi)\cdot\xi$,
we have $a_m(e_n)\neq0$.
Hence, in this case, we may assume that
$a(\xi)(>0)$ and $\partial_n a(\xi)$ are bounded away from $0$
for $\xi\in\Gamma$.
\item[(ii)]
$\partial_n a_m(e_n)=0$.
Then there exits $j\neq n$ such that
$\partial_j a_m(e_n)\neq0$, say $\partial_1a_m(e_n)\neq0$.
Hence, in this case, we may assume
$\partial_1a(\xi)$ is bounded away from $0$
for $\xi\in\Gamma$.
We remark $a(e_n)=0$ by Euler's identity.
\end{description}
\medskip
\par
\begin{proof}[Proof of Theorem \ref{M:H1}]
The estimate with the case $n=1$ is given by estimate \eqref{model:1}
in Corollary \ref{Th:typeI}.
In fact, we have $a(\xi)=a(1)|\xi|^m$ for $\xi>0$ in this case.
Hence we may assume $n\geq2$.
We remark that it is sufficient to show theorem with $1/2<s<n/2$
because the case $s\geq n/2$ is easily reduced to this case.
\par
In the case (i), we take
\begin{equation}\label{EQ:proofc1a1}
\sigma(\eta)=|\eta_n|^m,\quad
\psi(\xi)=(\xi_1,\ldots,\xi_{n-1},a(\xi)^{1/m}).
\end{equation}
Then we have $a(\xi)=\p{\sigma\circ\psi}(\xi)$ and
\begin{equation}\label{EQ:proofc1a2}
\det\partial\psi(\xi)
=
\begin{vmatrix}
E_{n-1}&0
\\
*&  (1/m)a(\xi)^{1/m-1} \partial_n a(\xi)
\end{vmatrix},
\end{equation}
where $E_{n-1}$ is the identity matrix of order $n-1$.
We remark that \eqref{infty} is satisfied since
$\det\partial\psi(e_n)=(1/m)a(e_n)^{1/m-1} \partial_n a(e_n)\neq0$.
By estimate \eqref{model:1} in Corollary \ref{Th:typeI},
we have
estimate \eqref{red} in Theorem \ref{Th:reduction} with $\sigma(D_x)=|D_n|^m$,
$w(x)=\jp{x}^{-s}$, and $\rho(\xi)=|\xi_n|^{(m-1)/2}$.
Note here the trivial inequality $\jp{x}^{-s}\leq\jp{x_n}^{-s}$.
If we take $\zeta(\xi)=|\xi|^{(m-1)/2}$,
then $q(\xi)=\gamma(\xi)\p{|\xi|/a(\xi)^{1/m}}^{(m-1)/2}$ defined by
\eqref{bdd} is a bounded function.
On the other hand, $I_{\psi,\gamma}$ is $L^2_{-s}$--bounded for $1/2<s<n/2$
by Theorem \ref{Th:L'2k}.
Hence, by Theorem \ref{Th:reduction}, we have estimate \eqref{org},
that is, estimate \eqref{EQ:main1}.
\par
In the case (ii), we take
\[
\sigma(\eta)=\eta_1|\eta_n|^{m-1},\quad
\psi(\xi)=\p{a(\xi)|\xi_n|^{1-m},\xi_2,\ldots,\xi_n}
\]
Then we have $a(\xi)=\p{\sigma\circ\psi}(\xi)$ and
\[
\det\partial\psi(\xi)
=
\begin{vmatrix}
\partial_1 a(\xi)|\xi_n|^{1-m}&*
\\
0&E_{n-1}
\end{vmatrix}.
\]
Since $\det\partial\psi(e_n)=\partial_1 a(e_n)\neq0$,
\eqref{infty} is satisfied.
Similarly to the case (i), the estimate for $\sigma(D_x)=D_1|D_n|^{m-1}$
is given by estimate \eqref{model:2} in Corollary \ref{Th:typeI}, which implies
estimate \eqref{EQ:main1} again by Theorem \ref{Th:reduction}.
\end{proof}
\medskip
\begin{proof}[Proof of Theorem \ref{M:H2}]
In the case (i), which is the only possibility for the elliptic
$a(\xi)\neq0$ ($\xi\neq0$),
we take
\[
\sigma(\eta)=|\eta|^m,\quad
\psi(\xi)
=\p{\xi_1,\ldots,\xi_{n-1},\sqrt{a(\xi)^{2/m}-(\xi_1^2+\cdots+
\xi_{n-1}^2)}}.
\]
Then we have $a(\xi)=\p{\sigma\circ\psi}(\xi)$ and
\[
\det\partial\psi(\xi)
=
\begin{vmatrix}
E_{n-1}&0
\\
*&  (1/m)a(\xi)^{2/m-1}\partial_n a(\xi)
   /\sqrt{a(\xi)^{2/m}-(\xi_1^2+\cdots+\xi_{n-1}^2)}
\end{vmatrix}.
\]
Since $\det\partial\psi(e_n)=(1/m)a(e_n)^{1/m-1}\partial_n a(e_n)\neq0$,
\eqref{infty} is satisfied.
The estimate for $\sigma(D_x)=|D|^m$ is given
by estimate \eqref{model:5} in Corollary \ref{Th:typeII}.
In the case (ii), we take
\[
\sigma(\eta)=|\eta_1|^m-(\eta_2^2+\cdots\eta_n^2)^{m/2},\quad
\psi(\xi)
=\p{\p{a(\xi)+(\xi_2^2+\cdots+\xi_n^2)^{m/2}}^{1/m},\xi_2,\ldots,\xi_n}
\]
Then we have $a(\xi)=\p{\sigma\circ\psi}(\xi)$ and
\[
\det\partial\psi(\xi)
=
\begin{vmatrix}
(1/m)\p{a(\xi)+(\xi_2^2+\cdots+\xi_n^2)^{m/2}}^{1/m-1}
  \partial_1 a(\xi)&*
\\
0&E_{n-1}
\end{vmatrix}.
\]
Since $\det\partial\psi(e_n)=(1/m)\partial_1 a(e_n)\neq0$,
\eqref{infty} is satisfied.
The estimate for 
$\sigma(D_x)=|D_1|^m-(D_2^2+\cdots+D_n^2)^{m/2}$ is given by
estimate \eqref{model:6} in Corollary \ref{Th:typeII}.
By the same argument as used in the proof of Theorem \ref{M:H1},
we have Theorems \ref{M:H2}.
\end{proof}
\medskip
\par
As another advantage of the new method, we can also consider
the case that $a(\xi)$ has lower order terms,
and assume that $a(\xi)$ is dispersive in the following sense:
\medskip
\begin{equation}\tag{{\bf L}}
\begin{aligned}
&a(\xi)\in C^\infty(\R^n),\qquad \nabla a(\xi)\neq0 \quad(\xi\in\R^n),
\quad\nabla a_m(\xi)\neq0 \quad (\xi\in\R^n\setminus0), 
\\
&|\partial^\alpha\p{a(\xi)-a_m(\xi)}|\leq C_\alpha\abs{\xi}^{m-1-|\alpha|}
\quad\text{for all multi-indices $\alpha$ and all $|\xi|\geq1$}.
\end{aligned}
\end{equation}
\medskip
\par\noindent
We note that $a(\xi)=|\xi|^m$ does not satisfy (L) 
because $\nabla a(\xi)$
vanishes at the origin $\xi=0$, while it satisfies (H).
On the other hand, $a(\xi)=a_3(\xi)+\xi_1$
satisfies (L) with $m=3$, where $a_3(\xi)=\xi_1^3+\xi_2^3+\cdots+\xi_n^3$
and $\xi=(\xi_1,\xi_2,\ldots,\xi_n)$.
As will be discussed soon,
the ability to include the lower order terms and conditions
on them is very important in global problems. In fact, it is
known that low frequencies are often responsible for the
orders of decay of the solutions and their smoothing 
property for large times.
However,
the difference between the principal part and the lower order
terms becomes extinct in the low frequency part,
and one has to look at the properties of
the full symbol. Thus, if we want to have the dispersive
behaviour of the problem we need to look at the dispersiveness
of the full symbol in assumption (L). For large $\xi$ conditions
$\nabla a(\xi)\neq0$ and $\nabla a_m(\xi)\neq0$ are clearly
equivalent, while for small $\xi$ condition
$\nabla a_m(\xi)\neq0$ is not necessary (but it is satisfied
anyway due to the homogeneity of $a_m$).
Thus, condition (L) may be formulated also in the following way
\medskip
\begin{equation}\tag{{\bf L}}
\begin{aligned}
&a(\xi)\in C^\infty(\R^n),\qquad 
|\nabla a(\xi)|\geq C\jp{\xi}^{m-1}\quad(\xi\in\R^n)\quad
\textrm{for some}\; C>0,
\\
&|\partial^\alpha\p{a(\xi)-a_m(\xi)}|\leq C_\alpha\abs{\xi}^{m-1-|\alpha|}
\quad\text{for all multi-indices $\alpha$ and all $|\xi|\geq1$}.
\end{aligned}
\end{equation}
The last line of this assumption simply amounts to saying that
the principal part $a_m$ of $a$ is positively homogeneous
of order $m$ for $|\xi|\geq 1$.
\medskip
\par
The following result is also derived from Corollary \ref{Th:typeI}:
\medskip
\begin{thm}\label{M:L4}
Assume {\rm (L)}.
Suppose $n\geq1$, $m>0$, and $s>1/2$.
Then we have
\begin{equation}\label{EQ:main4}
\n{\jp{x}^{-s}\jp{D_x}^{(m-1)/2}e^{ita(D_x)}\varphi(x)}
_{L^2\p{\R_t\times\R^n_x}}
\leq C\n{\varphi}_{L^2\p{\R^n_x}}.
\end{equation}
\end{thm}
\medskip
Thus, by Theorem \ref{M:L4}, 
we can have better estimate than that in Corollary \ref{M:H3}
even under weaker conditions on $m$ and $n$ if we assume 
(L) instead of (H).
This fact does not contradict to the optimality of Corollary \ref{M:H3}
with the case $a(\xi)=|\xi|^m$ (see the remark below Corollary \ref{M:H3})
because it does not satisfy assumption (L). This does
emphasise once again the importance of the dispersiveness
assumption $\nabla a\not=0.$
\par
Note that the following result is a straightforward consequence
of Theorem \ref{M:L4}
and the $L^2$--boundedness
of $|D_x|^{(m-1)/2}\jp{D_x}^{-(m-1)/2}$
with $m\geq1$, which is an analog of Theorem \ref{M:H1} for $a(D_x)$
with lower order terms (assumption $m\geq 1$ is natural to be
able to talk about lower order terms):
\medskip
\begin{cor}\label{M:L5}
Assume {\rm{(L)}}.
Suppose $n\geq1$, $m\geq1$ and $s>1/2$.
Then we have
\begin{equation}\label{EQ:main9}
\n{\jp{x}^{-s}|D_x|^{(m-1)/2}e^{ita(D_x)}\varphi(x)}_{L^2\p{\R_t\times\R^n_x}}
\leq C\n{\varphi}_{L^2\p{\R^n_x}}.
\end{equation}
\end{cor}
\medskip
\begin{proof}[Proof of Theorem \ref{M:L4}]
We decompose the initial data $\varphi$ into the sum of
the high frequency part and the low frequency part.
For high frequency part,
the same argument as in the proof of Theorem \ref{M:H1} is valid.
(Furthermore, we can use Theorem \ref{Th:L2k} 
instead of Theorem \ref{Th:L'2k}
to assure the boundedness of $I_{\psi,\gamma}$, 
hence we need not assume
$n\geq2$.)
We show how to get the estimates for low frequency part.
Because of the compactness of it, we may assume
$\partial_ja(\xi)\neq0$ with some $j$, say $j=n$,
on a bounded set $\Gamma\subset\R^n$ and $\supp\widehat{\varphi}\subset\Gamma$.
Since we have $a(\xi)+c>0$ on $\Gamma$ with some constant $c>0$
and
\[
\n{\jp{x}^{-s}\jp{D_x}^{(m-1)/2}e^{ita(D_x)}\varphi}_{L^2(\R_t\times\R^n_x)}
=
\n{\jp{x}^{-s}\jp{D_x}^{(m-1)/2}
e^{it\p{a(D_x)+2c}}\varphi}_{L^2(\R_t\times\R^n_x)},
\]
we may assume $a(\xi)\geq c>0$ on $\Gamma$ without loss of generality.
We take a cut-off function
$\gamma(\xi)\in C_0^\infty(\Gamma)$,
and choose $\psi(\xi)$ and $\sigma(\eta)$ in the same way as
\eqref{EQ:proofc1a1}.
Assumption \eqref{infty} is also verified if we notice \eqref{EQ:proofc1a2}.
By estimate \eqref{model:1} in Corollary \ref{Th:typeI},
we have
estimate \eqref{red} in Theorem \ref{Th:reduction} with $\sigma(D_x)=|D_n|^m$,
$w(x)=\jp{x}^{-s}$ ($s>1/2$), and 
$\rho(\xi)=|\xi_n|^{(m-1)/2}$ as in
the proof of Theorem \ref{M:H1}.
If we take $\zeta(\xi)=\jp{\xi}^{(m-1)/2}$,
then $q(\xi)=\gamma(\xi)\p{\jp{\xi}/a(\xi)^{1/m}}^{(m-1)/2}$ 
defined by
\eqref{bdd} is a bounded function.
On the other hand, $I_{\psi,\gamma}$ is $L^2_{-s}$--bounded 
for all $s>1/2$
by Theorem \ref{Th:L2k}.
Hence, by Theorem \ref{Th:reduction}, we have estimate \eqref{org},
that is, estimate \eqref{EQ:main4}.
\end{proof}
\medskip
\par
\medskip
Recall that assumption (L) in 
Theorem \ref{M:L4}
requires the condition
$\nabla a(\xi)\neq0$ ($\xi\in\R^n$) for the full symbol, 
besides the same
one $\nabla a_m(\xi)\neq0$ ($\xi\neq0$) for the principal term.
We will now 
introduce an intermediate assumption between (H) and (L),
and discuss what happens if we do not have
the condition $\nabla a(\xi)\neq0$:
\medskip
\begin{equation}\tag{{\bf HL}}
\begin{aligned}
&a(\xi)=a_m(\xi)+r(\xi),\quad
\nabla a_m(\xi)\neq0 \quad (\xi\in\R^n\setminus0),
\quad r(\xi)\in C^\infty(\R^n)
 \\
&|\partial^\alpha r(\xi)|\leq C\jp{\xi}^{m-1-|\alpha|}
\quad\text{for all multi-indices $\alpha$}.
\end{aligned}
\end{equation}
In view of the proof of Theorem \ref{M:L4}, we see that
Theorems \ref{M:H1}, \ref{M:H2}, and Corollary \ref{M:H3}
remain valid if we
replace assumption (H) by (HL) and functions $\varphi(x)$ in 
the estimates by its 
(sufficiently large) high frequency part $\varphi_h(x)$.
However we cannot control the low frequency part $\varphi_l(x)$,
and so have only the time local estimates on the whole:
\medskip
\begin{thm}\label{Th:HL}
Assume {\rm (HL)}.
Suppose $n\geq1$, $m>0$, $s>1/2$, and $T>0$.
Then we have
\[
\int^T_0\n{\jp{x}^{-s}\jp{D_x}^{(m-1)/2}e^{ia(D_x)}}^2_{L^2(\R^n_x)}\,dt
\leq
 C\n{\varphi}_{L^2(\R^n)}^2,
\]
where $C>0$ is a constant depending on $T>0$.
\end{thm}
\medskip
\begin{proof}[Proof of Theorem \ref{Th:HL}]
We decompose $\varphi$ into the sum of low and high 
frequency parts.
For the high frequency part,
the same arguments as in the proof of Theorems \ref{M:H1} and
\ref{M:L4}
are valid (and furthermore we can have the estimate with $T=\infty$).
The estimate for the low frequency part is trivial.
In fact, if $\supp \FT\varphi\subset\{\xi:|\xi|\leq R\}$,
we have
\begin{align*}
\int^T_0
\n{\jp{x}^{-s}\jp{D_x}^{(m-1)/2}e^{ita(D_x)}\varphi(x)}^2_{L^2(\R^n_x)}\,dt
\leq &\int^T_0\n{\jp{D_x}^{(m-1)/2}e^{ita(D_x)}\varphi(x)}^2_{L^2(\R^n_x)}\,dt
\\
\leq &CT\n{\jp{\xi}^{(m-1)/2}\widehat\varphi(\xi)}^2_{L^2(\R^n)}
\\
\leq &CT\jp{R}^{m-1}\n{\varphi}_{L^2(\R^n)}^2
\end{align*}
by Plancherel's theorem.
\end{proof}
\medskip
\par
We remark that Theorem \ref{M:L4} is the time global version
(that is, the estimate with $T=\infty$) of Theorem \ref{Th:HL},
and the extra assumption $\nabla a(\xi)\neq0$ is needed for that.
Since the assumption $\nabla a(\xi)\neq0$ for large $\xi$ 
is automatically satisfied by assumption (HL),
Theorem \ref{M:L4} means that the condition $\nabla a(\xi)\neq0$
for small $\xi$ assures the time global estimate.
In this sense, the low frequency part have a responsibility for
the time global smoothing. 
\bigskip
\par
Finally, we will state a secondary comparison result,
similarly to Corollaries \ref{COR:dim1ex} and \ref{COR:dim2ex} in Section
\ref{SECTION:model}, especially for radially symmetric operators,
which will play various important roles in later sections
(Sections \ref{SECTION:nondisp} and \ref{SECTION:relative}).
In notation of Theorem \ref{prop:dim1eqmod}, setting
$\tau(\rho)=\rho^{(m-1)/2}$ and
$g(\rho)=\rho^{m}$, we have
$|\tau(\rho)|/|g'(\rho)|^{1/2}=m^{-1/2}$.
If we take
$\tau(\rho)=\rho^{\alpha}$ and
$g(\rho)=\rho^{2}$ instead, we have
$|\tau(\rho)|/|g'(\rho)|^{1/2}=2^{-1/2}\rho^{\alpha-1/2}$.
Then we obtain the following results from Theorem \ref{M:H1}
with $a(\xi)=|\xi|^m$ and Theorem \ref{M:H2} with $a(\xi)=|\xi|^2$,
that is, estimate \eqref{core2'} in Section \ref{SECTION:model}:
\medskip
\begin{cor}\label{COR:RStype}
Suppose $n\geq1$, $s>1/2$, and $1-n/2<\alpha<1/2$.
Let $\chi\in L^\infty(\R_+)$. 
Let $f\in C^1(\R_+)$ be real-valued and
strictly monotone on $\supp \chi$.
Let $\sigma\in C^0(\R_+)$
be such that for some $A>0$ we have
\begin{equation}\label{EQ:RStypeass}
|\sigma(\rho)|\leq A |f^\prime(\rho)|^{1/2}
\end{equation}
for all $\rho\in\supp\chi$.
Then we have
\begin{equation}\label{EQ:RStype1}
\n{\jp{x}^{-s}\chi(|D_x|)\sigma(|D_x|)
e^{itf(|D_x|)}\varphi(x)}_\L2tx
\leq C\n{\varphi}_\Lx, 
\end{equation}
\begin{equation} \label{EQ:RStypeh}
\n{|x|^{\alpha-1} \chi(|D_x|)|D_x|^{\alpha-1/2}\sigma(|D_x|)
e^{itf(|D_x|)}\varphi(x)}_\L2tx
\leq C\n{\varphi}_\Lx.
\end{equation}
\end{cor}
\medskip

\section{Invariant estimates and sharpness}
\label{SECTION:invariant}
Let us now suggest an invariant form of smoothing estimates
for the solution
$u(t,x)=e^{ita(D_x)}\varphi(x)$
to the equation
\begin{equation}\label{equation-inv}
\left\{
\begin{aligned}
\p{i\partial_t+a(D_x)}\,u(t,x)&=0\quad\text{in $\R_t\times\R^n_x$},\\
u(0,x)&=\varphi(x)\quad\text{in $\R^n_x$},
\end{aligned}
\right.
\end{equation}
which remain valid also in some
areas without dispersion $\nabla a(\xi)\neq0$, where standard
smoothing estimates are known to fail.
We will discuss these estimates in this section and then will
establish them in a variety of situations in Section \ref{SECTION:nondisp}.
\par
We can equivalently rewrite estimates in Section
\ref{SECTION:main} in the form
\begin{equation}\label{EQ:inv-form}
\n{w(x)\zeta\p{\abs{\nabla a(D_x)}}
e^{it a(D_x)}\varphi(x)}_{L^2\p{\R_t\times\R^n_x}}
\leq C\n{\varphi}_{L^2\p{\R^n_x}},
\end{equation}
where $w$ is a weight function of the form
$w(x)=|x|^{\delta},\,\jp{x}^{\delta}$
and $\zeta$ is a function on $\R_+$ of the form
$\zeta(\rho)=\rho^\eta,\,\p{1+\rho^2}^{\eta/2}$
with some $\delta,\eta\in\R$.
For example, we can rewrite estimate \eqref{EQ:main1}
of Theorem \ref{M:H1} as well as estimate \eqref{EQ:main9}
of Corollary \ref{M:L5} for the dispersive equations in the form
\begin{equation}\label{EQ:maininv}
\n{\jp{x}^{-s}|\nabla a(D_x)|^{1/2}
e^{it a(D_x)}\varphi(x)}_{L^2\p{\R_t\times\R^n_x}}
\leq C\n{\varphi}_{L^2\p{\R^n_x}}.
\end{equation}
Similarly we can rewrite estimate \eqref{EQ:main5} of Theorem \ref{M:H2}
in the form
\begin{equation}\label{EQ:maininv2}
\n{\abs{x}^{\alpha-m/2}|\nabla a(D_x)|^{\alpha/(m-1)}e^{ita(D_x)}\varphi(x)}_
{L^2\p{\R_t\times\R^n_x}}
\leq C\n{\varphi}_{L^2\p{\R^n_x}}
\qquad(m\neq1),
\end{equation}
and estimate \eqref{EQ:main3} of Corollary \ref{M:H3} ($s=-m/2$)
as well as estimate \eqref{EQ:main4} of Theorem \ref{M:L4} 
in the form
\begin{equation}\label{EQ:maininv3}
\n{\jp{x}^{-s}\jp{\nabla a(D_x)}^{1/2}e^{ita(D_x)}\varphi(x)}_
{L^2\p{\R_t\times\R^n_x}}
\leq C\n{\varphi}_{L^2\p{\R^n_x}}.
\end{equation}
Indeed, under assumption (H) we clearly have 
$|\nabla a(\xi)|\geq c|\xi|^{m-1}$, so the equivalence
between estimate \eqref{EQ:maininv} 
and estimate \eqref{EQ:main1} 
in Theorem \ref{M:H1}
follows from the fact that
$|\nabla a(D_x)|^{1/2}|D_x|^{-(m-1)/2}$ and
$|\nabla a(D_x)|^{-1/2}|D_x|^{(m-1)/2}$ are bounded
in $L^2(\Rn)$.
Under assumption (L) the same argument works for large 
frequencies, while for small frequencies both
$\jp{\xi}^{(m-1)/2}$ and $|\nabla a(\xi)|^{1/2}$ are
bounded away from zero.
Thus we have the equivalence between
estimate \eqref{EQ:maininv} and estimate \eqref{EQ:main9} 
in Corollary \ref{M:L5}.
The same is true for the other equivalences.
We may also look at the other type of global smoothing of the
form \eqref{EQ:maininv3}, but with homogeneous weight functions.
However, this follows from the type \eqref{EQ:maininv2}
as was also explained in the remark below Corollary \ref{M:H3}.
\par
Estimate \eqref{EQ:inv-form}, hence
estimates \eqref{EQ:maininv} -- \eqref{EQ:maininv3} are
invariant under canonical transformations by Theorem \ref{Th:reduction}.
More precisely, we have the following theorem: 
\medskip
\begin{thm}\label{Th:caninv}
Let $\zeta$ be a function on $\R_+$ of the form
$\zeta(\rho)=\rho^\eta$ or $\p{1+\rho^2}^{\eta/2}$
with some $\eta\in\R$.
Assume that the operators $I_{\psi,\gamma}$ and $I_{\psi,\gamma}^{-1}$
defined by \eqref{DefI0} are $L^2(\Rn;w)$--bounded.
Then the following two estimates
\begin{align*}
&\n{w(x)\zeta(|\nabla a(D_x)|)
e^{it a(D_x)}\varphi(x)}_{L^2\p{\R_t\times\R^n_x}}
\leq C\n{\varphi}_{L^2\p{\R^n_x}}
\quad (\supp\widehat\varphi\subset\supp\gamma),
\\
&\n{w(x) \zeta(|\nabla \sigma(D_x)|)
e^{it \sigma(D_x)}\varphi(x)}_{L^2\p{\R_t\times\R^n_x}}
\leq C\n{\varphi}_{L^2\p{\R^n_x}}
\quad (\supp\widehat\varphi\subset\supp\widetilde\gamma)
\end{align*}
are equivalent to each other, where $a=\sigma\circ\psi\in C^1$
on $\supp \gamma$ and $\widetilde\gamma=\gamma\circ\psi^{-1}$.
\end{thm}
\medskip
\begin{proof}
Note that $\nabla a(\xi)=\nabla\sigma(\psi(\xi))D\psi(\xi)$
and
$C|\nabla a(\xi)|\leq|\nabla\sigma(\psi(\xi))|\leq C'|\nabla a(\xi)|$
on $\supp\gamma$ with some $C,C'>0$,
which is assured by the assumption \eqref{infty}.
Then the result is obtained from Theorem \ref{Th:reduction}.
\end{proof}
\medskip
\par
On account of these argument,
we will call estimate \eqref{EQ:inv-form}
an {\it invariant estimate},
and indeed we expect invariant estimates
\eqref{EQ:maininv}, \eqref{EQ:maininv2},
and \eqref{EQ:maininv3} to hold for $s>1/2$, $(m-n)/2<\alpha<(m-1)/2$,
and $s=-m/2$ ($n>m>1$),
respectively in ordinally settings (elliptic case for example),
with $m>0$ being the order of $a(D_x)$.
\par
Let us briefly indicate that invariant
estimate \eqref{EQ:maininv} with $s>1/2$ 
is also a refinement of
another known estimate for non-dispersive equations.
If operator $a(D_x)$ has real-valued 
symbol $a(\xi)\in C^1(\Rn)$ 
which is positively homogeneous 
of order $m\geq 1$ and no dispersiveness assumption 
is made, Hoshiro \cite{Ho1} showed the estimate
\begin{equation}\label{EQ:Hoshiro1}
\n{\jp{x}^{-s}\jp{D_x}^{-s}|a(D)|^{1/2} e^{ita(D_x)}
\va(x)}_\L2tx\leq
C\n{\varphi}_\Lx \quad {\rm (}s>1/2{\rm)}.
\end{equation}
But once we prove \eqref{EQ:maininv} with $s>1/2$,
we can have better estimate
$$\n{\jp{x}^{-s}\jp{D_x}^{-1/2}|a(D)|^{1/2} e^{ita(D_x)}
\va(x)}_\L2tx\leq
C\n{\varphi}_\Lx \quad {\rm (}s>1/2{\rm)}$$
with respect to the number of derivatives.
In fact, using the Euler's identity 
$m a(\xi)=\xi\cdot\nabla a(\xi)$,
we see that this estimate trivially follows from
$$\n{\jp{x}^{-s}\jp{D_x}^{-1/2}|D_x|^{1/2}
|\nabla a(D)|^{1/2} e^{ita(D_x)}\va(x)}_\L2tx\leq
C\n{\varphi}_\Lx \quad {\rm (}s>1/2{\rm)},$$
which in turn follows from \eqref{EQ:maininv} with $s>1/2$
because $|D_x|^{1/2}\jp{D_x}^{-1/2}$ is
$L^2(\Rn)$--bounded. In fact, estimate 
\eqref{EQ:Hoshiro1} holds only because of the homogeneity
of $a$, since in this case by Euler's identity
zeros of $a$ contain
zeros of $\nabla a$. In general, estimate
\eqref{EQ:Hoshiro1} cuts off too much, and therefore
does not reflect the nature of the problem for non-homogeneous
symbols, as \eqref{EQ:maininv} still does.
\par
In terms of invariant estimates,
we can also give another explanation to the reason
why we do not have time global estimate in Theorem \ref{Th:HL}.
The problem is that the symbol of
the smoothing operator $\jp{D_x}^{(m-1)/2}$ does not vanish
where the symbol of $\nabla a(D_x)$ vanishes, as should
be anticipated by the invariant estimate \eqref{EQ:maininv}.
If zeros of $\nabla a(D_x)$ are not taken into account, the
weight should change to the one as in estimate \eqref{EQ:main3}.
\medskip
\par
{\bf Sharpness of smoothing estimates.}
Let us now discuss the scaling invariance and sharpness
properties of estimates 
\eqref{EQ:maininv}--\eqref{EQ:maininv3}
taking liberty of also
referring to results that will be established
in the sequel.
Let us restrict to the case when 
$a(\xi)\in C^\infty(\R^n\backslash 0)$ is elliptic and
positively homogeneous of order $m>0$.
Then it is easy to see that  estimate \eqref{EQ:maininv2} is scaling invariant
with respect to the natural scaling $u_\lambda(t,x)=u(\lambda^m t,\lambda x)$
to the solution of equation \eqref{equation-inv}.
If $a(\xi)$ is dispersive, that is if
$\nabla a(\xi)\not=0$ ($\xi\not=0$),
estimate \eqref{EQ:maininv2} holds for 
$(m-n)/2<\alpha<(m-1)/2$ by Theorem \ref{M:H2}.
Also, the validity of this
estimate for some value of $\alpha$ implies the validity of
the estimate for smaller $\alpha$'s (see the proof of this
given just before Lemma \ref{Lem:SW}). Thus, the critical case
of this estimate is for the largest value $\alpha=(m-1)/2$.
In the case of the Schr\"odinger equation ($m=2$) this is
the critical case of Kato--Yajima's estimate and it was
shown to fail in the critical case $\alpha=1/2$ by Watanabe
\cite{W} (although quite implicitly). 
\par
We will now give a more direct explicit argument for the
failure of this and other critical estimates.
We note that in the critical case $\alpha=(m-1)/2$
estimate \eqref{EQ:maininv2} (which we will show to fail)
becomes
\begin{equation}\label{EQ:critical1}
\n{\abs{x}^{-1/2}|\nabla a(D_x)|^{1/2}e^{ita(D_x)}\varphi(x)}_
{L^2\p{\R_t\times\R^n_x}}
\leq C\n{\varphi}_{L^2\p{\R^n_x}}.
\end{equation}
Such an estimate would be very useful for the well-posedness
analysis of derivative nonlinear equations or equations
with magnetic potentials since the
recovery of the loss of regularity would be sharp, 
so one wants to
repair it. One way is to locate and then cut-off the main
singularity. This was done by the authors in \cite{RS3} and
is briefly discussed in Section \ref{SECTION:other}. The other
way is to first observe that this estimate is equivalent 
to a weaker estimate
\begin{equation}\label{EQ:critical2}
\n{\jp{x}^{-1/2}|\nabla a(D_x)|^{1/2}e^{ita(D_x)}\varphi(x)}_
{L^2\p{\R_t\times\R^n_x}}
\leq C\n{\varphi}_{L^2\p{\R^n_x}}.
\end{equation}
Indeed, \eqref{EQ:critical2} follows from \eqref{EQ:critical1}
by the trivial inequality $\jp{x}^{-1/2}\leq |x|^{-1/2}$,
while \eqref{EQ:critical1} follows from 
\eqref{EQ:critical2} by the scaling argument 
(similar to the one just before Theorem \ref{prop:basic2}).
Now, for dispersive $a(\xi)$ by using the canonical transform
method of Section \ref{SECTION:main}, 
estimate \eqref{EQ:critical2} is equivalent to its
normal form. For example, in the case of elliptic $a(D_x)$,
it is equivalent to the one dimensional estimate
\begin{equation}\label{EQ:critical3}
\n{\jp{x}^{-1/2}|D_x|^{(m-1)/2}e^{it|D_x|^m}\varphi(x)}_
{L^2\p{\R_t\times\R_x}}
\leq C\n{\varphi}_{L^2\p{\R_x}}.
\end{equation}
Now, by the comparison principle of Section 
\ref{SECTION:comparison}, it is equivalent to its special case
with $m=1$, which is estimate
\begin{equation}\label{EQ:critical4}
\n{\jp{x}^{-1/2}e^{it|D_x|}\varphi(x)}_
{L^2\p{\R_t\times\R_x}}
\leq C\n{\varphi}_{L^2\p{\R_x}}.
\end{equation}
If $\supp\widehat\va\subset [0,\infty)$, we have
$e^{it|D_x|}\varphi(x)=\va(x+t)$, and so, finally, 
\eqref{EQ:critical4} is equivalent to
\begin{equation}\label{EQ:critical5}
\n{\jp{x}^{-1/2}\varphi(x+t)}_
{L^2\p{\R_t\times\R_x}}
\leq C\n{\varphi}_{L^2\p{\R_x}}.
\end{equation}
The last estimate clearly fails since $\jp{x}^{-1/2}$ is
not in $L^2(\R^1_x)$, thus implying that all the estimates
\eqref{EQ:critical1}--\eqref{EQ:critical5} fail. 
Note that we may talk about equivalence of (false) estimates
here since both the canonical transform method and the
comparison principle apply to expressions on the left 
hand side of these estimates and these arguments are of
equivalence, showing that estimates hold or fail simultaneously.
\par
Now, we can
try to repair \eqref{EQ:critical1}, or rather 
\eqref{EQ:critical2}, by taking a stronger weight 
$\jp{x}^{-s}$ for $s>1/2$. In this way we arrive at the
``almost'' scaling invariant estimate 
\begin{equation}\label{EQ:repair1}
\n{\jp{x}^{-s}|\nabla a(D_x)|^{1/2}e^{ita(D_x)}\varphi(x)}_
{L^2\p{\R_t\times\R^n_x}}
\leq C\n{\varphi}_{L^2\p{\R^n_x}} \quad {\rm(}s>1/2{\rm)},
\end{equation}
which is invariant estimate \eqref{EQ:maininv}.
We remark that, for high frequencies, this estimate implies
another type of invariant estimate
\eqref{EQ:maininv3} in some cases, for example under
assumption (H) of $a(\xi)$ being positively homogeneous of
of order $m\geq 1$.
\par
Let us discuss the third invariant estimate
\eqref{EQ:maininv3} with $s=-m/2$ ($n>m>1$).
For large frequencies
it is weaker than \eqref{EQ:repair1}, so we may
restrict ourselves to bounded frequencies, in which case
\eqref{EQ:maininv3} is equivalent to the estimate
$$
\n{\jp{x}^{-m/2}e^{ita(D_x)}\varphi(x)}_
{L^2\p{\R_t\times\R^n_x}}
\leq C\n{\varphi}_{L^2\p{\R^n_x}}.
$$
By Theorem \ref{prop:basic2} and especially the scaling
argument preceding it, we can conclude that this is in turn
equivalent to the estimate
$$
\n{|x|^{-m/2}e^{ita(D_x)}\varphi(x)}_
{L^2\p{\R_t\times\R^n_x}}
\leq C\n{\varphi}_{L^2\p{\R^n_x}}.
$$
But this estimate is scaling invariant (it is a special case
of \eqref{EQ:maininv2} with $\alpha=0$), which
justifies the sharpness of the order $-m/2$ of the weight.
Thus, the expected orders of the weights in 
invariant estimates \eqref{EQ:maininv}--\eqref{EQ:maininv3}
are sharp.
\par
A similar argument can be used to justify the optimality of the
smoothing operator $|\nabla a(D_x)|^{1/2}$ in estimate
\eqref{EQ:maininv}.
For example, in the case of elliptic $a(D_x)$, the weighted
estimate \eqref{EQ:maininv} for 
$|\nabla a(D_x)|^{1/2} e^{ita(D_x)}\va(x)$
will be reduced (by the canonical transform method)
to the weighted estimate for 
the model case $|D_x|^{(m-1)/2} e^{it|D_x|^m}\va(x)$. 
This, in turn,
by the comparison principle, can
be reduced to the pointwise estimate for its special case
$m=1$, that is, to the $L^2$--estimate for
$e^{it|D_x|}\va(x)=\va(x+t)$, with $\supp\widehat\va\subset
[0,\infty)$. Since there is no smoothing of a travelling wave,
operator $|\nabla a(D_x)|^{1/2}$ in \eqref{EQ:maininv}
is sharp. Similar arguments apply to non-elliptic 
dispersive 
$a(\xi)$ by reducing to models in two dimension, and to
non-homogeneous
symbols $a(\xi)$ 
by using assumption (L) in Section \ref{SECTION:main}.

%
%
%
%
%

\section{Equations with time-dependent coefficients}
\label{SECTION:time-dependent}
We now briefly discuss smoothing estimates for equations
with time-dependent coefficients:
\begin{equation}\label{equation-t}
\left\{
\begin{aligned}
\p{i\partial_t+b(t,D_x)}\,u(t,x)&=0\quad\text{in $\R_t\times\R^n_x$},\\
u(0,x)&=\varphi(x)\quad\text{in $\R^n_x$}.
\end{aligned}
\right.
\end{equation}
If the symbol $b(t,\xi)$ is independent of $t$, invariants estimates
\eqref{EQ:maininv}, \eqref{EQ:maininv2} and \eqref{EQ:maininv3} say that
$\nabla_\xi b(t,D_x)$ is responsible for the smoothing property.
The natural question here is what quantity replaces it 
if $b(t,\xi)$ depends on $t$.
\par
We can give an answer to this question
if $b(t,\xi)$ is of the product type
\[
b(t,\xi)=c(t)a(\xi),
\]
where we only assume that $c(t)>0$ is a continuous function.
In the case of dispersive and Strichartz estimates for
higher order (in time) equations the
situation may be very delicate and in general depends on
the rates of oscillations of $c(t)$ (see e.g. Reissig \cite{Rei},
for the wave equation).
For smoothing estimates, we will be able to state a very
general result in Theorem \ref{Th:time-dependent} below.
The final formulae show that a natural extension of the
invariant estimates of the
previous section still remain valid in this case.
In this special case, the equation \eqref{equation-t} can be 
transformed
to the equation with time-independent coefficients.
In fact, by the assumption for $c(t)$, the function
\[
C(t)=\int^t_0 c(s)\,ds
\] 
is strictly monotone and the inverse $C^{-1}(t)$ exists.
Then the function
\[
v(t,x)=u(C^{-1}(t),x)
\]
satisfies
\[
\partial_tv(t,x)=\frac1{c(C^{-1}(t))}(\partial_tu)(C^{-1}(t),x),
\]
hence $v(t,x)$ solves the equation
\[
\left\{
\begin{aligned}
\p{i\partial_t+a(D_x)}\,v(t,x)=&0,\\
v(0,x)=&\varphi(x),
\end{aligned}
\right.
\]
if $u(t,x)$ is a solution to equation \eqref{equation-t}.
By this argument,
invariant estimates 
for $v(t,x)=e^{ita(D_x)}\varphi(x)$
should imply some standard estimates for the solution
\[
u(t,x)=v(C(t),x)=e^{i\int^t_0b(s,D_x)\,ds}\varphi(x)
\]
to equation \eqref{equation-t}.
For example, if we notice the relations
\[
\n{v(\cdot,x)}_{L^2}=\n{|c(\cdot)|^{1/2}u(\cdot,x)}_{L^2}
\]
and
\[
c(t)\nabla a(D_x)=\nabla_\xi b(t,D_x),
\]
we obtain the estimate
\begin{equation}\label{EQ:maininv-t}
\n{\jp{x}^{-s}|\nabla_\xi b(t,D_x)|^{1/2}
e^{i\int^t_0b(s,D_x)\,ds}\varphi(x)}_{L^2\p{\R_t\times\R^n_x}}
\leq C\n{\varphi}_{L^2\p{\R^n_x}}
\end{equation}
from the invariant estimate \eqref{EQ:maininv}.
Estimate \eqref{EQ:maininv-t} is a natural extension of the invariant
estimate \eqref{EQ:maininv} to the case of time-dependent coefficients,
which says that
$\nabla_\xi b(t,D_x)$ is still
responsible for the smoothing property.
From this point of view, we may call it an invariant estimate
too. We can also note that estimate \eqref{EQ:maininv-t} may
be also obtained directly, by formulating an obvious extension of
the comparison principles to the time dependent setting.
We also have similar estimates from the invariant estimates
\eqref{EQ:maininv2} and \eqref{EQ:maininv3}. The same method
of the proof yields the following:

\medskip
\begin{thm}\label{Th:time-dependent}
Let $[\alpha,\beta]\subset [-\infty,+\infty]$. Assume that
function $c=c(t)$ is continuous on $[\alpha,\beta]$ and that
$c\not=0$ on $(\alpha,\beta)$. Let $u=u(t,x)$ be the solution
of equation \eqref{equation-t} with $b(t,\xi)=c(t)a(\xi)$,
where $a$ satisfies assumptions of any of 
Theorem or Corollary \ref{M:H1}, \ref{M:H2}, \ref{M:H3},
\ref{M:L4}, \ref{M:L5}, or \ref{Th:HL}. Then the smoothing
estimate of the corresponding Theorem or Corollary holds
provided we replace $L^2(\R_t,\Rnx)$ by 
$L^2([\alpha,\beta],\Rnx)$, and insert $|c(t)|^{1/2}$ in the
left hand side norms.
\end{thm}
We note that it is possible that $\alpha=-\infty$ and that
$\beta=+\infty$. By continuity of $c$ at such points we
simply mean that the limits of $c(t)$ exist as
$t\to\alpha+$ and as $t\to\beta-$.

To give an example of an estimate 
from Theorem \ref{Th:time-dependent}, 
let us look at the case of Theorem \ref{M:H1}.
In that theorem, we suppose that $a(\xi)$ satisfies assumption (H),
and we assume $n\geq 1$, $m>0$, and $s>1/2$. 
Theorem \ref{M:H1} assures that in this case we have the
smoothing estimate \eqref{EQ:main1}, which is 
\begin{equation}\label{EQ:main1d}
\n{\jp{x}^{-s}|D_x|^{(m-1)/2}e^{ita(D_x)}\varphi(x)}_
{L^2\p{\R_t\times\R^n_x}}
\leq C\n{\varphi}_{L^2\p{\R^n_x}}.
\end{equation}
Theorem \ref{Th:time-dependent} states that solution
$u(t,x)$ of equation \eqref{equation-t} 
satisfies this estimate
provided we replace $L^2(\R_t,\Rnx)$ by 
$L^2([\alpha,\beta],\Rnx)$, and insert $|c(t)|^{1/2}$ in the
left hand side norm. This means that $u$ satisfies
\begin{equation}\label{EQ:main1-time}
\n{\jp{x}^{-s}|c(t)|^{1/2}|D_x|^{(m-1)/2}u(t,x)}_
{L^2\p{[\alpha,\beta]\times\R^n_x}}
\leq C\n{\varphi}_{L^2\p{\R^n_x}}.
\end{equation}
The same is true with statements of any of Theorem or Corollary
\ref{M:H2}, \ref{M:H3}, \ref{M:L4}, \ref{M:L5}, or \ref{Th:HL}.

\section{Smoothing estimates for non-dispersive equations}
\label{SECTION:nondisp}
The most important application of the secondary comparison results
Corollaries \ref{COR:dim1ex}, \ref{COR:dim2ex}, and \ref{COR:RStype}
which were stated in Sections \ref{SECTION:model} and \ref{SECTION:main} is to
the smoothing estimates for the equations 
\[
\left\{
\begin{aligned}
\p{i\partial_t+a(D_x)}\,u(t,x)&=0\quad\text{in $\R_t\times\R^n_x$},\\
u(0,x)&=\varphi(x)\quad\text{in $\R^n_x$},
\end{aligned}
\right.
\]
where real-valued function $a(\xi)$ fails to
satisfy dispersive assumption (H) or (L) in Section \ref{SECTION:main}.
In Corollary \ref{COR:RStype} for example,
even if we lose the dispersiveness assumption at zeros
of $f^\prime$, the estimate is still valid because 
$\sigma$ must vanish at the same points with the order
determined by condition \eqref{EQ:RStypeass}.
The same is true in other comparison results
Corollaries \ref{COR:dim1ex} and \ref{COR:dim2ex}.
In this section, we will treat
the smoothing estimates of non-dispersive equations 
based on this observation.
\par
The following result states that we still have
estimate \eqref{EQ:maininv} of invariant form 
suggested in Section \ref{SECTION:invariant}
for {\em non-dispersive equations}
in a general setting of the radially symmetric case:
\medskip
\begin{thm}\label{Th:nondisprad}
Suppose $n\geq 1$ and $s>1/2$.
Let $a(\xi)=f(|\xi|)$, where $f\in C^1(\R_+)$ is real-valued.
Assume that $f'$ has only finitely many zeros.
Then we have
\[
\n{\jp{x}^{-s}|\nabla a(D_x)|^{1/2}
e^{it a(D_x)}\varphi(x)}_{L^2\p{\R_t\times\R^n_x}}
\leq C\n{\varphi}_{L^2\p{\R^n_x}}.
\]
\end{thm}
\begin{proof}
Noticing $|\nabla a(\xi)|=|f'(|\xi|)|$, use
Corollary \ref{COR:RStype} for $\sigma(\rho)=|f'(\rho)|^{1/2}$
in each interval where $f$ is strictly monotone.
\end{proof}
\medskip
\begin{example}\label{nondisprad}
As a consequence of Theorem \ref{Th:nondisprad},
we have the estimate of invariant form \eqref{EQ:maininv}
if $a(\xi)$ is a real polynomial of $|\xi|$.
For example, let $a(\xi)=f(|\xi|^2)^2$,
with $f(\rho)$ being a non-constant polynomial on $\R$,
The principal part $a_m(\xi)$ of $a(\xi)$
is a power of $|\xi|^2$ multiplied by
a constant, hence it satisfies $\nabla a_m(\xi)\neq0$ ($\xi\neq0$).
If $f(\rho)$ is a homogeneous polynomial,
then $a(\xi)$ satisfies assumption (H) in Section \ref{SECTION:main}
and we have estimate \eqref{EQ:maininv} by Theorem \ref{M:H1}.
In the case when $f(\rho)$ is not homogeneous, 
trivially $a(\xi)$ does not satisfy (H).
Furthermore $a(\xi)$ does not satisfies assumption (L)
in Section \ref{SECTION:main} either
since $\nabla a(\xi)=4f(|\xi|^2)f'(|\xi|^2)\xi$ vanishes
on the set $|\xi|^2=c$ such that $f(c)=0$ or $f'(c)=0$
as well as at the origin $\xi=0$.
Hence Corollary \ref{M:L5} does not assure the estimate \eqref{EQ:maininv},
but even in this case, we have it by Theorem \ref{Th:nondisprad}.
\end{example}
\medskip
If we use Corollaries \ref{COR:dim1ex} and \ref{COR:dim2ex},
we can obtain estimate \eqref{EQ:maininv}
for non-dispersive equations in the non-radially symmetric case,
as well.
We will not try to exhaust the general case, but give some typical
examples of the case when assumption (L) 
in Section \ref{SECTION:main} breaks.
Below, we use the notation
$x=(x_1,\ldots,x_n)$,
$\xi=(\xi_1,\ldots,\xi_n)$, and $D_x=(D_1,\ldots,D_n)$
as in Section \ref{SECTION:comparison}:
\begin{example}\label{nondispersive1}
Let $a(\xi)=\xi_1^2+\xi_2^2+\xi_1$ in $\R^2$, so
that $a(\xi)$ fails to satisfy
$\nabla a(\xi)\neq0$ at $(\xi_1,\xi_2)=(-1/2,0)$.
By taking $\chi$ in Corollaries \ref{COR:dim1ex} and \ref{COR:dim2ex}
to be characteristic functions of appropriate
sets and $s>1/2$, we obtain estimates
\begin{align*}
&\n{\jp{x_2}^{-s}|D_2|^{1/2}
e^{it D_2^2}\varphi(x)}_{L^2(\R_t\times\R_x^2)}
\leq C\n{\varphi}_{L^2(\R_x^2)}, \\
&\n{\jp{x_1}^{-s}|2D_1+1|^{1/2}
e^{it (D_1^2+D_1)}\varphi(x)}_{L^2(\R_t\times\R_x^2)}
\leq C\n{\varphi}_{L^2(\R_x^2)},
\end{align*}
which imply 
\[
\n{\jp{x}^{-s}\p{|2D_1+1|^{1/2}+|D_2|^{1/2}}
e^{it a(D_x)}\varphi(x)}_{L^2(\R_t\times\R_x^2)}
\leq C\n{\varphi}_{L^2(\R_x^2)}
\]
by trivial inequalities $\jp{x}^{-s}\leq\jp{x_k}^{-s}$ ($k=1,2$)
and Plancherel's theorem for partial Fourier transforms.
Substituting $\eta(D_x)\varphi$
for $\varphi$ in the estimate,
where
\[
\eta(\xi)=
|\nabla a(\xi)|^{1/2}\p{|2\xi_1+1|^{1/2}+|\xi_2|^{1/2}}^{-1},
\]
we have
estimate \eqref{EQ:maininv}
if we note the boundedness of $\eta(\xi)$ 
and use Plancherel's theorem again.
\end{example}
\begin{example}\label{nondispersive2}
Let $a(\xi)=\xi_1^3+\xi_2^3+\xi_3^2$ in $\R^3$, so that
$a(\xi)$ fails to satisfy $\nabla a(\xi)\neq0$
at $(\xi_1,\xi_2,\xi_3)=(0,0,0)$.
Choosing appropriate cut-off functions $\chi$ in
Corollary \ref{COR:dim1ex}, we get the
estimates
\begin{align*}
&\n{\jp{x_k}^{-s}|D_k|
e^{it D_k^3}\varphi(x)}_{L^2_{t,x}}
\leq C\n{\varphi}_{L^2(\R_x^3)}
\qquad(k=1,2),
\\
&\n{\jp{x_3}^{-s}|D_3|^{1/2}
e^{it D_3^2}\varphi(x)}_{L^2_{t,x}}
\leq C\n{\varphi}_{L^2(\R_x^3)}
\end{align*}
for $s>1/2$,
which implies 
\[
\n{\jp{x}^{-s}\p{|D_1|+|D_2|+|D_3|^{1/2}}
e^{it a(D_x)}\varphi(x)}_{L^2(\R_t\times\R_x^3)}
\leq C\n{\varphi}_{L^2(\R_x^3)},
\]
hence we get estimate \eqref{EQ:maininv}
in the same way as in Example \ref{nondispersive1}.
\end{example}
\begin{example} \label{nondispersive-shrira}
Equations of the third order often appear in applications to KdV
and nonlinear Schr\"odinger equation. For example, the Shrira
equation \cite{Sh} describing the propagation of a 
three-dimensional packet of weakly nonlinear internal gravity 
waves leads to third order polynomials in two dimensions.
Strichartz estimates for the corresponding solutions have 
been analysed by e.g. Ghidaglia and Saut \cite{GS} and by
Ben-Artzi, Koch and Saut \cite{BKS} by reducing the equations
to pointwise estimates for operators in normal forms given by
$$
a(\xi)=\xi_1^3+\xi_2^3,\quad a(\xi)=\frac16\xi_1^3+\frac12\xi_2^2
\quad\textrm{and}\quad a(\xi)=\frac12(\xi_1^2+\xi_1\xi_2^2).
$$
By the same argument
of Example \ref{nondispersive2}, we obtain estimate \eqref{EQ:maininv}
for the first two polynomials ($s>1/2$).
For the third polynomial, we use Corollary \ref{COR:dim2ex}
to obtain the estimates
\begin{align*}
&\n{\jp{x_1}^{-s}|D_1+D_2^2/2|^{1/2}
e^{it a(D_1,\,D_2)}\varphi(x)}_{L^2(\R_t\times\R_x^2)}
\leq C\n{\varphi}_{L^2(\R_x^2)},
\\
&\n{\jp{x_2}^{-s}|D_1D_2|^{1/2}
e^{it a(D_1,\,D_2)}\varphi(x)}_{L^2(\R_t\times\R_x^2)}
\leq C\n{\varphi}_{L^2(\R_x^2)}
\end{align*}
for $s>1/2$, which imply estimate \eqref{EQ:maininv}
by the same argument as in Example \ref{nondispersive1}.
\end{example}
\par
\bigskip
Now we will present two more 
approaches to treat non-dispersive equations.
Recall that, in Section \ref{SECTION:main},
the method of canonical transformation effectively works
to reduce smoothing estimates for dispersive equations
to standard estimates.
We explain here that this strategy works for non-dispersive 
cases as well.
\par
We will however look at the rank of the Hessian $\nabla^2a(\xi)$,
instead of the principal type assumption $\nabla a(\xi)\not=0$.
Assume now that $a(\xi)\in C^\infty(\Rn\setminus0)$ is 
real-valued and
positively homogeneous of order two.
It can be noted that 
from Euler's identity we obtain
\begin{equation}\label{id:euler2}
\nabla a(\xi)=\xi\nabla^2 a(\xi)
\end{equation}
since $\nabla a(\xi)$ is homogeneous
of order one (here $\xi$ is viewed as a row).
Then the condition $\rank \nabla^2 a(\xi)=n$ implies
$\nabla a(\xi)\not=0$ ($\xi\neq0$),
and as we have already explained,
we have estimates \eqref{EQ:maininv} and \eqref{EQ:maininv2} 
by Theorems \ref{M:H1} and \ref{M:H2}
in this favourable case.
We will show that in the non-dispersive situation
the rank of $\nabla^2 a(\xi)$ still has a responsibility for
smoothing properties.
We assume 
\begin{equation} \label{a1}
\rank \nabla^2 a(\xi)
\geq k\quad\text{whenever}\quad \nabla a(\xi)=0\quad
(\xi\neq0)
\end{equation}
with some $1\leq k\leq n-1$.
We note that condition \eqref{a1} is invariant under
the canonical transformation
in the following sense:
\medskip
\begin{lem}\label{l1}
Let $a=\sigma\circ \psi$, with $\psi:U\to\R^n$
satisfying $\det D\psi(\xi)\not=0$ on an open set
$U\subset\R^n$.
Then, for each $\xi\in U$,
$\nabla a(\xi)=0$ if and only if $\nabla \sigma(\psi(\xi))=0$.
Furthermore the ranks of $\nabla^2 a(\xi)$ and $\nabla^2\sigma(\psi(\xi))$
are equal on $\Gamma$ whenever $\xi\in U$ and $\nabla a(\xi)=0$.
\end{lem}
\begin{proof}
Differentiation gives
$\nabla a(\xi)=\nabla\sigma(\psi(\xi))D\psi(\xi)$ and
we have the first assertion.
Another differentiation gives
$\nabla^2 a(\xi)=\nabla^2\sigma (\psi(\xi)) D\psi D\psi$
when $\nabla a(\xi)=0$.
This implies the second assertion.
\end{proof}
\medskip
To fix the notation, we assume
\begin{equation}\label{a2}
\nabla a(e_n)=0
\quad\text{and}\quad
\rank\nabla^2 a(e_n)=k
\quad(1\leq k\leq n-1),
\end{equation}
where $e_n=(0,\ldots,0,1)$.
Then we have $a(e_n)=0$ by Euler's identity $2a(\xi)=\xi\cdot\nabla a(\xi)$.
We claim that there exists a conic neighbourhood
$\Gamma\subset\Rn\setminus0$ of $e_n$
and a homogeneous $C^\infty$-diffeomorphism $\psi:\Gamma\to\widetilde\Gamma$
(satisfying $\psi(\lambda\xi)=\lambda\psi(\xi)$ for all $\lambda>0$
and $\xi\in\Gamma$)
as appeared in Section \ref{SECTION:canonical}
such that we have the form
\begin{equation}\label{rf}
a(\xi)=(\sigma\circ\psi)(\xi),
\quad
\sigma(\eta)=
c_1\eta_1^2+\cdots+c_k\eta_k^2+r(\eta_{k+1},\ldots,\eta_n),
\end{equation}
where $\eta=(\eta_1,\ldots,\eta_n)$ and $c_j=\pm1$ ($j=1,2,\ldots,k$).
We remark that $r$ must be real-valued and 
positively homogeneous of order two.
\par
We will prove the existence of such $\psi$ that
will satisfy \eqref{rf}.
By \eqref{id:euler2}, \eqref{a2}, and the symmetricity,
all the entries of the matrix $\nabla^2 a (e_n)$ are zero
except for the (perhaps)
non-zero upper left $(n-1)\times(n-1)$ corner matrix.
Moreover, by a linear transformation involving only
the first $(n-1)$ variables of
$\xi=(\xi_1,\ldots,\xi_{n-1},\xi_n)$,
we may assume
$\partial^2a/\partial{\xi_1}^2(e_n)\not=0$.
We remark that \eqref{a2} still holds
under this transformation.
Then, by the Malgrange preparation theorem, we can write
\begin{equation} \label{ra1}
a(\xi)=\pm c(\xi)^2 (\xi_1^2+a_1(\xi^\prime)\xi_1+a_2(\xi^\prime)),
\quad
\xi^\prime=(\xi_2,\ldots,\xi_n).
\end{equation}
locally  in a neighbourhood of $e_n$,
where $c(\xi)>0$ is some strictly positive function, while
function $a_1$ and $a_2$ are smooth and real valued.
Restricting this expression to the hyperplane $\xi_n=1$, 
and using the homogeneity
\[
a(\xi)=\pm\xi_n^2a(\xi_1/\xi_n,\ldots,\xi_{n-1}/\xi_n,1),
\]
we can extend the expression \eqref{ra1} to a
conic neighbourhood $\Gamma$ of $e_n$, 
so that functions 
$c(\xi), a_1(\xi^\prime)$ and $a_2(\xi^\prime)$ 
are positively homogeneous of orders zero, one, and two, respectively.
Let us define $\psi_0(\xi)=c(\xi)\xi$ and 
$\tau(\eta)=\pm(\eta_1^2+a_1(\eta^\prime)\eta_1+a_2(\eta^\prime))$,
so that $a(\xi)=(\tau\circ\psi_0)(\xi)$,
where we write $\eta=(\eta_1,\eta')$, $\eta'=(\eta_2,\ldots,\eta_n)$.
Furthermore,
let us define
$\psi_1(\xi)=(\xi_1+\frac12 a_1(\xi^\prime),\xi')$,
so that $\tau(\xi)=(\sigma\circ\psi_1)(\xi)$ with
$\sigma(\eta)=\eta_1^2+r(\eta^\prime)$,
where $r(\eta^\prime)=
a_2(\eta^\prime)-\frac14 a_1(\eta^\prime)^2$ is
positively homogeneous of degree two.
Then we have $a=\sigma\circ\psi$, where $\psi=\psi_1\circ\psi_0$,
and thus we have the expression \eqref{rf} with $k=1$.
\par
We note that, by the construction, we have
$\psi(e_n)=c(e_n)(\frac12a_1(e_n'),e_n')$,
where $c(e_n)>0$ and $e_n'=(0,\ldots,0,1)\in\R^{n-1}$.
Then we can see that the function $r(\eta^\prime)$ of
$(n-1)$-variables is
defined on a conic neighbourhood of $e_n'$ in $\R^{n-1}$.
On account of this fact and Lemma \ref{l1},
we can apply the same argument above to $r(\eta^\prime)$,
and repeating the process $k$-times,
we have the expression \eqref{rf}.
\par
To complete the proof,
we check that $\det D\psi_0(\xi)=c(\xi)^n$, which clearly 
implies $\det D\psi(\xi)=c(\xi)^n$, and assures that 
it does not vanish on a sufficiently narrow $\Gamma$.
We observe first that
$D\psi_0(\xi)=c(\xi)I_n+{}^t\xi\nabla c(\xi)$, 
where $I_n$ is the identity $n$ by $n$ matrix. 
We note that if we consider the matrix
$A=(\alpha_i\beta_j)_{i,j=1}={}^t\alpha \beta$, where
$\alpha=(\alpha_1,\ldots,\alpha_n)$, $\beta=(\beta_1,\ldots,\beta_n)$,
then $A$ has rank one, so its eigenvalues are $n-1$ zeros and
some $\lambda$. But $\tr A$ is also the sum of the eigenvalues,
hence $\lambda=\tr A$. 
Now, let $\alpha=\xi$, $\beta=\nabla c(\xi)$, and $A={}^t\alpha \beta$.
Since $c(\xi)$ is homogeneous of order zero, by Euler's identity
we have $\tr A=\xi\cdot\nabla c(\xi)=0$, hence all eigenvalues of
$A$ are zero. It follows now that there is a non-degenerate matrix
$S$ such that $S^{-1} A S$ is strictly upper triangular.
But then $\det D\psi_0(\xi)=\det (c(\xi)I_n + S^{-1}A S)$, where
matrix $c(\xi)I_n + S^{-1}A S$ is upper 
triangular with $n$ copies of
$c(\xi)$ at the diagonal. Hence $\det D\psi_0(\xi)=c(\xi)^n$.
\par
\medskip
On account of the above observation,
we have the following result which states that estimates
\eqref{EQ:maininv} and \eqref{EQ:maininv2} 
with $m=2$ still holds for a
class of non-dispersive equations. It is an illustrations
of invariant estimates \eqref{EQ:main-invariant} and
\eqref{EQ:main-invariant3} with $m=2$.
\medskip
\begin{thm}\label{th1}
Let $a\in C^\infty(\Rn\setminus0)$ be real-valued and satisfy
$a(\lambda\xi)=\lambda^2 a(\xi)$ for all $\lambda>0$ and $\xi\neq0$.
Assume that $\rank\nabla^2a(\xi)\geq n-1$ whenever
$\nabla a(\xi)=0$ and $\xi\neq0$.
Suppose $n\geq2$ and $s>1/2$.
Then we have
\[
\n{\jp{x}^{-s}|\nabla a(D_x)|^{1/2}
e^{it a(D_x)}\varphi(x)}_{L^2\p{\R_t\times\R^n_x}}
\leq C\n{\varphi}_{L^2\p{\R^n_x}}.
\]
Suppose $(4-n)/2<\alpha<1/2$, or $(3-n)/2<\alpha<1/2$
in the elliptic case $a(\xi)\neq0$ {\rm (}$\xi\neq0${\rm)}.
Then we have
\[
\n{\abs{x}^{\alpha-1}|\nabla a(D_x)|^{\alpha}e^{ita(D_x)}\varphi(x)}_
{L^2\p{\R_t\times\R^n_x}}
\leq C\n{\varphi}_{L^2\p{\R^n_x}}.
\]
\end{thm}
\medskip
\begin{proof}
By microlocalisation and an appropriate rotation,
we may assume $\supp \widehat{\varphi}\subset \Gamma$,
where $\Gamma\subset\Rn\setminus0$ is a sufficiently 
narrow conic neighbourhood
of the direction $e_n=(0,\ldots,0,1)$.
Since everything is all right in the dispersive case 
$\nabla a(e_n)\neq0$ by Theorems \ref{M:H1} and \ref{M:H2},
we assume $\nabla a(e_n)=0$.
We may also assume $n\geq2$ since $\nabla a(e_n)=0$ 
implies $\nabla a(\xi)=0$
for all $\xi\neq0$ in the case $n=1$.
Then we have $\rank\nabla^2a(e_n)\neq n$ by the relation \eqref{id:euler2},
hence $\rank\nabla^2a(e_n)=n-1$ by the assumption
$\rank\nabla^2a(\xi)\geq n-1$.
In the setting \eqref{a2} and \eqref{rf} above, we have
\begin{equation}\label{rank:rem}
\rank\nabla^2\widetilde r(\psi(e_n))=0
\end{equation}
by Lemma \ref{l1}, where
$\widetilde r(\eta)=r(\eta_{k+1},\ldots,\eta_n)$.
Since $k=n-1$ in our case, we can see that $r$ is a function
of one variable and $r^{\prime\prime}$ vanishes
identically by \eqref{rank:rem} and the homogeneity of $r$.
Then $r$ is a polynomial of order one,
but is also positively homogeneous of order two.
Hence we can conclude that $r=0$ identically, and have the relation
\[
a(\xi)=(\sigma\circ\psi)(\xi),
\quad
\sigma(\eta)=
c_1\eta_1^2+\cdots+c_{n-1}\eta_{n-1}^2.
\]
Now, we have the estimates
\begin{align*}
&\n{\jp{x}^{-s} |\nabla \sigma(D_x)|^{1/2}
e^{it \sigma(D_x)}\varphi(x)}_{L^2\p{\R_t\times\R^n_x}}
\leq C\n{\varphi}_{L^2\p{\R^n_x}},
\\
&\n{\abs{x}^{\alpha-1}|\nabla \sigma(D_x)|^{\alpha}
e^{it\sigma(D_x)}\varphi(x)}_{L^2\p{\R_t\times\R^n_x}}
\leq C\n{\varphi}_{L^2\p{\R^n_x}},
\end{align*}
if we use the trivial inequalities $\jp{x}^{-s}\leq\jp{x'}^{-s}$
and $|x|^{\alpha-1}\leq|x'|^{\alpha-1}$,
Theorems \ref{M:H1} and \ref{M:H2} with respect to $x'$, and
the Plancherel's theorem in $x_n$,
where $x=(x',x_n)$ and $x'=(x_1,\ldots, x_{n-1})$.
On account of Theorem \ref{Th:caninv}
and the $L^2_{-s}$, $\Dot{L}^2_{\alpha-1}$--boundedness
of the operators $I_{\psi,\gamma}$ and $I_{\psi,\gamma}^{-1}$
for $(1/2<)s<n/2$, $-n/2<\alpha-1(<-1/2)$ (see Theorem \ref{Th:L'2k}),
we have the conclusion.
\end{proof}
\medskip
\begin{example}\label{nondispersive3}
The function $a(\xi)=b(\xi)^2$ satisfies
condition \eqref{a1} with $k=1$, where $b(\xi)$
is a positively homogeneous function of order one
such that $\nabla b(\xi)\neq0$ ($\xi\neq0$).
Indeed, if $b(\xi)$ is elliptic, then 
$\nabla a(\xi)=2b(\xi)\nabla b(\xi)\neq0$
($\xi\neq0$).
If $b(\xi_0)=0$ at a point $\xi_0\neq0$, then $\nabla a(\xi_0)=0$
and further differentiation immediately yields
$\nabla^2 a(\xi_0)=2^t\nabla b(\xi_0) \nabla b(\xi_0)$,
and clearly we have $\rank\nabla^2 a(\xi_0)\geq1$.
Especially in the case $n=2$, $a(\xi)$ meets the condition in
Theorem \ref{th1}.
As an example, we consider
\[
a(\xi)=\frac{\xi_1^2\xi_2^2}{\xi_1^2+\xi_2^2},
\]
where we write $\xi=(\xi_1,\xi_2)$.
Setting $b(\xi)=\xi_1\xi_2/|\xi|$, we clearly have $a(\xi)=b(\xi)^2$
and
\[
\nabla b(\xi)
=\p{\frac{\xi_2^3}{|\xi|^3},\frac{\xi_1^3}{|\xi|^3}},
\]
hence $\nabla b(\xi)\neq0$ ($\xi\neq0$).
Although $\nabla a(\xi)=0$ on the lines $\xi_1=0$ and $\xi_2=0$,
we have estimate \eqref{EQ:maininv} in virtue of
Theorem \ref{th1}.
This is an illustration of a smoothing estimate for the 
Cauchy problem for an equation like
$$
i\partial_t \Delta u+D_1^2D_2^2 u=0,
$$
which can be reduced to the second order non-dispersive
pseudo-differential equation with symbol $a(\xi)$ above.
Similarly, we have estimates \eqref{EQ:maininv} 
and \eqref{EQ:maininv2}
for more general case
\[
a(\xi)=\frac{\xi_1^2\xi_2^2}{\xi_1^2+\xi_2^2}
+\xi^2_3+\cdots+\xi_n^2,
\]
where we write $\xi=(\xi_1,\ldots,\xi_n)$
since we obtain $\rank\nabla^2a(\xi)\geq n-1$
from the observation above.
\end{example}
\medskip
Next we consider more general operators $a(\xi)$ of order $m$
which may have some lower order terms.
Then even the most favourable case $\det\nabla^2a(\xi)\not=0$
does not imply the dispersive assumption $\nabla a(\xi)\not=0$.
The method of canonical transforms, 
however, can also allow us to treat this problem by
obtaining localised estimates near points $\xi$ 
where $\nabla a(\xi)=0$.
\par
Assume that $\xi_0$ is a non-degenerate critical point
of $a(\xi)$, that is, that we have $\nabla a(\xi_0)=0$
and $\det \nabla^2 a(\xi_0)\not=0$. Let us microlocalise
around $\xi_0$, so that we only look at what happens around $\xi_0$.
In this case, the order of the symbol $a(\xi)$ does not play any role
and we do not distinguish between the main part and lower
order terms.
Let $\Gamma$ denote a sufficiently small
open bounded neighbourhood
of $\xi_0$ so that $\xi_0$ is the only critical point of
$a(\xi)$ in $\Gamma$. Since $\nabla^2 a(\xi_0)$ is symmetric and
non-degenerate, we may assume
$\nabla^2 a(\xi_0)={\rm diag}\{\pm 1,\cdots,\pm 1\}$
by a linear transformation.
By Morse lemma for $a(\xi)$, there exists a diffeomorphism
$\psi:\Gamma\to\widetilde\Gamma\subset\Rn$ 
with an open bounded neighbourhood
of the origin such that
\[
a(\xi)=(\sigma\circ\psi)(\xi),
\quad
\sigma(\eta)=
c_1\eta_1^2+\cdots+c_n\eta_n^2,
\]
where $\eta=(\eta_1,\ldots,\eta_n)$ and $c_j=\pm1$ ($j=1,2,\ldots,n$).
From Theorems \ref{M:H1} and \ref{M:H2} applied to operator
$\sigma(D_x)$, we obtain the estimates
\begin{equation}\label{EQ:criticalp1}
\n{\jp{x}^{-s} |\nabla \sigma(D_x)|^{1/2}
e^{it \sigma(D_x)}\varphi(x)}_{L^2\p{\R_t\times\R^n_x}}
\leq C\n{\varphi}_{L^2\p{\R^n_x}}
\qquad(s>1/2).
\end{equation}
Hence by Theorem \ref{Th:caninv}, together with
the $L^2_{-s}$--boundedness of the operators
$I_{\psi,\gamma}$ and $I_{\psi,\gamma}^{-1}$
(which is assured by Theorem \ref{Th:L2k}),
we have estimate \eqref{EQ:criticalp1} with $\sigma(D_x)$
replaced by $a(D_x)$ assuming $\supp\widehat \va\subset\Gamma$.
\par
On the other hand, we have the same estimate
for general $\va$ by Corollary \ref{M:L5} if we assume (L)
(see the argument after estimate \eqref{EQ:maininv}).
The above argument, however, assures that the 
following weak assumption
is also sufficient if $a(\xi)$ has finitely many 
critical points and
they are non-degenerate:
\medskip
\begin{equation}\tag{{\bf L$'$}}
\begin{aligned}
&a(\xi)\in C^\infty(\R^n),\qquad 
|\nabla a(\xi)|\geq C\jp{\xi}^{m-1}\quad(\text{for large $\xi\in\R^n$})\quad
\textrm{for some}\; C>0,
\\
&|\partial^\alpha\p{a(\xi)-a_m(\xi)}|\leq C_\alpha\abs{\xi}^{m-1-|\alpha|}
\quad\text{for all multi-indices $\alpha$ and all $|\xi|>>1$}.
\end{aligned}
\end{equation}
\medskip
Thus, we have established the following result:
\medskip
\begin{thm}\label{THM:isolated-critical}
Let $a\in C^\infty(\Rn)$ be real-valued and assume that it
has finitely many critical points,
all of which are non-degenerate.
Assume also {\rm{(L$'$)}}.
Suppose $n\geq1$, $m\geq1$, and $s>1/2$.
Then we have
\[
\n{\jp{x}^{-s}|\nabla a(D_x)|^{1/2}
e^{it a(D_x)}\varphi(x)}_{L^2\p{\R_t\times\R^n_x}}
\leq C\n{\varphi}_{L^2\p{\R^n_x}}.
\]
\end{thm}
\bigskip

\section{Relativistic Schr\"odinger, wave, and Klein--Gordon equations}
\label{SECTION:relative}
In Section \ref{SECTION:main}, we gave a criteria
Corollary \ref{COR:RStype}
for smoothing estimates to hold in the radially symmetric case.
Such subject has been also investigated
by Walther \cite{Wa2}, and he derived another type of criteria
based on certain integrals involving Bessel functions 
and their asymptotics. However, the approach presented in
this paper applies to such estimates in an
essentially different way in the sense that instead
of verifying convergence of infinitely many
integrals involving expressions based on special
functions we simply compare the estimate we want to have
to one that we already know to hold (in a model case or
otherwise).
\par
A typical direct application of Corollary \ref{COR:RStype}
is to the relativistic Schr\"odinger
type equations
\begin{equation}
\tag{{Relativistic Schr\"odinger}}
\label{eq:relative}
\left\{
\begin{aligned}
\p{i\partial_t-\sqrt{1-\Delta_x}}\,u(t,x)&=0,\\
u(0,x)&=\varphi(x).
\end{aligned}
\right.
\end{equation}
In \cite{BN}, Ben-Artzi and Nemirovsky proved the following
results.
Suppose first that $h\in C^1(\R_+)$ is real valued,
$h^\prime>0$, and $h^\prime$
is locally H\"older continuous.
Then, it follows
that $h(-\Delta_x)$ is self-adjoint in $L^2(\R^n)$ and
its spectrum is absolutely continuous and satisfies
$\sigma(h(-\Delta_x))=\overline{[h(0),h(\infty)]},$
where $h(\infty)=\lim_{\theta\to\infty} h(\theta).$
Suppose further that $h^\prime(\theta)$ satisfies a uniform 
H\"older condition near $\theta=0$ and that $h^\prime(0)>0$.
We remark that then we have
\begin{equation}\label{EQ:BNs}
h^\prime(\theta)\geq C
\;\textrm{ as }\;
\theta\searrow0
\end{equation}
for some $C>0$.
Assuming also $n\geq 3$ and
\begin{equation}\label{EQ:BNa}
h^\prime(\theta)\geq \frac{C}{\sqrt{\theta}}
\;\textrm{ as }\;
\theta\to+\infty
\end{equation}
for some $C>0$,
Ben-Artzi and Nemirovsky proved the estimate
\begin{equation}\label{EQ:BN}
\|\jp{x}^{-1} e^{-it h(-\Delta_x)} \varphi\|_
{L^2(\R_t\times\R_x^n)} \leq C\|\varphi\|_\Lx.
\end{equation}
for the solution $u(t,x)=e^{-it h(-\Delta_x)} \varphi$ to
the equation
\begin{equation}\label{EQ:hop}
\left\{
\begin{aligned}
\p{i\partial_t-h(-\Delta_x)}\,u(t,x)&=0,\\
u(0,x)&=\varphi(x).
\end{aligned}
\right.
\end{equation}
In particular, for $h(\theta)=\sqrt{1+\theta}$,
this leads to the time global estimate for the
relativistic Schr\"odinger equation: 
\begin{equation}\label{EQ:BNrs}
\|\jp{x}^{-1} e^{-it \sqrt{1-\Delta_x}} \varphi\|_
{L^2(\R_t\times\R_x^n)} \leq C\|\varphi\|_\Lx.
\end{equation}
We remark that the order of the weight $\jp{x}^{-1}$
in estimate \eqref{EQ:BNrs}
is sharp (see Walther \cite{Wa2}, for example, or Section
\ref{SECTION:invariant}). However, it can still be refined, and
this will be done in Theorem \ref{COR:RS}.
\par
The proof of \cite{BN} is based on the limiting absorption
principle for the resolvent of the operator $h(-\Delta_x)$.
But the comparison principle also allows us to 
get a simple proof of several refinements of estimate
\eqref{EQ:BN}.
Now we remark that, by looking at invariant
estimates \eqref{EQ:maininv} and \eqref{EQ:maininv2} in Section
\ref{SECTION:invariant} for equation \eqref{EQ:hop},
we should expect the estimates
of the form 
\begin{align*}
&\n{\jp{x}^{-s}|D_x h^\prime(-\Delta_x)|^{1/2}
e^{-it h(-\Delta_x)}\varphi(x)}_\L2tx
 \leq C\n{\varphi}_\Lx , \\
&\n{|x|^{\alpha-m/2}|D_x h^\prime(-\Delta_x)|^{\alpha/(m-1)}
e^{-it h(-\Delta_x)}\varphi(x)}_\L2tx
 \leq C\n{\varphi}_\Lx
\quad(m\neq1),
\end{align*}
where the order $m$ of the operator $h(-\Delta_x)$
has a different meaning for low frequency ($m=2$) and high frequency
($m=1$). In fact, these estimate can be shown using the
comparison principle in Theorem \ref{prop:dim1eqmod} for
radially symmetric operators, also without assumptions
\eqref{EQ:BNs} and \eqref{EQ:BNa}. In a special case with
conditions
\eqref{EQ:BNs} and \eqref{EQ:BNa}, we get
the following realisation of these estimates:
\medskip
\begin{thm}\label{COR:RS}
Suppose $n\geq1$, $s>1/2$, and $1-n/2<\alpha<1/2$.
Let $h\in C^1(\R_+)$ be a real-valued and strictly increasing function
which satisfies \eqref{EQ:BNs} and \eqref{EQ:BNa}.
Let $\chi\in C_0^\infty(\Rn)$ be equal to one in a 
neighbourhood of the origin.
Then we have
\begin{align}
\label{EQ:RStype1slf}
&\n{\jp{x}^{-s}|D_x|^{1/2}
e^{-it h(-\Delta_x)}\varphi_l(x)}_\L2tx
\leq C\n{\varphi_l}_\Lx, \\ 
\label{EQ:RStype1s}
&\n{\jp{x}^{-s}
e^{-it h(-\Delta_x)}\varphi_h(x)}_\L2tx
\leq C\n{\varphi_h}_\Lx, \\ 
\label{EQ:RStypehslf}
&\n{|x|^{\alpha-1}|D_x|^{\alpha}
e^{-it h(-\Delta_x)}\varphi_l(x)}_\L2tx
\leq C\n{\varphi_l}_\Lx, \\
\label{EQ:RStypehs}
&\n{|x|^{\alpha-1}|D_x|^{\alpha-1/2}
e^{-it h(-\Delta_x)}\varphi_h(x)}_\L2tx
\leq C\n{\varphi_h}_\Lx,
\end{align}
where $\varphi_l=\chi(D_x)\varphi$ and $\varphi_h=(1-\chi(D_x))\varphi$.
Consequently, if $n\geq3$, we have
\begin{align}\label{EQ:RS3n3}
\n{\jp{x}^{-1}e^{-ith(-\Delta_x)}
\varphi(x)}_\L2tx\leq C\n{\varphi}_\Lx.
\end{align}
If $n=2$ and $r>1$, we have
\begin{align}\label{EQ:RS3n2}
\n{\jp{x}^{-r}e^{-ith(-\Delta_x)}
\varphi(x)}_\L2tx\leq C\n{\varphi}_\Lx.
\end{align}
\end{thm}
\medskip
We note straight away that estimates \eqref{EQ:RStype1slf}
and \eqref{EQ:RStypehslf}
improve Ben-Artzi and Nemirovsky's estimate
\eqref{EQ:BN} for the low frequency part,
while \eqref{EQ:RStype1s} also improves the
weight given in \eqref{EQ:BN}
for the high frequency part.
From this point of view, we can see that estimate \eqref{EQ:BN}
does only capture estimate \eqref{EQ:RStypehslf} with
$\alpha=0$ for the low frequency part of the smoothing.
In fact, \eqref{EQ:RStypehslf} with $\alpha=0$ improves
the low frequency part of \eqref{EQ:BN} to the better
weight $|x|^{-1}$ in \eqref{EQ:RStypehslf}, compared to
$\jp{x}^{-1}$ in \eqref{EQ:BN}.
\medskip
\begin{proof}
Taking $f(\rho)=-h(\rho^2)$, condition
\eqref{EQ:BNs} implies that $|f^\prime(\rho)|\geq C\rho$
as $\rho\searrow0$.
At the same time, condition \eqref{EQ:BNa} implies that
$|f^\prime(\rho)|=2\rho h^\prime(\rho^2)\geq C$ as $\rho\to+\infty$.
It follows that we can take $\sigma(\rho)$ to be
$\sigma(\rho)=\rho^{1/2}$ for small $\rho$ and $\sigma(\rho)=1$
for large $\rho$ to meet condition \eqref{EQ:RStypeass} in Corollary
\ref{COR:RStype}.
Then estimates \eqref{EQ:RStype1} and \eqref{EQ:RStypeh} imply
estimates \eqref{EQ:RStype1slf}--\eqref{EQ:RStypehs}.
Estimate \eqref{EQ:RS3n3} is just a 
consequence estimate \eqref{EQ:RStype1s} and 
\eqref{EQ:RStypehslf} with $\alpha=0$,
which we can take to meet $1-n/2<\alpha<1/2$ if $n\geq3$.
In the case $n=2$, instead of \eqref{EQ:RStypehs} with $\alpha=0$,
we alternatively use the estimate
\[
\n{\jp{x}^{-r}e^{-ith(-\Delta_x)}
\varphi_l(x)}_\L2tx\leq C\n{\varphi_l}_\Lx
\qquad (r>1),
\]
which can be easily given by the comparison
(use Theorem \ref{prop:dim1eqmod}) with the estimate
\[
\n{\jp{x}^{-r}e^{it\Delta_x}
\varphi_l(x)}_\L2tx\leq C\n{\varphi_l}_\Lx
\qquad (r>1)
\]
for Schr\"odinger equations in the case $n=2$.
This type of estimate can be found in
Ben-Artzi and Klainerman \cite{BK} or Walther \cite{Wa1}.
\end{proof}
\medskip
Taking $h(\theta)=\sqrt{1+\theta}$ in Theorem \ref{COR:RS} as a special case,
we obtain estimates
\eqref{EQ:RStype1slf}--\eqref{EQ:RS3n2} for solutions
to the relativistic Schr\"odinger equation.
For example, estimate \eqref{EQ:BNrs} is a special case of
estimate \eqref{EQ:RS3n3} (that is, estimate \eqref{EQ:BN}).
We can also observe the refinement of the weight
in \eqref{EQ:BNrs} for both high and low frequencies, given
by \eqref{EQ:RStype1s} and \eqref{EQ:RStypehslf} with $\alpha=0$,
to $\jp{x}^{-s}$ and $|x|^{-1}$, respectively.
We also remark that by the comparison principle for
radially symmetric operators,
all of these estimates are
equivalent to corresponding estimates for Schr\"odinger or
wave equation, which can be also derived from pointwise estimates in
one dimension as was explained in Section \ref{SECTION:main}.
More precisely, by Theorem \ref{prop:dim1eqmod}, we have the equalities
\begin{equation}\label{equiv-mult}
\begin{aligned}
\n{e^{-it\sqrt{1- \Delta_x}} \varphi(x)}_{L^2(\R_t)} 
&=\sqrt{2}\,\n{\jp{D_x}^{1/2}e^{it \Delta_x} \varphi(x)}_{L^2(\R_t)}
\\
&=\n{|D_x|^{-1/2}\jp{D_x}^{1/2}e^{\pm it\sqrt{- \Delta_x}} \varphi(x)}_
{L^2(\R_t)}
\end{aligned}
\end{equation}
for almost all $x\in\R^n$.
If fact, since $f(\rho)=-\sqrt{1+\rho^2}$ and $g(\rho)=-\rho^2$
satisfy $1/|f'(\rho)|^{1/2}=|2f(\rho)|^{1/2}/|g'(\rho)|^{1/2}$,
we have the first equality.
The proof of the second one is similar.
Then multiplying appropriate weight functions to the both sides of equalities
\eqref{equiv-mult} and integrating them in $x$ imply the equivalence
of the estimates.
\par
For example, by \eqref{equiv-mult}, we have the equivalence of the estimate
\begin{equation}\label{eq:schrod}
\|\jp{x}^{-1} \jp{D_x}^{1/2}e^{it \Delta_x} \varphi\|_
{L^2(\R_t\times\R_x^n)} \leq C\|\varphi\|_\Lx
\end{equation}
for the standard Schr\"odinger equation
and estimate \eqref{EQ:BNrs} 
for the relativistic Schr\"odinger equation.
We remark that Corollary \ref{M:H3} also assures 
estimate \eqref{eq:schrod}
in the case $n\ge3$, so we have
\medskip
\begin{thm}\label{equiv:rel-sch}
Let $n\geq3$.
Then we have equivalent estimates \eqref{EQ:BNrs} and \eqref{eq:schrod}.
We also have the equality
\begin{equation}\label{EQ:SHequiv}
\|\jp{x}^{-1} e^{-it\sqrt{1- \Delta_x}} \varphi\|_
{L^2(\R_t\times\R_x^n)} 
=\sqrt{2}
\|\jp{x}^{-1} \jp{D_x}^{1/2}e^{it \Delta_x} \varphi\|_
{L^2(\R_t\times\R_x^n)}.
\end{equation}
\end{thm}
\medskip
Such equivalence as in Theorem \ref{equiv:rel-sch}
was shown by Walther \cite{Wa2} (but without equivalence nor
without $\sqrt{2}$), who used an explicit
calculation using spherical harmonics and Bessel functions,
specific for the radially symmetric case,
but it is easy to see it if we use the comparison method.
Similar equivalence between the relativistic Schr\"odinger equation
and the wave equation can be also given by equality \eqref{equiv-mult}:
\begin{align*}
&\n{\jp{x}^{-s}
e^{\pm it \sqrt{-\Delta_x}}\varphi_l(x)}_\L2tx
\sim
\n{\jp{x}^{-s}|D_x|^{1/2}
e^{-it \sqrt{1-\Delta_x}}\varphi_l(x)}_\L2tx,
\\
&\n{\jp{x}^{-s}
e^{\pm it \sqrt{-\Delta_x}}\varphi_h(x)}_\L2tx
\sim
\n{\jp{x}^{-s}
e^{-it \sqrt{1-\Delta_x}}\varphi_h(x)}_\L2tx,
\\
& \n{|x|^{\beta-1/2}|D_x|^{\beta}
e^{\pm it \sqrt{-\Delta_x}}\varphi_l(x)}_\L2tx
\sim
\n{|x|^{\beta-1/2}|D_x|^{\beta+1/2}
e^{-it \sqrt{1-\Delta_x}}\varphi_l(x)}_\L2tx,
\\
& \n{|x|^{\beta-1/2}|D_x|^{\beta}
e^{\pm it \sqrt{-\Delta_x}}\varphi_h(x)}_\L2tx
\sim
\n{|x|^{\beta-1/2}|D_x|^{\beta}
e^{-it \sqrt{1-\Delta_x}}\varphi_h(x)}_\L2tx.
\end{align*}
As another consequence of Theorem \ref{COR:RS}, we have the estimates
\begin{equation}\label{wave-est}
\begin{aligned}
&\n{\jp{x}^{-s}
e^{\pm it \sqrt{-\Delta_x}}\varphi(x)}_\L2tx
\leq C\n{\varphi}_\Lx
\qquad(s>1/2),
 \\ 
&\n{|x|^{\beta-1/2}|D_x|^{\beta}
e^{\pm it \sqrt{-\Delta_x}}\varphi(x)}_\L2tx
\leq C\n{\varphi}_\Lx
\quad((1-n)/2<\beta<0)
\end{aligned}
\end{equation}
for $n\geq1$.
Indeed, we also obtain the first
estimate from Theorem \ref{M:H1} and
the second estimate from Theorem \ref{M:H2}, or
from \eqref{core2g} with $m=1$.
We note that contrary to the relativistic Schr\"odinger equation,
here we get the same estimates for low and high frequencies.
The critical case of the second estimate with $\beta=0$
was analysed by the authors in \cite{RS3} and it was shown
that its modification still holds by introducing
the Laplace-Beltrami operator on the sphere into the estimate.
In fact, that analysis was done for general 
second order strictly hyperbolic
equations with homogeneous symbols with critical sets 
associated to some sets related to the classical orbits.
See Section \ref{SECTION:other} for further information.
\par
Now we apply estimate \eqref{wave-est} to 
the wave equation
\begin{equation}
\tag{{Wave Equation}}
\left\{
\begin{aligned}
\partial_t^2 u -\Delta u&=0,\\
u(0,x)&=u_0(x), \\
\partial_t u(0,x)&=v_0(x).
\end{aligned}
\right.
\end{equation}
Then we have estimates
\begin{equation}\label{WEest}
\begin{aligned}
&\n{\jp{x}^{-s}u}_\L2tx
  \leq C\p{||u_0||_{L^2(\R^n_x)}+||D_x|^{-1}v_0||_{L^2(\R^n_x)}},
\\
&\n{|x|^{\beta-1/2}|D_x|^{\beta}u}_\L2tx
  \leq C\p{||u_0||_{L^2(\R^n_x)}+|||D_x|^{-1}v_0||_{L^2(\R^n_x)}},
\end{aligned}
\end{equation}
where we can take any $n\geq 1$,
$s>1/2$, and $(1-n)/2<\beta<0$.
These estimates have been previously established
for $n\geq 3$ and $-1<\beta<0$ (see Ben-Artzi \cite{Be},
where spectral methods were used).
These estimates follow now from the smoothing estimates for
propagators $e^{\pm it\sqrt{-\Delta_x}}$, which can be obtained by
the comparison principle.
We note that
the usual way of relating smoothing estimates of wave
and Schr\"odinger equation goes via a change of
variables in the corresponding restriction theorems
(see, for example, \cite{RS3}). Now we can relate them 
directly by 
the comparison principle in Theorem \ref{prop:dim1eqmod}.
We also note that in the case of $n\geq 3$ and $\beta=-1/2$ 
the best constant $\sqrt{\frac{2\pi}{n-2}}$ 
in the second inequality is given by
\eqref{EQ:Simon-const} with $m=1$.
\par
Let us finally state smoothing estimates for the Klein--Gordon
equation
\begin{equation}
\tag{{Klein--Gordon}}
\left\{
\begin{aligned}
\partial_t^2 u -\Delta u+\mu^2 u&=0,\\
u(0,x)&=u_0(x), \\
\partial_t u(0,x)&=v_0(x),
\end{aligned}
\right.
\end{equation}
for $\mu>0$. In the case $n\geq 3$ the estimate
\begin{equation}\label{EQ:KGest}
\n{\jp{x}^{-1}u}_\L2tx \leq C\p{||u_0||_{L^2(\R^n_x)}
+||(\mu^2-\Delta)^{-1/2}v_0||_{L^2(\R^n_x)}}
\end{equation}
was given in \cite{Be}. Since propagators here are 
of the form $e^{\pm it\sqrt{\mu^2-\Delta_x}}$, we can
apply Theorem \ref{COR:RS} with $h(\theta)=\sqrt{\mu^2+\theta}.$
In particular, this implies estimate \eqref{EQ:KGest},
as well as all of its refinements given by
Theorem \ref{COR:RS}. In particular, we get
the weight $\jp{x}^{-s}$
with $s>1$ in the case of $n=2$, and better
weights for high frequencies in all dimensions $n\geq 1$.
\medskip

\section{Model estimates for inhomogeneous equations}
\label{SECTION:model-inh}
We now turn to deal with inhomogeneous equations, for which
we also have similar smoothing estimates. Such estimates are
necessary for nonlinear applications, and they can be
obtained by further developments of the presented methods.
Note that
\[
u(t,x)
=-i\int^t_0e^{i(t-\tau)a(D_x)}f(\tau,x)\,d\tau
\]
solves the equation
\[
\left\{
\begin{aligned}
\p{i\partial_t+a(D_x)}\,u(t,x)&=f(t,x)\quad\text{in $\R_t\times\R^n_x$},\\
u(0,x)&=0\quad\text{in $\R^n_x$}.
\end{aligned}
\right.
\]
We will give model estimates for it below,
where we write $x=(x_1,x_2,\ldots,x_n)\in\R^n$ and $D_x=(D_1,D_2\ldots,D_n)$.
We also write $x=x_1$, $D_x=D_1$ in the case $n=1$, 
and $(x,y)=(x_1,x_2)$, $(D_x,D_y)=(D_1,D_2)$ in the case $n=2$.
\medskip
\begin{thm}\label{Prop:inhom}
Suppose $n=1$ and $m>0$.
Let $a(\xi)\in C^\infty\p{\R\setminus0}$ be a real-valued
function which satisfies $a(\lambda\xi)=\lambda^m a(\xi)$ for
all $\lambda>0$ and $\xi\neq0$.
Then we have
\begin{equation}\label{eq:inhomdim1}
\n{a'(D_x)
\int^t_0e^{i(t-\tau)a(D_x)}f(\tau,x)\,d\tau}_{L^2(\R_t)}
\leq C\int_\R \n{f(t,x)}_{L^2(\R_t)}\,dx
\end{equation}
for all $x\in\R$.
Suppose $n=2$ and $m>0$.
Then we have
\begin{multline}\label{eq:inhomdim2}
\n{|D_x|^{m-1}\int^t_0e^{i(t-\tau)|D_x|^{m-1}D_y}
f(\tau,x,y)\,d\tau}_{L^2(\R_t\times\R_x)}
\\
\leq
C\int_\R \n{f(t,x,y)}_{L^2(\R_t\times\R_x)}\,dy
\end{multline}
for all $y\in\R$.
\end{thm}
\begin{cor}\label{Th:model-hom}
Suppose $n\geq1$, $m>0$, and $s>1/2$.
Let $a(\xi)\in C^\infty\p{\R\setminus0}$ be a real-valued
function which satisfies $a(\lambda\xi)=\lambda^m a(\xi)$ for
all $\lambda>0$ and $\xi\neq0$.
Then we have
\[
\n{\jp{x_n}^{-s}a'(D_n)
\int^t_0e^{i(t-\tau)a(D_n)}f(\tau,x)\,d\tau}
_{L^2(\R_t\times\R^n_x)}
\leq
 C\n{\jp{x_n}^{s}f(t,x)}_{L^2(\R_t\times\R_x^n)}.
\]
Suppose $n\geq2$, $m>0$, and $s>1/2$.
Then we have
\[
\n{\jp{x_1}^{-s}|D_n|^{m-1}\int^t_0e^{i(t-\tau)D_1|D_n|^{m-1}}f(\tau,x)\,d\tau}
_{L^2(\R_t\times\R^n_x)}
\leq
 C\n{\jp{x_1}^{s}f(t,x)}_{L^2(\R_t\times\R_x^n)}.
\]
\end{cor}
\medskip
Theorem \ref{Prop:inhom} with the case $n=1$
is a unification of the results
by Kenig, Ponce and Vega who treated the
cases $a(\xi)=\xi^2$ (\cite[p.258]{KPV3}), $a(\xi)=|\xi|\xi$
(\cite[p.160]{KPV4}),
and $a(\xi)=\xi^3$ (\cite[p.533]{KPV2}).
Corollary \ref{Th:model-hom} is a straightforward result of 
Theorem \ref{Prop:inhom} and Cauchy--Schwarz's inequality.
\par
Since we unfortunately do not know the comparison principle
for inhomogeneous equations, we cannot reduce
Theorem \ref{Prop:inhom} to more elementary estimates
as we can successfully do that for homogeneous equations
in Section \ref{SECTION:model}.
Hence we will give a direct proof to Theorem \ref{Prop:inhom}.
Note that we have another expression 
of the solution to inhomogeneous equation
\[
\left\{
\begin{aligned}
\p{i\partial_t+a(D_x)}\,u(t,x)&=f(t,x)\quad\text{in $\R_t\times\R^n_x$},\\
u(0,x)&=0\quad\text{in $\R^n_x$},
\end{aligned}
\right.
\]
using the weak limit $R(\tau\pm i0)$ 
of the
resolvent $R(\tau\pm i\varepsilon)$ as $\varepsilon\searrow0$, where
$R(\lambda)=\p{a(D_x)-\lambda}^{-1}$:
\begin{equation}\label{resolv}
\begin{aligned}
u(t,x)
&=-i\int^t_0e^{i(t-\tau)a(D_x)}f(\tau,x)\,d\tau.
\\
&=\FT_\tau^{-1}R(\tau-i0)\FT_tf^+
+\FT_\tau^{-1}R(\tau+i0)\FT_tf^-
\end{aligned}
\end{equation}
(see Sugimoto \cite {Su1} and Chihara \cite{Ch}).
Here $\FT_t$ denotes the Fourier Transformation in $t$ and 
$\FT_\tau^{-1}$
its inverse, and $f^{\pm}(t,x)=f(t,x)Y(\pm t)$ is
the characteristic function $Y(t)$ of the set $\b{t\in\R:\,t>0}$.
\begin{proof}[Proof of Estimate \eqref{eq:inhomdim1}]
Let us use a variant of 
the argument of Chihara \cite[Section 4]{Ch}.
We set
$R(\lambda)=\p{a(D_x)-\lambda}^{-1}$
and show the estimate
\[
\abs{a^\prime (D_x)R(s\pm i0)g(x)}\leq C\int_\R|g(x)|\,dx,
\]
where $C>0$ is a constant independent of $s\in\R$, 
$x\in\R$ and $g\in L^1(\R)$.
Then, on account of the expression \eqref{resolv},
Plancherel's theorem, and Minkowski's inequality,
we have the desired result.
For this purpose, we consider the kernel
\[
k(s,x)=\FT^{-1}\br{a'(\xi)\p{a(\xi)-(s\pm i0)}^{-1}}(x)
\]
and show its uniform boundedness.
By the scaling argument, everything is reduced to show
the estimates
\[
\sup_{x\in \R}|k(\pm1,x)|\leq C
\qquad\text{and}\qquad \sup_{x\in \R}|k(0,x)|\leq C.
\]
By using an appropriate partition of
unity $\widehat{\phi}_1(\xi)+\widehat{\phi}_2(\xi)+\widehat{\phi}_3(\xi)=1$,
we split $k(\pm1, x)$ into the corresponding three parts
$k=k_1+k_2+k_3$,
where $\widehat{\phi}_1$ has its support near the origin,
$\widehat{\phi}_2$ near the point
$\xi^m=\pm1$, and $\widehat{\phi}_3$ away from these points.
The estimate for $k_1$ is trivial.
The other estimates are reduced to the boundedness of
\begin{equation}\label{heaviside}
k_0^\pm(x)=\FT^{-1}\br{\p{\xi\pm i0)}^{-1}}(x)=
\mp i\sqrt{2\pi}Y(\pm x).
\end{equation}
In fact, 
\[
k_2(\pm1,x)=\FT^{-1}\br{\p{\xi-(\alpha\pm i0)}^{-1}\widehat{\psi}(\xi)}(x)
=(e^{i\alpha x}k_0^\mp)*\psi(x)
\]
where $\alpha\in\R$ is a point which solves $a(\alpha)=\pm1$, and
\[
\widehat{\psi}(\xi)
=a'(\xi)\frac{\xi-\alpha}{a(\xi)-(\pm1)}
\widehat{\phi}_2(\xi)\in C^\infty_0(\R).\]
Furthermore, if we notice
\[
\frac{a'(\xi)}{a(\xi)-s}
=m\p{\frac{s}{(a(\xi)-s)\xi}+\frac1{\xi}},
\]
we have
\[
\frac1m
k_3(\pm1,x)=\pm\FT^{-1}\br{\frac{\widehat{\phi}_3(\xi)}
{(a(\xi)\mp1)\xi}}(x)
+k^\pm_0(x)-k^\pm_0*\p{\phi_1(x)+\phi_2(x)}.
\]
It is easy to deduce the estimates for $k_2$ and $k_3$.
It is also easy to verify
\[
\frac{a'(\xi)}{a(\xi)\pm i0}
=\frac{m}{\xi\pm i0}+c\delta
\]
with a constant $c$ and Dirac's delta function $\delta$,
and have the estimate for $k(0,x)$.
\end{proof}
\begin{proof}[Proof of Estimate \eqref{eq:inhomdim2}]
We set
$R(\lambda)=\p{|D_x|^{m-1}D_y-\lambda}^{-1}$
and show the estimate
\[
\n{|D_x|^{m-1}R(s\pm i0)g(x,y)}_{L^2(\R_x)}
\leq C\int\n{g(x,y)}_{L^2(\R_x)}\,dy,
\]
where $C>0$ is a constant independent of $s\in\R$,
$y\in\R$ and $g\in L^1(\R^2)$.
Then, by the expression \eqref{resolv}, Plancherel's theorem,
and Minkowski's inequality again,
we have the desired result.
\par
First we note, we may assume $\widehat{g}(\xi,\eta)=0$ for $\xi<0$.
Then we have
\begin{align*}
&|D_x|^{m-1}R(s\pm i0)g(x,y)
\\
=&\p{2\pi}^{-2}\int^\infty_0\int^\infty_{-\infty}
 e^{i(x\xi+y\eta)} |\xi|^{m-1}\p{|\xi|^{m-1}\eta-(s\pm i0)}^{-1}
 \widehat{g}(\xi,\eta)\,d\xi d\eta
\\
=&\p{2\pi}^{-2}\int^\infty_0\int^\infty_{-\infty}
 e^{ix\xi} |\xi|^{m-1}\p{|\xi|^{m-1}\eta-(s\pm i0)}^{-1}
 \widehat{g_y}(\xi,\eta)\,d\xi d\eta
\\
=&\p{2\pi}^{-2}\int^\infty_{-\infty}\int^\infty_{0}
 e^{ixb} \p{a-(s\pm i0)}^{-1}\widehat{g_y}(b,ab^{-(m-1)})\,dadb
\\
=&\p{2\pi}^{-2}\int^\infty_{-\infty}\int^\infty_{0}
 e^{ixb} \FT_a\br{\p{a-(s\pm i0)}^{-1}}
 \FT_a^{-1}\br{\widehat{g_y}(b,ab^{-(m-1)})}\,dadb
\\
=&\p{2\pi}^{-1}\int^\infty_{-\infty}\int^\infty_{0}
 e^{ixb} e^{-isa}k_0^\mp(-a)
 b^{m-1}\widetilde{g_y}(b,ab^{m-1})\,dadb,
\end{align*}
hence we have
\[
\FT_x\br{|D_x|^{m-1}R(s\pm i0)g(x,y)}(b)
=\int^\infty_{-\infty}
  e^{-isa}k_0^\mp(-a)
 b^{m-1}\widetilde{g_y}(b,ab^{m-1})\,da
\]
for $b\geq0$, and it vanishes for $b<0$.
Here $g_y(x,\,\cdot\,)=g(x,\,\cdot\,+y)$,
and $\widetilde{g_y}$ denotes its
partial Fourier transform with respect to the first variable.
We have also used here the change of variables
$a=\xi^{m-1}\eta$, $b=\xi$ and Parseval's formula.
Note that $\partial(a,b)/\partial(\xi,\eta)=b^{m-1}$ and
$k_0^\mp$ is a bounded function defined by \eqref{heaviside}.
Then we have the estimate
\begin{align*}
\abs{
\FT_x\br{|D_x|^{m-1}R(s\pm i0)g(x,y)}(b)
}
&\leq\sqrt{2\pi}
\int^\infty_{-\infty}
 \abs{b^{m-1}\widetilde{g_y}(b,ab^{m-1})}\,da
\\
&=\sqrt{2\pi}\int^\infty_{-\infty}
 \abs{\widetilde{g_y}(b,a)}\,da,
\end{align*}
and, by Plancherel's theorem and Minkowski's inequality, we have
\begin{align*}
\n{|D_x|^{m-1}R(s\pm i0)g(x,y)}_{L^2(\R_x)}
&\leq
\sqrt{2\pi}\int^\infty_{-\infty}
 \n{g_y(x,a)}_{L^2(\R_x)}\,da
\\
&=\sqrt{2\pi}\int^\infty_{-\infty}
 \n{g(x,y)}_{L^2(\R_x)}\,dy,
\end{align*}
which is the desired estimate.
\end{proof}
\medskip

\section{Smoothing estimates for dispersive inhomogeneous equations}
\label{SECTION:inhomogeneous}
Let us consider the inhomogeneous equation
\[
\left\{
\begin{aligned}
\p{i\partial_t+a(D_x)}\,u(t,x)&=f(t,x)\quad\text{in $\R_t\times\R^n_x$},\\
u(0,x)&=0\quad\text{in $\R^n_x$},
\end{aligned}
\right.
\]
where we always assume that 
function $a(\xi)$ is real-valued.
Let the principal part $a_m(\xi)\in C^\infty(\R^n\setminus0)$,
be a positively homogeneous function
of order $m$.
Recall the dispersive conditions we used in Section \ref{SECTION:main}:
\medskip
\begin{equation}\tag{{\bf H}}
a(\xi)=a_m(\xi),\qquad\nabla a_m(\xi)\neq0 \quad(\xi\in\R^n\setminus0),
\end{equation}
\begin{equation}\tag{{\bf L}}
\begin{aligned}
&a(\xi)\in C^\infty(\R^n),\qquad \nabla a(\xi)\neq0 \quad(\xi\in\R^n),
\quad\nabla a_m(\xi)\neq0 \quad (\xi\in\R^n\setminus0), 
\\
&|\partial^\alpha\p{a(\xi)-a_m(\xi)}|\leq C_\alpha\abs{\xi}^{m-1-|\alpha|}
\quad\text{for all multi-indices $\alpha$ and all $|\xi|\geq1$}.
\end{aligned}
\end{equation}
\medskip
The following is a counterpart of Theorem \ref{M:H1}
which treated homogeneous equations:
\medskip
\begin{thm}\label{Th:inhom}
Assume {\rm (H)}.
Suppose $m>0$ and $s>1/2$.
Then we have
\begin{equation}\label{eq:inhom:1}
\n{\jp{x}^{-s}|D_x|^{m-1}\int^t_0e^{i(t-\tau)a(D_x)}f(\tau,x)\,d\tau}_
{L^2(\R_t\times\R^n_x)}
\leq
 C\n{\jp{x}^s f(t,x)}_{L^2(\R_t\times\R^n_x)}
\end{equation}
in the case $n\geq2$, and
\begin{equation}\label{eq:inhom:2}
\n{\jp{x}^{-s}a'(D_x)\int^t_0e^{i(t-\tau)a(D_x)}f(\tau,x)\,d\tau}_
{L^2(\R_t\times\R_x)}
\leq
 C\n{\jp{x}^s f(t,x)}_{L^2(\R_t\times\R_x)}
\end{equation}
in the case $n=1$.
\end{thm}
\medskip
Chihara \cite{Ch} proved Theorem \ref{Th:inhom} with $m>1$
under the assumption \rm{(H)}.
As was pointed out in \cite[p.1958]{Ch},
we cannot replace $a'(D_x)$ by $|D_x|^{m-1}$ in
estimate \eqref{eq:inhom:2} for the case
$n=1$, but there is another explanation for this obstacle.
If we decompose
$f(t,x)=\chi_+(D_x)f(t,x)+\chi_-(D_x)f(t,x)$,
where $\chi_{\pm}(\xi)$ is a characteristic function of the 
set $\b{\xi\in\R\,:\,\pm\xi\geq0}$,
then we easily obtain
\begin{align*}
&\n{\jp{x}^{-s}|D_x|^{(m-1)/2}\int^t_0e^{i(t-\tau)a(D_x)}f(\tau,x)\,d\tau}_
{L^2(\R_t\times\R_x)}
\\
\leq
 &C\Bigl(
\n{\jp{x}^s |D_x|^{-(m-1)/2}f_+(t,x)}_{L^2(\R_t\times\R_x)}
 +\n{\jp{x}^s |D_x|^{-(m-1)/2}f_-(t,x)}_{L^2(\R_t\times\R_x)}
\Bigr)
\end{align*}
from Theorem \ref{Th:inhom}.
But we cannot justify the estimate
\begin{equation}\label{EQ:dim1problem}
 \n{\jp{x}^s|D_x|^{-(m-1)/2} f_\pm(t,x)}_{L^2(\R_t\times\R_x)}
\leq
 \n{\jp{x}^s|D_x|^{-(m-1)/2} f(t,x)}_{L^2(\R_t\times\R_x)}
\end{equation}
for $s>1/2$ by Lemma \ref{Prop:wtbdd}
because it requires $s<n/2$ and it is impossible for $n=1$.
\par
Similarly, as a counterpart of Theorem \ref{M:L4}, we have
\medskip
\begin{thm}\label{Th:inhom2}
Assume {\rm (L)}.
Suppose $n\geq1$, $m>0$, and $s>1/2$.
Then we have
\begin{equation}\label{eq:inhom:3}
\n{\jp{x}^{-s}\jp{D_x}^{m-1}\int^t_0e^{i(t-\tau)a(D_x)}f(\tau,x)\,d\tau}_
{L^2(\R_t\times\R^n_x)}
\leq
 C\n{\jp{x}^s f(t,x)}_{L^2(\R_t\times\R^n_x)}.
\end{equation}
\end{thm}
\medskip
The following result is a straightforward consequence
of Theorem \ref{Th:inhom2} and the $L^2_s$--boundedness
of $|D_x|^{(m-1)/2}\jp{D_x}^{-(m-1)/2}$
with $(1/2<)s<n/2$ and $m\geq1$ (which is assured by Lemma \ref{Prop:wtbdd}):
\medskip
\begin{cor}\label{Cor:inhom2}
Assume {\rm (L)}.
Suppose $n\geq2$, $m\geq1$, and $s>1/2$.
Then we have
\[
\n{\jp{x}^{-s}|D_x|^{m-1}\int^t_0e^{i(t-\tau)a(D_x)}f(\tau,x)\,d\tau}_
{L^2(\R_t\times\R^n_x)}
\leq
 C\n{\jp{x}^s f(t,x)}_{L^2(\R_t\times\R^n_x)}.
\]
\end{cor}
\medskip
\par
We remark that the same argument of canonical transformations
as used for homogeneous equations in
Section \ref{SECTION:main} works for inhomogeneous ones, as well.
That is, the proofs of Theorems \ref{Th:inhom} and \ref{Th:inhom2}
are carried out by reducing them to model estimates in
Corollary \ref{Th:model-hom}.
We omit the details because the argument is essentially the same,
but we just remark that we use the following slight 
modification of
Theorem \ref{Th:reduction}.
The only difference is that we need the weighted $L^2$--boundedness
of the operator $I_{\psi,q}^{-1}$
instead of just the $L^2$--boundedness of it induced by the boundedness
of $q(\xi)$:
\medskip
\begin{thm}\label{Th:redhom}
Assume that the operator $I_{\psi,\gamma}$ defined by \eqref{DefI0}
is $L^2(\R^n;w)$--bounded.
Suppose that we have the estimate
\[
\n{w(x)\rho(D_x)\int^t_0
e^{i(t-\tau)\sigma(D_x)}f(\tau,x)\,d\tau}_{L^2\p{\R_t\times\R^n_x}}
\leq C\n{v(x)f(t,x)}_{L^2\p{\R_t\times\R^n_x}}
\]
for all $f$ such that $\supp \FT_x f(t,\cdot)\subset\supp\widetilde\gamma$,
where $\widetilde\gamma=\gamma\circ\psi^{-1}$.
Also assume that
the operator
$I_{\psi,q}^{-1}$
defined by \eqref{DefI0}
with $q(\xi)=\p{\gamma\cdot\zeta}/\p{\rho\circ \psi}(\xi)$
is $L^2(\R^n;v)$--bounded.
Then we have
\[
\n{w(x)\zeta(D_x)\int^t_0
e^{i(t-\tau)a(D_x)}f(\tau,x)\,d\tau}_{L^2\p{\R_t\times\R^n_x}}
\leq C\n{v(x)f(t,x)}_{L^2\p{\R_t\times\R^n_x}}
\]
for all $f$ such that $\supp \FT_x f(t,\cdot)\subset\supp\gamma$,
where $a(\xi)=(\sigma\circ\psi)(\xi)$.
\end{thm}
\medskip
\par
The following is a counterpart of Theorem \ref{Th:HL}:
\medskip
\begin{thm}\label{Th:HLinh}
Assume {\rm (HL)}.
Suppose $n\geq1$, $m>0$, $s>1/2$, and $T>0$.
Then we have
\begin{multline*}
\int^T_0\n{\jp{x}^{-s}\jp{D_x}^{m-1}
\int^t_0e^{i(t-\tau)a(D_x)}f(\tau,x)\,d\tau}^2_{L^2(\R^n_x)}\,dt
\\
\leq
 C\int^T_0\n{\jp{x}^sf(t,x)}^2_{L^2(\R^n_x)}\,dt,
\end{multline*}
where $C>0$ is a constant depending on $T>0$.
\end{thm}
\medskip
\begin{proof}
By multiplying $\chi(D_x)$ and $(1-\chi)(D_x)$ to $f(t,x)$,
we decompose it into the sum 
of low frequency part and high frequency part,
where $\chi(\xi)$ is an appropriate cut-off function.
As in the proof of Theorem \ref{M:L4},
the estimate for the high frequency part can be reduced to
Corollaries \ref{Th:typeI} and \ref{Th:model-hom}
by using Theorems \ref{Th:reduction}
and \ref{Th:redhom},
together with the boundedness result Theorem \ref{Th:L2k}.
Here we note that, for $t\in[0,T]$,
\[
\int^t_0 e^{i(t-\tau)a(D_x)}f(\tau,x)\,d\tau
=
\int^t_0 e^{i(t-\tau)a(D_x)}\chi_{[0,T]}(\tau)f(\tau,x)\,d\tau,
\]
where $\chi_{[0,T]}$ denotes the characteristic 
function of the interval $[0,T]$.
The estimate for the low frequency part is trivial.
In fact, if $\supp_\xi \FT_xf(t,\xi)\subset
\br{\xi;|\xi|\leq R}$, we have
\begin{align*}
&\int^T_0
\n{\jp{x}^{-s}\jp{D_x}^{m-1}\int^t_0e^{i(t-\tau)a(D_x)}
 f(\tau,x)\,d\tau}^2_{L^2(\R^n_x)}
\,dt
\\
\leq
&\int^T_0
\p{
\int^T_0
\n{
\jp{D_x}^{m-1}e^{i(t-\tau)a(D_x)}f(\tau,x)
}_{L^2(\R^n_x)} \,d\tau
}^2
\,dt
\\
\leq &CT^2\jp{R}^{2(m-1)}
 \int^T_0\n{\jp{x}^sf(t,x)}^2_{L^2(\R^n_x)}\,dt.
\end{align*}
by Plancherel's theorem.
\end{proof}
\par
\medskip
If we combine Theorem \ref{Th:inhom} with Theorem \ref{M:H1},
we have a result for the equation
\begin{equation}\label{eq:inhom2}
\left\{
\begin{aligned}
\p{i\partial_t+a(D_x)}\,u(t,x)&=f(t,x)\quad\text{in $\R_t\times\R^n_x$},\\
u(0,x)&=\varphi(x)\quad\text{in $\R^n_x$}.
\end{aligned}
\right.
\end{equation}
\medskip
\begin{cor}\label{Cor:inhom+hom}
Assume {\rm (H)}.
Suppose $m>0$ and $s>1/2$.
Then the solution $u$ to equation \eqref{eq:inhom2}
satisfies
\begin{multline*}
\n{\jp{x}^{-s}|D_x|^{-(m-1)/2} a'(D_x)u(t,x)}
_{L^2(\R_t\times\R_x)}
\\
\leq
 C\p{\n{\varphi}_{L^2(\R)}
+\n{\jp{x}^s|D_x|^{-(m-1)/2}f(t,x)}_{L^2(\R_t\times\R_x)}}
\end{multline*}
in the case $n=1$, and
\begin{multline*}
\n{\jp{x}^{-s}|D_x|^{(m-1)/2}u(t,x)}
_{L^2(\R_t\times\R^n_x)}
\\
\leq
 C\p{\n{\varphi}_{L^2(\R^n)}
+\n{\jp{x}^s|D_x|^{-(m-1)/2}f(t,x)}_{L^2(\R_t\times\R^n_x)}}
\end{multline*}
in the case $n\geq2$.
\end{cor}
\medskip
If we combine Theorem \ref{Th:inhom2} with Theorem \ref{M:L4},
we have the following:
\begin{cor}\label{Cor:inhom+hom2}
\medskip
Assume {\rm (L)}.
Suppose $n\geq1$, $m>0$, and $s>1/2$.
Then the solution $u$ to equation \eqref{eq:inhom2}
satisfies
\begin{multline*}
\n{\jp{x}^{-s}\jp{D_x}^{(m-1)/2}u(t,x)}
_{L^2(\R_t\times\R_x)}
\\
\leq
 C\p{\n{\varphi}_{L^2(\R)}
+\n{\jp{x}^s\jp{D_x}^{-(m-1)/2}f(t,x)}_{L^2(\R_t\times\R_x)}}.
\end{multline*}
\end{cor}
\medskip
If we combine Theorem \ref{Th:HLinh} with Theorem \ref{Th:HL},
we have the following:
\medskip
\begin{cor}\label{Cor:HL}
Assume {\rm (HL)}.
Suppose $n\geq1$, $m>0$, $s>1/2$, and $T>0$.
Then the solution $u$ to equation \eqref{eq:inhom2} satisfies
\begin{multline*}
\int^T_0\n{\jp{x}^{-s}\jp{D_x}^{(m-1)/2}u(t,x)}^2_{L^2(\R^n_x)}\,dt
\\
\leq
 C\p{\n{\varphi}_{L^2(\R^n)}^2
+\int^T_0\n{\jp{x}^s\jp{D_x}^{-(m-1)/2}f(t,x)}^2_{L^2(\R^n_x)}\,dt},
\end{multline*}
where $C>0$ is a constant depending on $T>0$.
\end{cor}
\medskip
Corollary \ref{Cor:HL} is an extension of the result
by Hoshiro \cite{Ho2}, which treated the case that $a(\xi)$ 
is a polynomial.
The proof relied on Mourre's method, which is known in spectral
and scattering theories. Here we use the argument of canonical 
transformations, extending the result and simplifying the proof.
\bigskip

\section{Trace theorems}
\label{SECTION:other}
Another consequence of the proposed method of canonical
transforms is that we can carry out the geometric 
analysis of the
smoothing estimates leading to results relating the symbols 
with the location of the main singularities for solutions.
To exemplify this and to give an idea of how to use it in
problems at hand let us briefly mention the result that was
established by the authors in \cite{RS3} partly relying on
a variant of the method developed here.
\par
This concerns the critical case
($\alpha=1/2$) of the Kato--Yajima's estimate
\begin{equation}\label{EQ:critical-est1}
\n{|x|^{\alpha-1}|D_x|^{\alpha}e^{it\Delta_x}\varphi(x)}_
{L^2\p{\R_t\times\R^n_x}}
\leq C\n{\varphi}_{L^2\p{\R^n_x}},
\end{equation}
which holds for $1-n/2<\alpha<1/2$. In Section
\ref{SECTION:model} we argued that this estimate for
values of $\alpha$ close to $1/2$ implies the same estimate
for smaller $\alpha$ (see \eqref{core2''}). Thus, the
critical case of this estimate with $\alpha=1/2$ is important,
especially since it can be applied to the well-posedness problems
of the derivative nonlinear Schr\"odinger equations
(see \cite{RS5}). However, the estimate fails in the critical
case (see Watanabe \cite{W}, or Chapter 
\ref{SECTION:invariant} for more general results) 
and it is known that it is necessary to cut-off the
radial derivatives for the estimate to hold in the critical
case as well (see \cite{Su2}). This can be done by
replacing operator
$|D_x|^\alpha$ by the Laplace-Beltrami operator on
the sphere at the level $\alpha=1/2$. The method of canonical
transforms allows one to use any operator there as long as
its symbol vanishes on a certain set related to the symbol
of the Laplace operator (the sphere is this case).
\par
To explain
this precisely, let us formulate it for the equation
\begin{equation}\label{eq-cr1}
\left\{
\begin{aligned}
\p{i\partial_t+a(D_x)}\,u(t,x)&=0,\\
u(0,x)&=\varphi(x)\in L^2\p{\R^n_x},
\end{aligned}
\right.
\end{equation}
where real-valued function
$a(\xi)\in C^\infty\p{\R^n\setminus0}$ 
is elliptic and positively homogeneous of order two,
that is, it satisfies
$a(\xi)>0$ and $a(\lambda \xi)=\lambda^2 a(\xi)$ for $\lambda>0$
and $\xi\neq0$.
We remark that these condition assure assumption (H) 
with $m=2$ in
Section \ref{SECTION:main} since we have
$\nabla a(\xi)\neq0$ by the Euler's identity
$a(\xi)=1/2\nabla a(\xi)\cdot\xi$ and the ellipticity of $a(\xi)$.
The case $a(\xi)=|\xi|^2$ corresponds to the usual Laplacian
$a(D_x)=-\Delta_x$.
\par
Let us define
$\b{(x(t),y(t)):\,t\in\R}$ to be the classical orbit, 
that is, the solution
of the ordinary differential equation
\[
\left\{
\begin{aligned}
\dot{x}(t)&=\p{\nabla_\xi a}(\xi(t)), \quad\dot{\xi}(t)=0,
\\
x(0)&=0,\quad \xi(0)=\xi_0,
\end{aligned}
\right.
\]
and consider the set of the paths of all
classical orbits
\begin{equation*}
\begin{aligned}
\Gamma_a
&=\b{\p{x(t),\xi(t)}\,:\, t\in\R,\, \xi_0\in\R^n\setminus 0}
\\
&=\b{\p{\lambda\nabla a(\xi),\xi}\,:\,
 \xi\in\R^n\setminus 0,\,\lambda\in\R}.
\end{aligned}
\end{equation*}
Let pseudo-differential operator $\sigma(X,D)$
have symbol $\sigma(x,\xi)$ which is  smooth in $x\neq0$, $\xi\neq0$,
and which is 
positively homogeneous of order $-1/2$ with respect to $x$,
and of order $1/2$ with respect to $\xi$.
Suppose also the structure condition
\begin{equation}\label{EQ:gamma-a}
 \sigma(x,\xi)=0\quad \text{if}\quad (x,\xi)\in
 \Gamma_a\quad \text{and}\quad x\neq0.
\end{equation}
Then it was shown in \cite{RS3} that the solution 
$u=e^{ita(D_x)}\varphi$ to 
$(\ref{eq-cr1})$ satisfies
\begin{equation}\label{math-ann}
\n{\sigma(X,D_x) e^{ita(D_x)}\varphi(x)}_{L^2\p{\R_t\times\R^n_x}}
\leq C
\n{\varphi}_{L^2(\R^n)}
\end{equation}
if $n\geq2$ and the Gaussian curvature of
the hypersurface
\begin{equation}\label{hyper}
\Sigma_a=\b{\xi\in\Rn\,:\, a(\xi)=1}
\end{equation}
never vanishes.
The typical example for such critical operator $\sigma(X,D_x)$ is given by
the elements of 
\begin{equation}\label{eqex1}
\Omega_1=|x|^{-1/2}
\p{\frac x{|x|}\wedge\frac{\nabla a(D_x)}{|\nabla a(D_x)|}}|D_x|^{1/2},
\end{equation}
where the outer product $p\wedge q$ of vectors
$p=(p_1,p_2,\ldots,p_n)$ and $q=(q_1,q_2,\ldots,q_n)$
is defined by $p\wedge q=(p_iq_j-p_jq_i)_{i<j}$.
Another interesting example is the element of
\begin{equation}\label{EQ:LB-dual}
\Omega_2=|x|^{-1/2}
\p{\frac{\nabla a^*(x)}{|\nabla a^*(x)|}\wedge\frac{D_x}{|D_x|}}|D_x|^{1/2},
\end{equation}
where $a^*(x)$ is the {\it dual function} of $a(\xi)$ which is positively
homogeneous of order two and is characterised by the
relation $a^*(\nabla a(\xi))=1$.
We remark that the sum of the squares of all
elements of $\Omega_2$ forms the
main factor of the homogeneous
extension of the Laplace-Beltrami operator on the {\it dual hypersurface}
$\Sigma^*_a=\b{\nabla a(\xi):\xi\in\Sigma_a}$.
The dual function $a^*(x)$ can be also determined by the 
relation $\Sigma_{a^*}=\Sigma^*_a$.
\par
The proof of this result relies on the critical case of the
limiting absorption principle which can be proved by reducing
its statement to a model situation by the canonical transform
method combined with weighted estimates for the transform
operators.
On the other hand, it can be reduced to a corresponding
smoothing estimate for the Laplace operator with {\em any}
critical operator, for example to the homogeneous extension
of the Laplace-Beltrami operator on the sphere, 
recovering, in particular, the result of \cite{Su2}. 
This result has been extended to
include small perturbations by Barcel\'o, Bennett and Ruiz \cite{BBR}.
For further details on these arguments we refer to authors'
paper \cite{RS3}. On the other hand, the set $\Gamma_a$ 
corresponds to the Hamiltonian flow of $a(D_x)$, which is known
to play a role in such problems also in a more general setting
of manifolds. There, non-trapping conditions also enter
(e.g. Doi \cite{Do1, Do2} in the case of Schr\"odinger operators
on manifolds, using Egorov theorem, or 
Burq \cite{Bu} and Burq, G\'erard and 
Tzvetkov \cite{BGT} in the case of Schr\"odinger
boundary value problems, using propagation properties of
Wigner measures), and such conditions can be also
expressed in terms of properties of the set $\Gamma_a$. In our
case this simply corresponds to the dispersiveness of $a(D_x)$.
\par
One important topic related to this argument
are the Fourier restriction and trace theorems.
Below, we frequently quote the description in \cite[Section 5]{RS3}
which the reader may consult for the details.
First note that the formal adjoint
$T^*:\mathcal{S}(\R_t\times\R^n_x)\to\mathcal{S}'(\R^n_x)$
of the solution operator
\[
T=e^{ita(D_x)}:\mathcal{S}(\R^n_x)\to\mathcal{S}'(\R_t\times\R^n_x)
\]
to equation \eqref{eq-cr1} is expressed as
\begin{equation}\label{rel:fr}
T^*\left[v(t,x)\right]=\FT^{-1}_{\xi}
  \left[
    \p{\FT_{t,x}v}\p{a(\xi),\xi}
  \right].
\end{equation}
Then, for any operator $A=A(X,D_x)$ acting on the variable $x$, 
the estimate
\begin{equation}\label{Smest}
\n{Ae^{ita(D_x)}\varphi}_{L^2\p{\R_t\times\R^n_x}}
\le
C\,
\Vert \,\varphi\,\Vert_{L^2\p{\R^n_x}}
\end{equation}
implies the estimate
\begin{equation}\label{FRes}
\n{\widehat{A^*f}_{|\Sigma_a}}_{L^2\p{\Sigma_a\,;\,d\omega/|\nabla a|}}
\le
C \n{f}_{L^2(\R^n_x)},
\end{equation}
where $d\omega$ is the standard surface element of the hypersurface $\Sigma_a$
defined by \eqref{hyper}.
Indeed, by \eqref{rel:fr} and Plancherel's theorem,
we have for $v(t,x)=g(t)f(x)$
\begin{align*}
\n{T^*A^*v}^2_{L^2(\R^n)}
&=(2\pi)^{-n}\n{\p{\FT_{t,x}A^*v}(a(\xi),\xi)}^2_{L^2(\R^n_\xi)}
\\
&=(2\pi)^{-n}\int^\infty_0
  \p{\int_{\Sigma_a}
   \left|
     \p{\FT_{t,x}A^*v}\p{\rho^2,\rho\omega}
   \right|^2
    \,\frac{2\rho^{n-1}d\omega}{|\nabla a(\omega)|}}\,d\rho
\\
&=(2\pi)^{-n}\int^\infty_0\abs{\widehat g(\rho^2)\sqrt\rho}^2
  \p{\int_{\Sigma_a}
   \left|\frac1{\sqrt{\rho}}\p{\widehat{A^*f}}\p{\rho\omega}
   \right|^2
    \,\frac{2\rho^{n-1}d\omega}{|\nabla a(\omega)|}}\,d\rho.
\end{align*}
Here we have used the change of variables $\xi\mapsto\rho\omega$
($\rho>0,\omega\in\Sigma_a$).
At the same time, by \eqref{Smest}, we have
\[
\n{T^*A^*v}^2_{L^2(\R^n)}
\leq C\n{v}^2_{L^2(\R_t\times\R^n_x)}
=C\n{g}_{L^2(\R)}^2\n{f}_{L^2(\R^n)}^2.
\]
Note that we have by Plancherel's theorem
\[
\n{g}_{L^2(\R)}^2=\frac1{2\pi}\n{\widehat g}_{L^2(\R)}^2=
\frac1{4\pi}\int^\infty_0\abs{\widehat 
g(\rho^2)\sqrt\rho}^2  \,d\rho,
\]
if $\supp\widehat g\subset[0,\infty)$.
Combining all these relations and taking arbitrary $g$,
we have estimate
\begin{equation}\label{est:restrict}
\n{\widehat{A^*f}_{|\rho\Sigma_a}}_
{L^2\p{\rho\Sigma_a\,;\,\rho^{n-1}d\omega/|\nabla a|}}
\le
C\sqrt{\rho}\,
\n{f}_{L^2(\R^n_x)},
\end{equation}
where $\rho>0$, $\rho\Sigma_a=\{\rho\omega:\,\omega\in\Sigma_a\}$.
Taking $\rho=1$, we have estimate \eqref{FRes}.
We remark that, conversely, estimate \eqref{est:restrict} implies
estimate \eqref{Smest}.
\par
In this paper we
have already obtained the examples of operators $A$ which satisfy
smoothing estimate \eqref{Smest}, hence the Fourier restriction
estimate \eqref{FRes}.
For example, by Theorems \ref{M:H1} and \ref{M:H2},
we can take
\begin{equation}\label{specialop}
\begin{aligned}
&A_1=\jp{x}^{-s}|D_x|^{1/2}\qquad (s>1/2), \\
&A_2=|x|^{\alpha-1}|D_x|^{\alpha}\qquad (1-n/2<\alpha<1/2).
\end{aligned}
\end{equation}
We can also take $A=\sigma(X,D_x)$ which appeared
in estimate \eqref{math-ann}, especially the elements of the
operators $\Omega_1$ or $\Omega_2$ defined by \eqref{eqex1} or
\eqref{EQ:LB-dual},
but in this case we also
need the non-degenerate Gaussian
curvature condition on the hypersurface
$\Sigma_a$ defined by \eqref{hyper},
which is equivalent to $\det \nabla^2a(\xi)\not=0$ ($\xi\neq0$)
(see Miyachi \cite{Mi}, for example).
Their formal adjoints are given by
\begin{align*}
&A_1^*=|D_x|^{1/2}\jp{x}^{-s}\qquad (s>1/2),
\quad
A_2^*=|D_x|^{1-s}|x|^{-s}\qquad (1/2<s<n/2),
\\
& \Omega_1^*=
|D_x|^{1/2}
\p{\frac{\nabla a(D_x)}{|\nabla a(D_x)|}\wedge\frac x{|x|}}|x|^{-1/2},
\quad
\Omega_2^*=|D_x|^{1/2}
\p{\frac{D_x}{|D_x|}\wedge\frac{\nabla a^*(x)}{|\nabla a^*(x)|}}|x|^{-1/2}.
\end{align*}
Note that we have
$|\nabla a(\xi)|\geq C>0$ on $\Sigma_a$
since $\nabla a(\xi)\neq0$ ($\xi\neq0$) in our case.
From the construction, we have the same property
for $a^*$, as well.
We also note that
$\n{f}_{L^2_{s}(\R^n)}=\|\widehat f\|_{H^s(\Rn)}$ and
$\n{f}_{\Dot{L}^2_{s}(\R^n)}=\|\widehat f\|_{\Dot{H}^s(\Rn)}$,
where $H^s(\Rn)$ and $\Dot{H}^s(\Rn)$ are (homogeneous) Sobolev spaces
with the norms 
$\n{g}_{H^s(\Rn)}=\n{\jp{D_x}^sg}_{L^2(\Rn)}$ and
$\n{g}_{\Dot{H}^s(\Rn)}=\n{|D_x|^sg}_{L^2(\Rn)}$ respectively.
Then we can conclude the following trace results:
\medskip
\begin{thm}\label{Th:Fres}
Let
$a(\xi)\in C^\infty\p{\R^n\setminus0}$ be real-valued and satisfy
$a(\xi)>0$ and $a(\lambda \xi)=\lambda^2 a(\xi)$ for $\lambda>0$
and $\xi\neq0$. Let $\Sigma_a=\b{\xi\in\Rn: a(\xi)=1}$.
Suppose $s>1/2$. 
Then we have
\begin{equation}\label{Fres:1}
\n{f_{\,\,|\Sigma_a}}_{L^2\p{\Sigma_a\,;\,d\omega}}
\le
C \n{f}_{H^{s}(\R^n)}.
\end{equation}
Moreover, for $1/2<s<n/2$, we have
\begin{equation}\label{Fres:1.5}
\n{f_{\,\,|\Sigma_a}}_{L^2\p{\Sigma_a\,;\,d\omega}}
\le
C \n{f}_{\Dot{H}^{s}(\R^n)}.
\end{equation}
If we in addition assume that the Gaussian curvature of
$\Sigma_a$ is non-vanishing, then we have also
\begin{equation}\label{Fres:2}
\n{\p{\frac{\nabla a(x)}{|\nabla a(x)|}\wedge\frac {D_x}{|D_x|}}
f_{\,\,|\Sigma_a}}_{L^2\p{\Sigma_a\,;\,d\omega}}
\le
C \n{f}_{\Dot{H}^{1/2}(\R^n)}
\end{equation}
and
\begin{equation}\label{Fres:3}
\n{\p{\frac x{|x|}\wedge\frac {\nabla a^*(D_x)}{|\nabla a^*(D_x)|}}
f_{\,\,|\Sigma_a}}_{L^2\p{\Sigma_a\,;\,d\omega}}
\le
C \n{f}_{\Dot{H}^{1/2}(\R^n)},
\end{equation}
where $a^*(x)$ is the dual function of $a(\xi)$.
\end{thm}
\medskip
The third and fourth
estimates \eqref{Fres:2} and \eqref{Fres:3} in Theorem \ref{Th:Fres} say that
we can attain the critical order $s=1/2$ in the first and second estimates
\eqref{Fres:1} and \eqref{Fres:1.5} under a structure condition. 
In fact, we get a complete range of operators for the critical
smoothing if we use the restriction estimate
\eqref{FRes} with $A=\sigma(X,D_x)$ satisfying
\eqref{EQ:gamma-a}. Estimates \eqref{Fres:2} and \eqref{Fres:3} are
the interesting special cases of them. 
\par
We note finally, that the results on the global smoothing imply the
growth rates of the restriction norms. For example,
smoothing with operator $A_1$ in \eqref{specialop} implies the
uniform trace estimate
$$
\n{f_{|\rho\Sigma_a}}_{L^2(\rho\Sigma_a,\rho^{n-1}
d\omega)} \leq C \n{f}_{H^s(\Rn)}
\qquad (s>1/2)
$$
on account of \eqref{est:restrict}.
If we use $A_2$ in \eqref{specialop} instead, 
we get
$$
\n{f_{|\rho\Sigma_a}}_{L^2(\rho\Sigma_a,\rho^{n-1}
d\omega)} \leq C \rho^{s-1/2} \n{f}_{\Dot{H}^s(\Rn)}
\qquad (n/2>s>1/2).
$$ 
However, in the critical cases, in addition to
\eqref{Fres:2} and \eqref{Fres:3}, we obtain
\begin{align*}
&\n{\p{\frac{\nabla a(x)}{|\nabla a(x)|}\wedge\frac {D_x}{|D_x|}}
f_{|\rho\Sigma_a}}_{L^2(\rho\Sigma_a,\rho^{n-1}
d\omega)}  \leq
C \n{f}_{\Dot{H}^{1/2}(\Rn)},
\\
&\n{\p{\frac x{|x|}\wedge\frac {\nabla a^*(D_x)}{|\nabla a^*(D_x)|}}
f_{|\rho\Sigma_a}}_{L^2(\rho\Sigma_a,\rho^{n-1}
d\omega)}  \leq
C \n{f}_{\Dot{H}^{1/2}(\Rn)}.
\end{align*}



\end{document}